\documentclass[11pt]{amsart}
\usepackage{amsmath,amsfonts}
\usepackage{amssymb}
\usepackage{graphicx,color}
\usepackage{epic,eepic}
\usepackage{algpseudocode}
\usepackage{longtable}
\algblockdefx[Foreach]{Foreach}{EndForeach}[1]{\textbf{for each} #1 \textbf{do}}{\textbf{end for}}

\renewcommand{\labelenumi}{(\theenumi)}
\newcommand{\C}{\mathbb{C}}
\newcommand{\CP}{\mathbb{CP}}
\newcommand{\R}{\mathbb{R}}
\newcommand{\Z}{\mathbb{Z}}
\newcommand{\N}{\mathbb{N}}
\newcommand{\htop}{h_{\mathrm{top}}}

\renewcommand{\Im}{\mathrm{Im}}
\renewcommand{\Re}{\mathrm{Re}}

\numberwithin{equation}{section}
\newtheorem{dfn}{Definition}[section]
\newtheorem{thm}[dfn]{Theorem}
\newtheorem{prp}[dfn]{Proposition}
\newtheorem{rmk}[dfn]{Remark}
\newtheorem{lmm}[dfn]{Lemma}
\newtheorem{cor}[dfn]{Corollary}

\begin{document}

\title{On parameter loci of the H\'enon family}


\author{Zin Arai}
\address{Chubu University Academy of Emerging Sciences, Kasugai, Aichi 487--8501, Japan.}
\email{zin@isc.chubu.ac.jp}

\author{Yutaka Ishii}
\address{Department of Mathematics, Kyushu University, Motooka, Fukuoka 819--0395, Japan.}
\email{yutaka@math.kyushu-u.ac.jp}

\date{March 7, 2018}

\begin{abstract}
The purpose of the current article is to investigate the dynamics of the H\'enon family $f_{a, b} : (x, y) \mapsto (x^2-a-by, x)$, where $(a, b)\in \R\times\R^{\times}$ is the parameter~\cite{H}. We are interested in certain geometric and topological structures of two loci of parameters $(a, b)\in\R\times\R^{\times}$ for which $f_{a, b}$ share common dynamical properties; one is the \textit{hyperbolic horseshoe locus} where the restriction of $f_{a, b}$ to its non-wandering set is hyperbolic and topologically conjugate to the full shift with two symbols, and the other is the \textit{maximal entropy locus} where the topological entropy of $f_{a, b}$ attains the maximal value $\log 2$ among all H\'enon maps. 

The main result of this paper states that these two loci are characterized by the graph of a real analytic function from the $b$-axis to the $a$-axis of the parameter space $\R\times\R^{\times}$, which extends in full generality the previous result of Bedford and Smillie~\cite{BS2} for $|b|<0.06$. As consequences of this result, we show that (i) the two loci are both connected and simply connected in $\{b>0\}$ and in $\{b<0\}$, (ii) the closure of the hyperbolic horseshoe locus coincides with the maximal entropy locus, (iii) the boundaries of both loci are identical and piecewise analytic with two analytic pieces. Among others, the consequence (i) indicates a weak form of monotonicity of the topological entropy as a function of the parameter $(a, b)\mapsto \htop(f_{a, b})$ at its maximal value. 

The proof consists of both theoretical and computational parts. In the theoretical part we extend both the dynamical and the parameter spaces over $\C$, investigate their complex dynamical and complex analytic properties, and reduce them to obtain the conclusion over $\R$ as in~\cite{BS2}. One of our new ingredients is to employ a flexible family of ``boxes'' in $\C^2$ that are intrinsically two-dimensional and works for all values of $b$. In the computational part we use interval arithmetic together with some numerical algorithms such as set-oriented computations and the interval Krawczyk method to verify certain numerical criteria which imply analytic, combinatorial and dynamical consequences.
\end{abstract}

\maketitle

\newpage

\tableofcontents

\newpage

\section{Introduction and Statements of Results}\label{section1}

\subsection{Preliminaries}\label{subsection1.1}

In his celebrated paper~\cite{H} published in 1976, the French mathematician/astronomer Michel H\'enon introduced a two-parameter family of polynomial automorphisms of the plane, now called the \textit{H\'enon family}:
\[f_{a, b}\, : \, (x, y) \longmapsto (x^2-a-by, x),\]
where $(a, b)\in \R\times\R^{\times}$ is the parameter with $b\ne 0$. He obtained this family of maps as an algebraic reduction of a Poincar\'e section of the Lorenz system~\cite{L} in which chaos in the sense of sensitive dependence on initial conditions was first discovered. Among other things in the paper, H\'enon numerically demonstrated the existence of a so-called strange attractor for the parameter $(a, b)=(1.4, -0.3)$. Since then, the H\'enon family has been regarded as one of the most fundamental classes of nonlinear systems and much work has been done for this family. However, the understanding of the dynamics is still far from being complete to this day.

In this article we are interested in certain geometric and topological structures of two loci of parameters $(a, b)\in \R\times\R^{\times}$ for which $f_{a, b}$ share common dynamical properties. To motivate them, let us recall some basic terminologies in the theory of dynamical systems.

First, let $X$ be a compact metrizable space and $f : X\to X$ be a continuous map. Take a metric $d$ on $X$. For $n\in\N$ and $\varepsilon>0$, a subset $E\subset X$ is called \textit{$(n, \varepsilon)$-separated} if for any distinct $x, y\in E$, there exists $0\leq k<n$ so that $d(f^k(x), f^k(y))\geq \varepsilon$. The \textit{topological entropy} of $f$ is given by
\[\htop(f)\equiv \sup_{\varepsilon>0}\limsup_{n\to \infty}\frac{1}{n}\log \sup\bigl\{\mathrm{card}(E) : E\mbox{ is $(n, \varepsilon)$-separated}\bigr\},\]
where $\mathrm{card}(E)$ denotes the cardinality of $E$. It is known that $\htop(f)$ is a topological conjugacy invariant and, in particular, it does not depend on the choice of a metric. Moreover, when $f$ is a homeomorphism, we have $\htop(f)=\htop(f^{-1})$. A point $x\in X$ is \textit{non-wandering} if for any neighborhood $U$ of $x$ there is $N$ so that $f^N(U)\cap U\ne\emptyset$ holds. Let $\mathrm{\Omega}(f)$ be the set of non-wandering points of $f$, called the \textit{non-wandering set} of $f$. Then, it is known that $\htop(f)=\htop(f|_{\mathrm{\Omega}(f)})$, i.e. the topological entropy is concentrated in $\mathrm{\Omega}(f)$.

Next, let $\{0, 1\}^{\Z}$ be the space of bi-infinite symbol sequences with two symbols $0$ and $1$ equipped with the metric:
\[d(\underline{\varepsilon}, \underline{\varepsilon}')\equiv \sum_{n\in\Z}\frac{|\varepsilon_n-\varepsilon_n'|}{2^{|n|}}\]
for $\underline{\varepsilon}=(\varepsilon_n)_{n\in\Z}, \underline{\varepsilon}'=(\varepsilon'_n)_{n\in\Z}\in\{0, 1\}^{\Z}$. The \textit{shift map} is defined by
\[\sigma : \{0, 1\}^{\Z}\ni \cdots\varepsilon_{-1}\cdot\varepsilon_0\varepsilon_1\cdots\longmapsto \cdots\varepsilon_{-1}\varepsilon_0\cdot\varepsilon_1\cdots \in\{0, 1\}^{\Z},\]
where $\cdot$ is placed at the left of the $0$-th digit. It is easy to see that $(\{0, 1\}^{\Z}, d)$ is a compact metric space and $\sigma$ is a continuous map. One can moreover compute that $\htop(\sigma)= \log 2$. 

Finally, let $M$ be a smooth manifold and $f : M\to M$ be a smooth diffeomorphism. An invariant set $\mathrm{\Lambda}$ for $f$ is called \textit{hyperbolic} if there exist a continuous splitting $T_{\mathrm{\Lambda}}M=E^s \oplus E^u$ of the tangent bundle over $\mathrm{\Lambda}$ and two constants $C>0$ and $0<\lambda<1$ so that the following conditions are satisfied
\begin{enumerate}
\renewcommand{\labelenumi}{(\roman{enumi})}
\item $Df_p(E^s_p)=E^s_{f(p)}$ and $Df_p(E^u_p)=E^u_{f(p)}$ for all $p\in\mathrm{\Lambda}$,
\item $\| Df^n(v)\|\leq C\lambda^n \|v\|$ for all $v\in E^s$ and $n>0$,
\item $\| Df^{-n}(v)\|\leq C\lambda^n \|v\|$ for all $v\in E^u$ and $n>0$,
\end{enumerate}
with respect to some Riemannian metric $\|\cdot\|$ on $M$. Let us say that $f$ is a \textit{hyperbolic horseshoe on $M$} if the non-wandering set $\mathrm{\Omega}(f)$ is a hyperbolic set and the restriction $f|_{\mathrm{\Omega}(f)} : \mathrm{\Omega}(f)\to\mathrm{\Omega}(f)$ is topologically conjugate to the shift map $\sigma : \{0, 1\}^{\Z}\to\{0, 1\}^{\Z}$.

\begin{figure}
  \includegraphics[width=15cm]{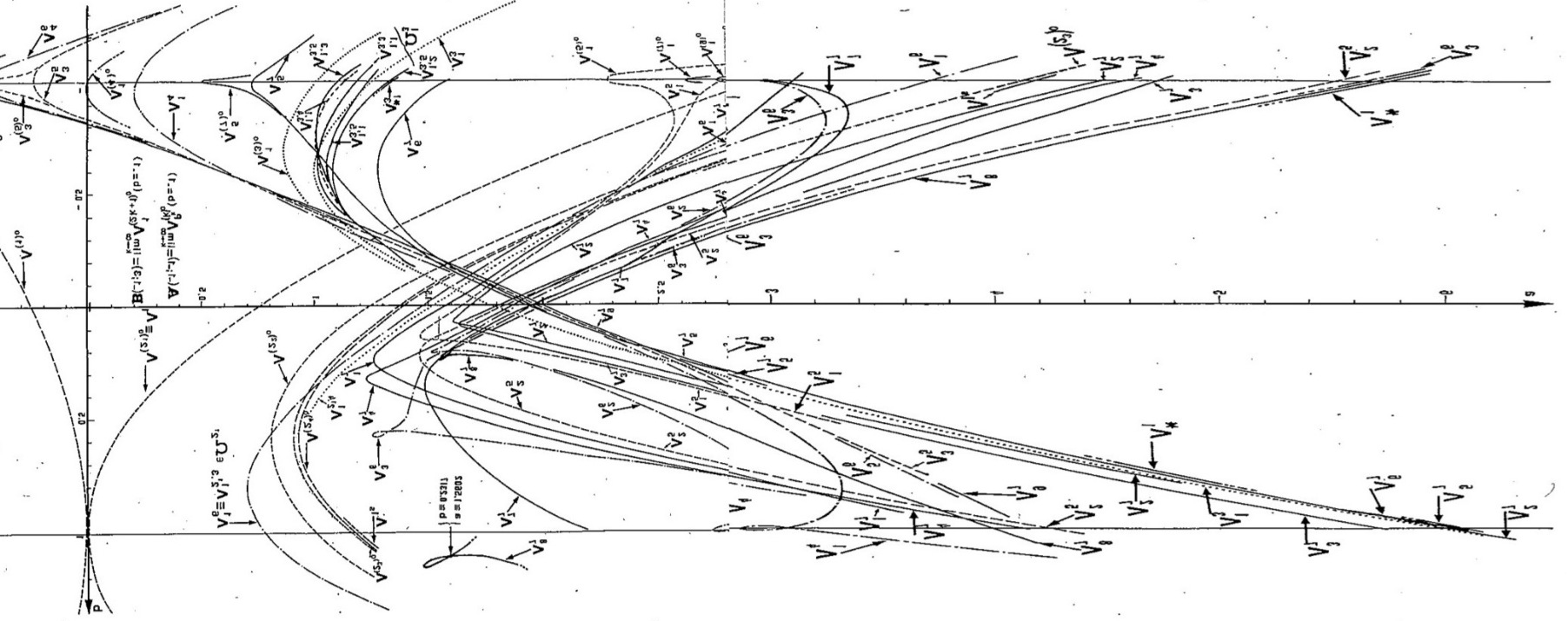}
  \caption{Bifurcation curves of the H\'enon family~\cite{EM}.}
  \label{FIG:mira}
\end{figure}

\subsection{Main results}\label{subsection1.2}

The dynamics of a H\'enon map $f_{a, b}$ depends on the choice of $(a, b)$. Let us glimpse how the dynamics of $f_{a, b}$ changes when $(a, b)$ varies. Suppose first that $b\ne 0$ is fixed and $a$ is small enough. Then, an easy computation shows that $(a, b)$ does not have periodic points of period at most two. By Brouwer's translation theorem we know that the dynamics of $f^2_{a, b}$ is topologically conjugate to a translation. It follows from this that the non-wandering set $\mathrm{\Omega}(f^2_{a, b})$ is empty, and hence the topological entropy of $f_{a, b}$ is zero (here we compactify $\R^2$ by adding a point at infinity $\infty$ and set $f_{a, b}(\infty)=\infty$). Suppose next that $b\ne 0$ is fixed and $a$ is large enough. Then, it was shown in~\cite{DN} that $f_{a, b}$ becomes a hyperbolic horseshoe on $\R^2$. Since the topological entropy of $f_{a, b}$ satisfies $0\leq \htop(f_{a, b})\leq\log 2$ for any $(a, b)\in\R\times\R^{\times}$ (see~\cite{FM}), this yields that \textit{$f_{a, b}$ attains the maximal entropy on $\R^2$} among the H\'enon maps, i.e. $\htop(f_{a, b})=\log 2$.

The notion of a horseshoe has been first introduced by Stephen Smale~\cite{S} and is regarded as one of the simplest models of a chaotic dynamical system. For several decades one of the central problems in the study of dynamical systems is to understand how a horseshoe is created through a bifurcation process. The discussion in the previous paragraph tells that the H\'enon family contains a transition from a translation to a horseshoe, i.e. a route from trivial dynamics to chaos. In this paper we focus on the last bifurcation problem among several aspects of the creation of horseshoes, which asks when and how the creation of horseshoes is completed. Equivalently, the problem is to investigate the topological and geometric structure of the the locus in the parameter space where the maps exhibit horseshoes, and to determine how the horseshoe structure is destroyed for maps in the locus boundary.

We are thus led to introduce the \textit{hyperbolic horseshoe locus}:
\[\mathcal{H}_{\R}\equiv \bigl\{(a, b) \in\R\times\R^{\times} : \mbox{$f_{a, b}$ is a hyperbolic horseshoe on $\R^2$} \bigr\}\]
as well as the \textit{maximal entropy locus}:
\[\mathcal{M}_{\R}\equiv \bigl\{(a, b) \in\R\times\R^{\times} : \mbox{$f_{a, b}$ attains the maximal entropy on $\R^2$} \bigr\}.\]
Note that $\mathcal{H}_{\R}$ is an open subset of $\R\times\R^{\times}$ and, since the topological entropy $\htop(f_{a, b})$ is a continuous function of $(a, b)$ by combining the results of~\cite{K} and~\cite{N, Y} (see page 110 of~\cite{M}), $\mathcal{M}_{\R}$ is a closed subset in $\R\times\R^{\times}$ and hence $\overline{\mathcal{H}_{\R}}\subset \mathcal{M}_{\R}$ holds. In~\cite{BS2} Bedford and Smillie have shown that these two parameter loci are characterized by a real analytic curve for $|b|<0.06$ (see also~\cite{CLR} on a weaker result for a wider class of families called the H\'enon-like families). The goal of this paper is to extend this result in full generality. Namely,

\medskip

\noindent
\textbf{Main Theorem.}
\textit{There exists a real analytic function $a_{\mathrm{tgc}} : \R^{\times} \to \R$ from the $b$-axis to the $a$-axis of the parameter space $\R\times\R^{\times}$ for the H\'enon family $f_{a, b}$ with $\lim_{b\to 0}a_{\mathrm{tgc}}(b)=2$ so that
\begin{enumerate}
\renewcommand{\labelenumi}{(\roman{enumi})}
\item $(a, b)\in \mathcal{H}_{\R}$ iff $a>a_{\mathrm{tgc}}(b)$, 
\item $(a, b)\in \mathcal{M}_{\R}$ iff $a\geq a_{\mathrm{tgc}}(b)$.
\end{enumerate}
Moreover, the map $f_{a, b}$ with $a=a_{\mathrm{tgc}}(b)$ has exactly one orbit of homoclinic (resp. heteroclinic) tangencies of stable and unstable manifolds of suitable fixed points when $b>0$ (resp. $b<0$).
}

\medskip

The statements described in the Main Theorem justify what were numerically computed at the beginning of 1980's by El Hamouly and Mira, Tresser, Ushiki and others. Figure~\ref{FIG:mira} is obtained by joining two figures in the numerical work of El Hamouly and Mira~\cite{EM} and turning it upside down. There, the graph of the function $a_{\mathrm{tgc}}$ is implicitly figured out by the right-most wedge-shaped curve. 

The Main Theorem in particular yields that the maps in $\mathcal{M}_{\R}$ lose their hyperbolicity exactly at the boundary of $\mathcal{M}_{\R}$ and the hyperbolicity persists over the interior of $\mathcal{M}_{\R}$. The proof of this persistence of hyperbolicity heavily depends on the deep dichotomy result for H\'enon maps with maximal entropy on $\R^2$ by Bedford and Smillie~\cite{BS1}. The existence of an orbit of homoclinic/heteroclinic tangencies (modulo the uniqueness) for the map with $a=a_{\mathrm{tgc}}(b)$ in the Main Theorem has been already obtained in~\cite{BS1}, and we give an alternative proof of this fact together its uniqueness.

A crucial step in~\cite{BS2} was to construct a family of ``boxes'' in $\C^2$ for $|b|<0.06$. This kind of boxes were first used in~\cite{HO} and later in~\cite{BS2, I1, I2, I3, ISm}. In the current paper, we introduce a new family of flexible boxes in $\C^2$ which is intrinsically two-dimensional and works for all values of $b$. This enables us to understand the global topology of the two loci. To state it, let us put\footnote{For a claim $X(\pm)$ containing the symbol $\pm$, the statement ``\textit{$X(\pm)$ holds}'' means ``\textit{both $X(+)$ and $X(-)$ hold}''. This convention applies when $X(\pm)$ is a definition as well, e.g. $\mathcal{H}_{\R}^{\pm}\equiv \mathcal{H}_{\R}\cap\{\pm b>0\}$ means $\mathcal{H}_{\R}^+\equiv \mathcal{H}_{\R}\cap\{b>0\}$ and $\mathcal{H}_{\R}^-\equiv \mathcal{H}_{\R}\cap\{b<0\}$.} 
\[\mathcal{H}_{\R}^{\pm}\equiv \mathcal{H}_{\R}\cap\{\pm b>0\} \quad \mbox{and} \quad \mathcal{M}_{\R}^{\pm}\equiv \mathcal{M}_{\R}\cap\{\pm b>0\}.\] 
Below, we take the closure and the boundary of the loci $\mathcal{H}_{\R}^{\pm}$ and $\mathcal{M}_{\R}^{\pm}$ in $\{\pm b>0\}$.

\medskip

\noindent
\textbf{Main Corollary.}
\textit{Both loci $\mathcal{H}_{\R}^{\pm}$ and $\mathcal{M}_{\R}^{\pm}$ are connected and simply connected in $\{\pm b>0\}$. Moreover, we have $\overline{\mathcal{H}_{\R}^{\pm}}=\mathcal{M}_{\R}^{\pm}$ and $\partial \mathcal{H}_{\R}^{\pm}=\partial \mathcal{M}_{\R}^{\pm}$.}

\medskip

As far as we know, this is the first result which determines global topological properties of parameter loci for the real H\'enon family. Moreover, this result can be regarded as a first step towards the understanding of an ``ordered structure'' in the H\'enon parameter space. Recall that in~\cite{MT} the monotonicity of the topological entropy for the cubic family (which has two parameters) is formulated as the connectivity of isentropes. In this sense, the Main Corollary indicates a weak form of monotonicity of the function $(a, b)\mapsto \htop(f_{a, b})$ at its maximal value. 

It is interesting to compare our results to the so-called \textit{anti-monotonicity theorem} in~\cite{KKY}. To be precise, we let $h_t : \R^2\to\R^2$ ($t\in\R$) be a one-parameter family of dissipative $C^3$-diffeomorphisms of the plane and assume that $h_{t_0}$ has a non-degenerate homoclinic tangency for certain $t=t_0$. The theorem states that there are both infinitely many orbit-creation and infinitely many orbit-annihilation parameters in any neighborhood of $t_0\in\R$. It has been shown in~\cite{BS2} that for the one-parameter family of H\'enon maps $\{f_{a, b_{\ast}}\}_{a\in\R}$ with a fixed $b_{\ast}>0$ close to zero, the homoclinic tangency of $f_{a, b_{\ast}}$ at $a=a_{\mathrm{tgc}}(b_{\ast})$ mentioned above is non-degenerate, hence the anti-monotonicity theorem applies. Of course, anti-monotonicity of some orbits does not necessarily imply anti-monotonicity of topological entropy or creation/destruction of horseshoes. Nonetheless, this theorem suggests that, a priori, $\mathcal{H}_{\R}$ and $\mathcal{M}_{\R}$ could have holes or other connected components separated from the ones described in the Main Corollary.

\begin{figure}
  \vspace*{3cm}
  \includegraphics[width=15.8cm]{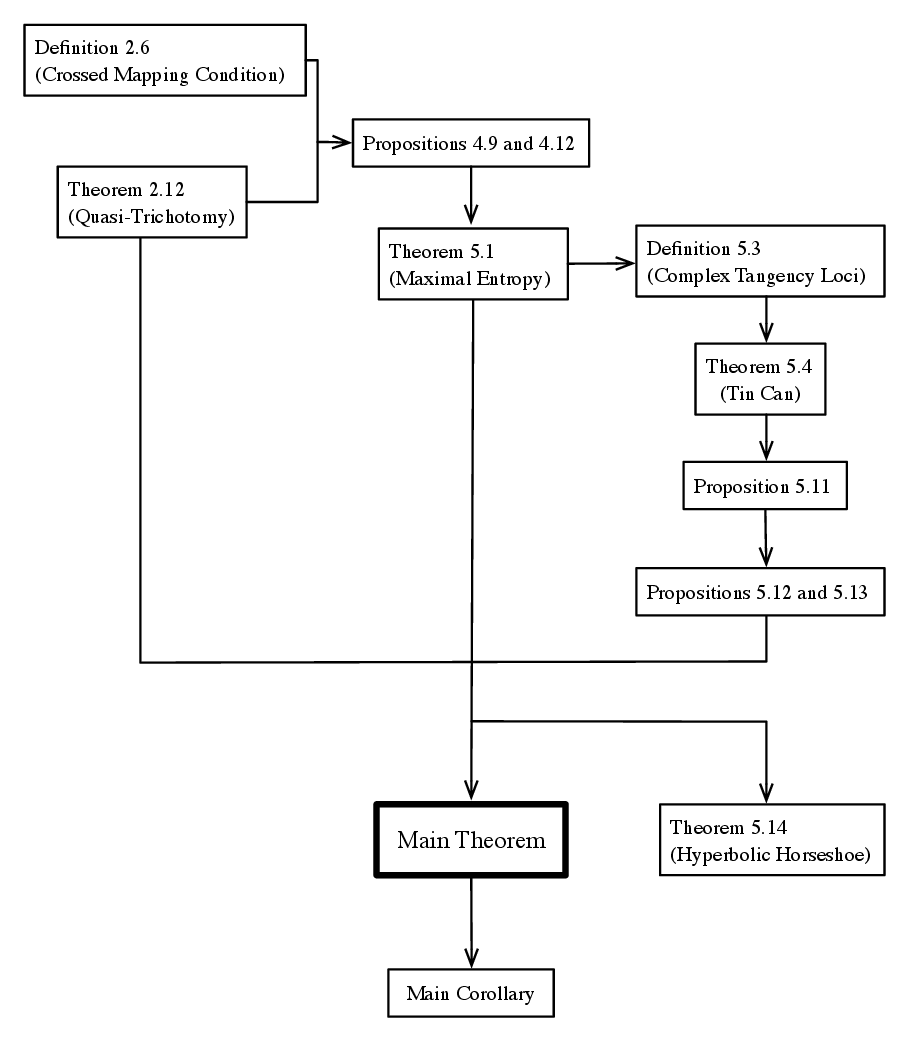}
  \caption{The flowchart of the proof of the Main Theorem.}
  \vspace*{2.5cm}
  \label{FIG:chart}
\end{figure}

\subsection{Open questions}\label{subsection1.3}

Let us discuss some open questions and remarks related to our results.

First, as clearly seen in Figure~\ref{FIG:mira}, the function $a_{\mathrm{tgc}}$ looks monotone both on $\{b>0\}$ and on $\{b<0\}$. It would be interesting to give a rigorous proof of this observation. Indeed, in a forthcoming paper~\cite{AIT} we apply the framework of this article to estimate the slope of the function $a_{\mathrm{tgc}}$ near $b=0$. As a consequence of this estimate, we obtain a variational characterization of equilibrium measures at ``temperature zero'' for real H\'enon maps at the last bifurcation parameter $(a, b)\in\mathcal{H}_{\R}^+$ with $b>0$ close to zero.

As the second question one may ask if an analogy of the Main Corollary holds for the \textit{complex} H\'enon family $f_{a, b} : \C^2\to\C^2$ with $(a, b)\in \C\times\C^{\times}$. For this family we define the locus $\mathcal{H}_{\C}$ as the set of parameters $(a, b)\in \C\times\C^{\times}$ for which the restriction of $f_{a, b}$ to $\mathrm{\Omega}(f_{a, b})$ in $\C^2$ is hyperbolic and is topologically conjugate to the shift map $\sigma : \{0, 1\}^{\Z}\to\{0, 1\}^{\Z}$. It is easy to see that $\mathcal{H}_{\C}$ is not simply connected. In fact, the two fixed points of $f_{a, b}$ are interchanged by changing the parameter along the loop $\gamma(t)=(a(t), b_0)$ where $|b_0|$ is small and $a(t)=Re^{2\pi i t}$ is a large circle with $a(0)=a_0$. In particular, the image of $\gamma$ by the monodromy representation $\rho : \pi_1(\mathcal{H}_{\C}, (a_0, b_0))\to \mathrm{Aut}(\{0, 1\}^{\Z}, \sigma)$ is non-trivial and hence $\mathcal{H}_{\C}$ is not simply connected (see also Proposition 6.1 in~\cite{BS3}). Moreover, Arai~\cite{A2} found a loop $\gamma\in\pi_1(\mathcal{H}_{\C}, (a_0, b_0))$ so that $\rho(\gamma)$ has infinite order in $\mathrm{Aut}(\{0, 1\}^{\Z}, \sigma)$. It is however an open question if $\mathcal{H}_{\C}$ is connected. On the other hand, the topological entropy of $f_{a, b}$ on $\C^2$ is always $\log 2$ and independent of the parameter~\cite{Sm}. Therefore, there is no analogous locus to $\mathcal{M}_{\R}$ in the complex setting. 

In this article we have analyzed the two parameter loci where the dynamics is ``maximal''. As the third problem we propose to investigate the opposite side of the parameter space, i.e. the \textit{zero-entropy locus} for the H\'enon family $\mathcal{Z}\equiv\{(a, b)\in\R\times\R^{\times} : \htop(f_{a, b})=0\}$. Recall that Katok~\cite{K} has shown that for a $C^{1+\alpha}$ diffeomorphism $f$ on a compact surface, its topological entropy is strictly positive if and only if $f^n$ contains a hyperbolic horseshoe for some $n\geq 1$. Therefore, the boundary of the zero-entropy locus $\partial\mathcal{Z}$ is often referred to as the ``boundary of chaos''. We conjecture that $\partial\mathcal{Z}$ is piecewise real analytic (see also page 19 of~\cite{GT}). Notice that for $b$ close to zero, this conjecture has been already solved in Theorem 2.2 of~\cite{GST} (see also Corollary 4.5 of~\cite{CLM}). 

Indeed, this conjecture is motivated by the comparison with a piecewise affine model of the H\'enon family called the \textit{Lozi family} $L_{a, b}\, : \, (x, y) \mapsto (1-a|x|+by, x)$. In~\cite{I4, ISa} it has been proved that both the hyperbolic horseshoe locus and the maximal entropy locus for the Lozi family are characterized by an \textit{algebraic} curve, similar to the Main Theorem. As a consequence, we have shown that exactly the same statement of the Main Corollary holds for the Lozi family. We also conjectured that the boundary of the zero-entropy locus for the Lozi family would be piecewise algebraic with countably many algebraic pieces (this conjecture has been also proposed by C.~Tresser) and proposed a strategy of its proof in~\cite{ISa}. Although there is a negative result on the conjugacy problem between H\'enon maps and Lozi maps~\cite{T}, we expect that it would be fruitful to compare the dynamics of these two families.

\subsection{Outline of proof}\label{subsection1.4}

The proof of our results consists of computational part and theoretical part. In the theoretical part we extend both the dynamical and the parameter spaces over $\C$, investigate their complex dynamical and complex analytic properties, and then reduce them to obtain the conclusion over $\R$ as in~\cite{BS2}. The idea of applying complex method to real dynamics in dimension two goes back to the earlier papers~\cite{BLS, HO, BS1}. In the computational part we employ interval arithmetic together with some numerical algorithms to verify numerical criteria which imply analytic, combinatorial and dynamical consequences (see Section~\ref{section6} for the idea of a computer-assisted proof). Below we discuss an outline of the proof with an emphasis on the new ingredients. Figure~\ref{FIG:chart} is a flowchart describing the implications between principal statements. In Table~\ref{TAB:notations} at the end of this section we summarize the notations in this article.

The starting point of our discussion is to classify any H\'enon map into the following three types (Theorem~\ref{THM:quasitrichotomy}); either (i) $\htop(f_{a, b})<\log 2$, (ii) $f_{a, b}$ is a hyperbolic horseshoe on $\R^2$, or (iii) $f_{a, b}$ for $(a, b)$ in a complex neighborhood of $\partial\mathcal{H}_{\R}^{\pm}=\partial\mathcal{M}_{\R}^{\pm}$ satisfies the \textit{crossed mapping condition} (see Definition~\ref{DFN:CMC}) with respect to a family of projective bidisks $\{\mathcal{B}_i^{\pm}\}_i$ (note that they are not exclusive). Thanks to this classification we can focus on the case (iii). In this case the family of projective bidisks allows us to partition the complex stable/unstable manifolds of $f_{a, b}$ into several pieces in terms of symbolic dynamics. By restricting the parameter $(a, b)$ to be real and the stable/unstable manifolds of $f_{a, b}$ to $\R^2$, certain plane topology arguments together with the crossed mapping condition implies that these pieces are properly configured in the bidisks (Propositions~\ref{PRP:pieces_positive} and~\ref{PRP:pieces_negative}). This enables us to detect which pieces are responsible for the last bifurcation for the creation of horseshoes and hence to characterize $\mathcal{M}_{\R}$ (Theorems~\ref{THM:maximalentropy}) as well as $\mathcal{H}_{\R}$ (Theorem~\ref{THM:hyperbolichorseshoe}). 

We are thus led to define the complex tangency loci $\mathcal{T}^{\pm}$ to be the complex parameters for which the corresponding complex special pieces have tangencies (Definition~\ref{DFN:tangency_loci}). Since $\mathcal{T}^{\pm}$ form complex subvarieties~\cite{BS0}, our problem is to show that they are non-singular. For this, we first verify a certain condition (Theorem~\ref{THM:tincan}) to prove that the projection from $\mathcal{T}^{\pm}$ to the $b$-axis is a proper map. The transversality of the quadratic family $p_a(x)=x^2-a$ at $a=2$ yields that its degree is one. Therefore, a version of the Weierstrass preparation theorem yields that $\mathcal{T}^{\pm}$ are complex submanifolds (Proposition~\ref{PRP:complex_loci}). This allows us to define the real analytic function $a_{\mathrm{tgc}}$ so that its graph coincides with the real part of $\mathcal{T}^{\pm}$ (Propositions~\ref{PRP:real_loci_positive} and~\ref{PRP:real_loci_negative}), which finishes the proof.

The first significant ingredient in our proof is a new construction of projective bidisks $\{\mathcal{B}_i^{\pm}\}_i$ in Theorem~\ref{THM:quasitrichotomy}. The proof of~\cite{BS2} employed a family of three bidisks in $\C^2$ called \textit{boxes} based on the Yoccoz puzzle partition for $p(z)=z^2-2$. In this paper we show that these boxes satisfy the crossed mapping condition only when $-0.5<b<0.4$ (see Appendix~\ref{appendix:comparison}). We therefore need to introduce a new  family of boxes which are intrinsically two-dimensional and are constructed based on the trellis formed by invariant manifolds in $\R^2$. This enables us to verify the necessary criteria for \textit{all} values of $b$, which is the basis of our discussion. However, there are two trade-offs of this new choice; one is that the new boxes cannot be computed algebraically in terms of the parameter and another is that the combinatorics of the transitions between the new boxes is more complicated than in~\cite{BS2}. Because of this, the numerical criteria on the behavior of boxes become impossible to verify by hand. To overcome this difficulty we use \textit{rigorous interval arithmetic}~\cite{Mo} and check several numerical criteria. 

The second significant ingredient is the introduction of numerical algorithms; \textit{set-oriented computations}~\cite{DJ} and the \textit{interval Krawczyk method}~\cite{Nm}. The former is an algorithm to generate a sequence of outer approximations of an invariant set in terms of the map and its iterates. It is used to compute the rigorous enclosure of invariant manifolds with very high accuracy, which is the key to excluding the occurrence of unnecessary tangencies. The latter is a modification of the well-known Newton's root-finding algorithm. It is used to guarantee the existence of non-real periodic orbits of $f_{a, b}$ for certain real parameter $(a, b)$. In the process of our proof, the fourth iteration of the H\'enon map is considered. This amounts to a polynomial of degree 16 and its large expansion factor increases computational error drastically. Therefore, the rigorous computation of invariant manifolds and the zeros of such polynomial with respect to projective coordinates, where its parameter varies over a small region in the parameter space, is not at all an immediate task. Without the two algorithms described above, the proof of the main results in this paper would not be accomplished.

\bigskip

\begin{longtable}{lcl} 

\caption{List of notations.}
\label{TAB:notations}
\endhead

$f_{a, b}$ & & H\'enon family (Subsection~\ref{subsection1.1}) \\

$\htop(f)$ & & topological entropy of $f$ (Subsection~\ref{subsection1.1}) \\

$\mathrm{\Omega}(f)$ & &non-wandering set of $f$ (Subsection~\ref{subsection1.1}) \\

$\sigma$ & & shift map on $\{0, 1\}^{\Z}$ (Subsection~\ref{subsection1.1}) \\

$\mathcal{H}_{\R}$ & & hyperbolic horseshoe locus (Subsection~\ref{subsection1.2}) \\

$\mathcal{M}_{\R}$ & & maximal entropy locus (Subsection~\ref{subsection1.2}) \\

$a_{\mathrm{tgc}}$ & & analytic function in the Main Theorem (Subsection~\ref{subsection1.2}) \\

$\mathcal{H}_{\R}^{\pm}$ & & intersection of $\mathcal{H}_{\R}$ with $\{\pm b>0\}$ (Subsection~\ref{subsection1.2}) \\

$\mathcal{M}_{\R}^{\pm}$ & & intersection of $\mathcal{M}_{\R}$ with $\{\pm b>0\}$ (Subsection~\ref{subsection1.2}) \\

$I^{\pm}$ & & complex neighborhood of $\{b\in\R : 0\leq \pm b\leq 1\}$ (Subsection~\ref{subsection2.1}) \\

$I_{\R}^{\pm}$ & & real part of $I^{\pm}$ (Subsection~\ref{subsection2.1}) \\

$a^{\pm}_{\mathrm{aprx}}$ & & function approximating $a_{\mathrm{tgc}}$ (Subsection~\ref{subsection2.1} and Subsection~\ref{subsection6.3}) \\

$\chi^{\pm}$ & & width of $\mathcal{F}_{\R}^{\pm}$ in the $a$-direction (Subsection~\ref{subsection2.1}) \\

$\mathcal{F}^{\pm}$ & & complex neighborhood of $\partial \mathcal{H}^{\pm}_{\R}=\partial \mathcal{M}^{\pm}_{\R}$ (Subsection~\ref{subsection2.1}) \\

$\mathcal{F}_{\R}^{\pm}$ & & real part of $\mathcal{F}^{\pm}$ (Subsection~\ref{subsection2.1}) \\

$W^{u/s}(p)$ & & real unstable/stable manifolds of $p$ (Subsection~\ref{subsection2.1}) \\

$W^{u/s}_{\mathrm{loc}}(p)$ & & local real unstable/stable manifolds of $p$ (Subsection~\ref{subsection2.1}) \\

$(\pi_u, \pi_v)$ & & projective coordinates (Subsection~\ref{subsection2.2}) \\

$D_u$, $D_v$ & & topological disks in the $u$-axis and the $v$-axis (Subsection~\ref{subsection2.2}) \\

$\times_{\mathrm{pr}}$ & & product with respect to projective coordinates (Subsection~\ref{subsection2.2}) \\

$\mathcal{B}_{\mathcal{Q}}$ & & projective box associated with $\mathcal{Q}$ (Subsection~\ref{subsection2.2}) \\

$\mathcal{Q}^{\pm}_i$ & & quadrilaterals associated with the trellis (Subsection~\ref{subsection2.2}) \\

$\mathcal{B}^{\pm}_i$ & & projective boxes associated with the trellis (Subsection~\ref{subsection2.2}) \\

$\mathfrak{T}^{\pm}$ & & set of admissible transitions (Subsection~\ref{subsection2.3}) \\

$\mathfrak{S}^{\pm}_{\mathrm{fwd}}$ & & forward admissible sequences (Subsection~\ref{subsection3.1}) \\

$\mathfrak{S}^{\pm}_{\mathrm{bwd}}$ & & backward admissible sequences (Subsection~\ref{subsection3.1}) \\

$\mathfrak{S}^{\pm}$ & & intersection of $\mathfrak{S}^{\pm}_{\mathrm{fwd}}$ and $\mathfrak{S}^{\pm}_{\mathrm{bwd}}$ (Subsection~\ref{subsection3.1}) \\

$V^{u/s}(p)$ & & complex unstable/stable manifolds at $p$ (Subsection~\ref{subsection3.2}) \\

$V^{u/s}_{\mathrm{loc}}(p)$ & & local complex unstable/stable manifolds at $p$ (Subsection~\ref{subsection3.2}) \\

$V^s_I(a, b)^{\pm}$ & & part of $V^s(p)$ with the itinerary $I$ (Subsection~\ref{subsection3.2}) \\

$V^u_J(a, b)^{\pm}$ & & part of $V^u(p)$ with the itinerary $J$ (Subsection~\ref{subsection3.2}) \\

$f_{\R}$ & & restriction of $f_{a, b}$ to $\R^2$ (beginning of Section~\ref{section4}) \\

$\mathcal{B}^{\pm}_{i, \R}$ & & real part of $\mathcal{B}^{\pm}_i$ (beginning of Section~\ref{section4}) \\

$W^s_I(a, b)^{\pm}$ & & part of $W^s(p)$ with the itinerary $I$ (Subsection~\ref{subsection4.1}) \\

$W^u_J(a, b)^{\pm}$ & & part of $W^u(p)$ with the itinerary $J$ (Subsection~\ref{subsection4.1}) \\

$W^u_{\overline{43}4124}(a, b)^-_{\mathrm{inner}}$ & & inner part of $W^u_{\overline{43}4124}(a, b)^-$ (Subsection~\ref{subsection4.1}) \\

$W^u_{\overline{43}4124}(a, b)^-_{\mathrm{outer}}$ & & outer part of $W^u_{\overline{43}4124}(a, b)^-$ (Subsection~\ref{subsection4.1}) \\

$\mathrm{upper}(\mathcal{B}^{\pm}_{i, \R})$ & & upper part of $\mathcal{B}^{\pm}_{i, \R}$ (Subsection~\ref{subsection4.2}) \\

$\mathrm{lower}(\mathcal{B}^{\pm}_{i, \R})$ & & lower part of $\mathcal{B}^{\pm}_{i, \R}$ (Subsection~\ref{subsection4.2}) \\

$\mathrm{right}(\mathcal{B}^{\pm}_{i, \R})$ & & right part of $\mathcal{B}^{\pm}_{i, \R}$ (Subsection~\ref{subsection4.2}) \\

$\mathrm{left}(\mathcal{B}^{\pm}_{i, \R})$ & & left part of $\mathcal{B}^{\pm}_{i, \R}$ (Subsection~\ref{subsection4.2}) \\

$\mathrm{outer}(\mathcal{B}^{\pm}_{i, \R})$ & & outer part of $\mathcal{B}^{\pm}_{i, \R}$ (Subsection~\ref{subsection4.2}) \\

$\mathrm{inner}(\mathcal{B}^{\pm}_{i, \R})$ & & inner part of $\mathcal{B}^{\pm}_{i, \R}$ (Subsection~\ref{subsection4.2}) \\

$(\varepsilon_u, \varepsilon_v)$ & & sign pair (Subsection~\ref{subsection4.2}) \\

$\mathcal{T}^{\pm}$ & & complex tangency loci (Subsection~\ref{subsection5.2}) \\

$\partial^v\mathcal{F}^{\pm}$ & & vertical boundaries of $\mathcal{F}^{\pm}$ (Subsection~\ref{subsection5.2}) \\

$\mathcal{V}^s_I(a, b)^{\pm}$ & & complex neighborhood of $V^s_I(a, b)^{\pm}$ (Subsection~\ref{subsection5.2}) \\

$\mathcal{V}^u_J(a, b)^{\pm}$ & & complex neighborhood of $V^u_J(a, b)^{\pm}$ (Subsection~\ref{subsection5.2}) \\

$\Psi_{a, b}$ & & uniformization of $V^u(p_3)$ (Subsection~\ref{subsection5.3}) \\

$\mathrm{\Omega}_{\mathrm{loc}}(a, b)$ & & pullback of $V^u_{\mathrm{loc}}(p_3)$ by $\Psi_{a, b}$ (Subsection 5.3) \\

$\mathrm{\Omega}_J(a, b)$ & & points in $\mathrm{\Omega}_{\mathrm{loc}}(a, b)$ with itinerary $J$ (Subsection 5.3) \\

$\varphi_a$ & & linearization of $p_a$ at $z_3$ (Subsection 5.3) \\

$p_a$ & & quadratic map $z^2-a$ (Subsection~\ref{subsection5.3}) \\

$\mathrm{\Gamma}_a$ & & parabola $\{(x, y)\in\C^2 : x=y^2-a\}$ (Subsection~\ref{subsection5.3}) \\

$\mathcal{T}^-_i$ & & irreducible components of $\mathcal{T}^-$ (Subsection~\ref{subsection5.3}) \\

$\mathcal{T}^{\pm}_{\R}$ & & real part of $\mathcal{T}^{\pm}$ (Subsection~\ref{subsection5.4}) \\

$\mathcal{T}^-_{i, \R}$ & & real part of $\mathcal{T}^-_i$ (Subsection~\ref{subsection5.4}) \\

$\kappa^{\pm}_{\R}$ & & function whose graph is $\mathcal{T}^{\pm}_{\R}$ (Subsection~\ref{subsection5.4}) \\

$\kappa^-_{i, \R}$ & & function whose graph is $\mathcal{T}^-_{i, \R}$ (Subsection~\ref{subsection5.4}) \\

$K_{g, x_0, A}$ & & the interval Krawczyk operator for $g$ (Subsection~\ref{subsection6.2}) \\

$\mathtt{F}$ & & the cubical representation of $f$ (Subsection~\ref{subsection6.4}) \\

$|\mathtt{C}|$ & & the union of cubical sets in $\mathtt{C}$ (Subsection~\ref{subsection6.4}) \\

\end{longtable}

\bigskip

\noindent
\textit{Acknowledgment.}
Y.I. thanks Eric Bedford and John Smillie for offering him the unpublished manuscript~\cite{BS0} (it has been eventually published as~\cite{BS2} except for Section 1 of~\cite{BS0}, which is now described in Appendix~\ref{appendix:regularity} of this article) when he was visiting Cornell in 2001, and for their fruitful suggestions and discussions during the conference ``New Directions in Dynamical Systems'' in 2002 at Ryukoku and Kyoto Universities as well as during their three-month stay for the International Research Project ``Complex Dynamical Systems'' at the RIMS, Kyoto University in 2003. Both of the authors thank them for allowing us to present the missing content of~\cite{BS0} in Appendix~\ref{appendix:regularity} of this article. They are also grateful to the anonymous referees for their fruitful comments which substantially improved the article. Z.A. is partially supported by JSPS KAKENHI Grant Number 23684002 and JST CREST funding program, and Y.I. is partially supported by JSPS KAKENHI Grant Numbers 25287020 and 25610020.

\newpage

\section{Quasi-Trichotomy in Parameter Space}\label{section2}

\subsection{Parameter space}\label{subsection2.1}

We first note that $f_{a, b}$ is a hyperbolic horseshoe on $\R^2$ if and only if $f_{a, b}^{-1}$ is a hyperbolic horseshoe. Similarly, $f_{a, b}$ attains the maximal entropy on $\R^2$ if and only if $f_{a, b}^{-1}$ attains the maximal entropy on $\R^2$. Since the inverse map $f^{-1}_{a, b}$ is affinely conjugate to $f_{a/b^2, 1/b}$, it is sufficient to consider the parameter region $\{(a, b)\in\R\times\R^{\times} : 0<|b|\leq 1\}$. We choose small constants $\varepsilon>0$ and $\delta>0$ \footnote{These constants are chosen so small that the results of our computer assisted proofs for the case $0 \leq \Re(b) \leq 1$ and $\Im(b) = 0$ also hold in $I^{\pm}$. See the beginning of Subsection~\ref{subsection6.4} for more details.} and define
\[I^{\pm}\equiv\bigl\{b\in\C : -\varepsilon \leq \Re(\pm b) \leq 1+\varepsilon,\, |\Im(b)| \leq \delta\bigr\}\]
and $I^{\pm}_{\R}\equiv I^{\pm}\cap \R$, where $\Re(b)$ (resp. $\Im(b)$) denotes the real (resp. imaginary) part of $b\in\C$. We note that both $I^{\pm}$ and $I^{\pm}_{\R}$ contain the degenerate case $b=0$ as well. 

\begin{figure}
  \includegraphics[width=7.5cm]{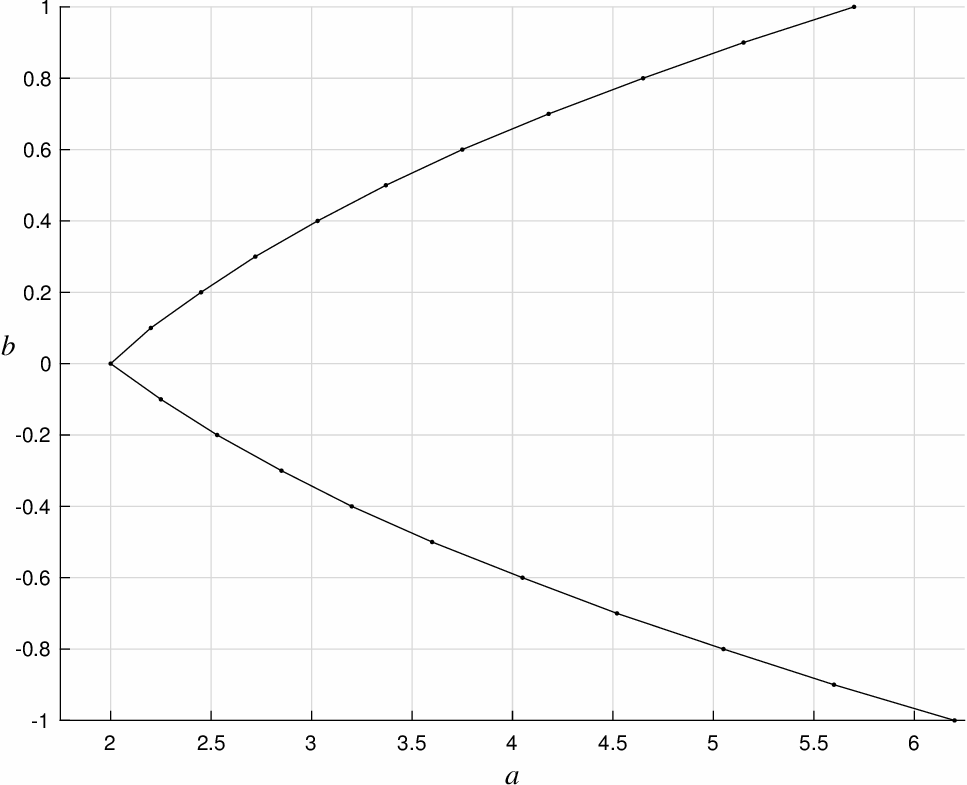}\quad
  \includegraphics[width=7.5cm]{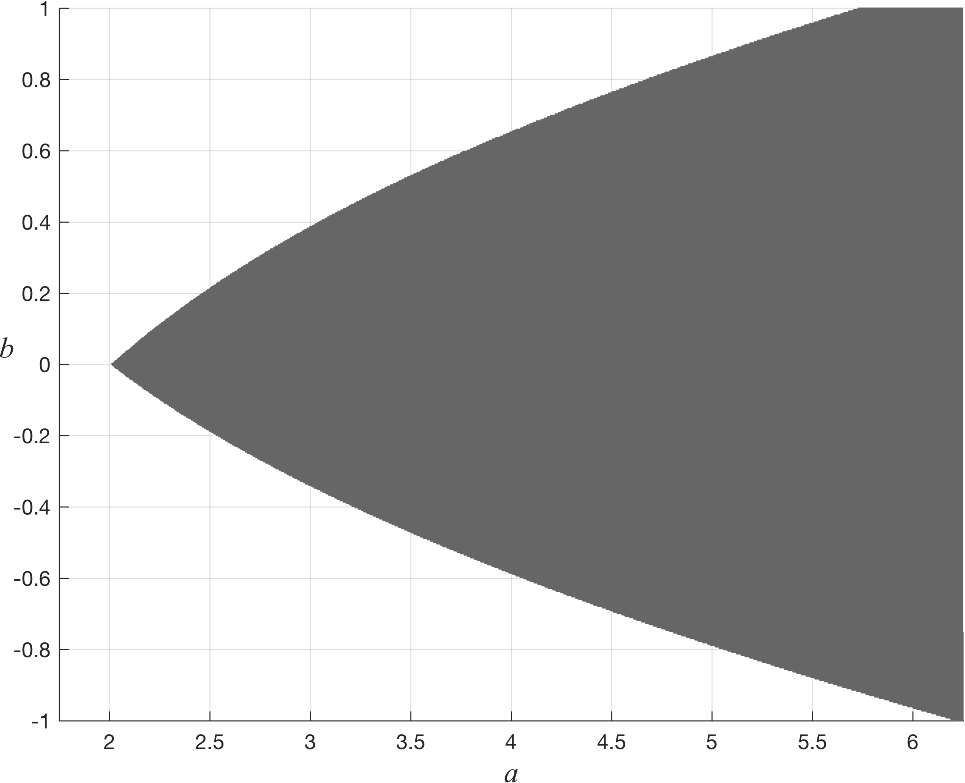}
  \caption{The graph of $a^{\pm}_{\mathrm{aprx}}$ (left) and the locus $\mathcal{H}^{\pm}_{\R}$ (right). The graph of $a^{\pm}_{\mathrm{aprx}}$ is almost identical to $\partial\mathcal{H}^{\pm}_{\R}$ and they are not distinguishable in these figures.}
\label{FIG:aprx}
\end{figure}

Let us define piecewise affine functions:
\[a^{\pm}_{\mathrm{aprx}} \, : \, I_{\R}^{\pm} \longrightarrow \R\]
to be the piecewise affine interpolations of the data given in Table~\ref{TAB:aprx} in Subsection~\ref{subsection6.3}. These are piecewise affine approximations of the function $a_{\mathrm{tgc}}$. See Figure~\ref{FIG:aprx} where we compare the graphs of $a^{\pm}_{\mathrm{aprx}}$ with $\partial \mathcal{H}^{\pm}_{\R}=\partial \mathcal{M}^{\pm}_{\R}$. The functions $a^{\pm}_{\mathrm{aprx}}$ extend to $I^{\pm}$ by letting $a^{\pm}_{\mathrm{aprx}}(b)\equiv a^{\pm}_{\mathrm{aprx}}(\Re(b))$. Put $\chi^+(b)\equiv 0.1$ for $b\in I^+$ and $\chi^-(b)\equiv 7/128 + 5 \times |\Re(b)| / 16$ for $b\in I^-$. Consider
\[\mathcal{F}^{\pm} \equiv \bigl\{(a, b) \in\C\times I^{\pm} : |a-a^{\pm}_{\mathrm{aprx}}(b)|\leq \chi^{\pm}(b) \bigr\}\]
and $\mathcal{F}^{\pm}_{\R}\equiv\mathcal{F}^{\pm}\cap \R^2$. We will see in Theorem~\ref{THM:quasitrichotomy} (Quasi-Trichotomy) that $\mathcal{F}^{\pm}$ form ``complex neighborhoods'' of $\partial \mathcal{H}^{\pm}_{\R}=\partial \mathcal{M}^{\pm}_{\R}$, and $\mathcal{F}^{\pm}_{\R}$ form ``real neighborhoods'' of $\partial \mathcal{H}^{\pm}_{\R}=\partial \mathcal{M}^{\pm}_{\R}$.

For $(a, b)\in \mathcal{F}^{\pm}_{\R}$, let $p_1\in\R^2$ (resp. $p_3\in\R^2$) be the unique fixed point in the first (resp. third) quadrant and let $p_2\in\R^2$ (resp. $p_4\in\R^2$) be the unique periodic point of period two in the second (resp. fourth) quadrant. We note that these points are well-defined in the case $b=0$ as well. The points $p_i$ then analytically continue into $\C^2$ for all $(a, b) \in \mathcal{F}^{\pm}$ which we denote again by $p_i\in\C^2$. When $(a, b)\in \mathcal{F}^{\pm}_{\R}\cap \{b\ne 0\}$, we define the real invariant manifolds $W^u(p_i)$ and $W^s(p_i)$ of $f_{a, b}|_{\R^2} : \R^2\to\R^2$ in the usual sense. When $(a, 0)\in \mathcal{F}^{\pm}_{\R}\cap \{b=0\}$, we set $W^u(p_i)\equiv\{(x, y)\in \R^2 : x=y^2-a\}$ and $W^s_{\mathrm{loc}}(p_i)\equiv\{(x, y)\in \R^2 : x=x_i\}$ where $p_i=(x_i, y_i)$.

\subsection{Projective boxes}\label{subsection2.2}

In this section we introduce the notion of projective boxes in $\C^2$. It is a generalization of coordinate bidisks, but more flexible and more useful for our purposes.

\begin{figure}
  \includegraphics[width=9cm]{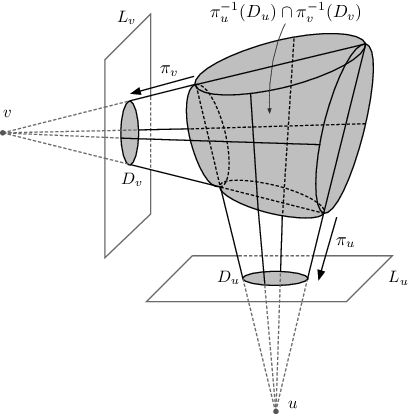}
  \caption{Projective coordinates and a projective box.}
  \label{FIG:projective_coordinates}
\end{figure}

Let us take $u\in\CP^2$ and let $L_u$ be a complex projective line in $\CP^2$ so that $u\notin L_u$. Let $L'_u$ be the unique complex line through $u$ parallel to $L_u$. Define the projection $\pi_u : \CP^2 \setminus L'_u \to L_u$ with respect to the focus $u\in \CP^2$, i.e. for $z\in \CP^2 \setminus L'_u$ we let $L$ be the unique complex line containing both $u$ and $z$, then $\pi_u(z)$ is defined as the unique point $L\cap L_u$. We call $u$ the \textit{focus} of the projection $\pi_u$ (see Figure~\ref{FIG:projective_coordinates}). 

Let $u$ and $v$ be two foci and let $L_u$ and $L_v$ be two complex lines in general position in $\CP^2$ such that $u\notin L_u$ and $v\notin L_v$. We call the pair of projections $(\pi_u, \pi_v)$ the \textit{projective coordinates} with respect to $u$, $v$, $L_u$ and $L_v$. Note that the Euclidean coordinates coincide with the projective coordinates in $\C^2$ corresponding to $u=[0 : 1 : 0]$, $v=[1 : 0 : 0]$, $L_u=\{y=0\}$ and $L_v=\{x=0\}$ under the standard identification $\CP^2 \cong\C^2\sqcup \CP^1$ by the map:
\[\CP^2 \ni [x: y: z]\longmapsto 
\begin{cases}
(x/z, y/z) \in \C^2 & \mbox{ if \, } z\ne 0, \\
[x : y] \in \CP^1 & \mbox{ if \, } z=0.
\end{cases}\]

In practice, it is sufficient to consider only the case where the foci $u$ and $v$ belong to $\C^2$ and we may assume that the complex projective lines $L_u$ and $L_v$ belong to $\C^2$ and are isomorphic to $\C$. Take two bounded topological disks $D_u\subset L_u$ and $D_v\subset L_v$ so that the following condition holds: $\pi_u^{-1}(x)\cap \pi_v^{-1}(D_v)$ is a bounded topological disk for any $x\in D_u$ and $\pi_u^{-1}(D_u)\cap \pi_v^{-1}(y)$ is a bounded topological disk for any $y\in D_v$. 

\begin{prp}
Under the assumption above, $\pi^{-1}_u(D_u)\cap \pi^{-1}_v(D_v)$ is biholomorphic to a coordinate bidisk in $\C^2$ (see Figure~\ref{FIG:projective_coordinates} again).
\label{PRP:bidisk}
\end{prp}

For a proof, see Proposition 4.6 in~\cite{I1}.

\begin{dfn}
We call $\pi_u^{-1}(D_u)\cap \pi_v^{-1}(D_v)$ a \textit{projective box} and write $D_u\times_{\mathrm{pr}} D_v$.
\label{DFN:projective_box}
\end{dfn}

\begin{figure}
  \includegraphics[height=7.7cm]{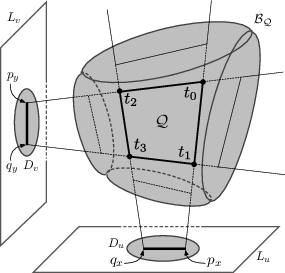}
  \caption{Projective box associated with $\mathcal{Q}$.}
  \label{FIG:projective_box}
\end{figure}

Given a quadrilateral $\mathcal{Q}$ in $\R^2$ and some additional data (such as the disks $D_u$ and $D_v$ which we shall explain shortly), we can construct a projective box as follows. Let $t_0$, $t_1$, $t_2$ and $t_3$ be the vertices of $\mathcal{Q}$ (named as in Figure~\ref{FIG:projective_box}) and assume that the segments $\overline{t_0t_1}$ and $\overline{t_2t_3}$ are close to vertical and $\overline{t_0t_2}$ and $\overline{t_1t_3}$ are close to horizontal. Let $u$ be the focus obtained as the unique intersection point of the lines containing $\overline{t_0t_1}$ and $\overline{t_2t_3}$ respectively, and let $v$ be the unique focus obtained as the unique intersection point of the lines containing $\overline{t_0t_2}$ and $\overline{t_1t_3}$ respectively. Let $L_u\equiv\{y=0\}$ be the $x$-axis of $\C^2$ and $L_v\equiv\{x=0\}$ be the $y$-axis of $\C^2$. 

\begin{dfn}
We call $(\pi_u, \pi_v)$ the \textit{projective coordinates associated with a quadrilateral $\mathcal{Q}$}. 
\end{dfn}

Let $p_x\in \R$ (resp. $q_x\in \R$) be the $x$-coordinate of the intersection of the real line containing $\overline{t_0t_1}$ (resp. $\overline{t_2t_3}$) and the $x$-axis, and $p_y\in \R$ (resp. $q_y\in \R$) be the $y$-coordinate of the intersection of the real line containing $\overline{t_0t_2}$ (resp. $\overline{t_1t_3}$) and the $y$-axis. We may assume $p_x>q_x$ and $p_y>q_y$. Then, $\pi_u(\mathcal{Q})=[q_x, p_x]$ and $\pi_v(\mathcal{Q})=[q_y, p_y]$ form intervals in $L_u$ and $L_v$ respectively. Let us choose a topological disk $D_u$ in $L_u\cong\C$ containing the interval $[q_x, p_x]\subset L_u$ and a topological disk $D_v$ in $L_v\cong\C$ containing the interval $[q_y, p_y]\subset L_v$. 

\begin{dfn}
We write $\mathcal{B}_{\mathcal{Q}}\equiv D_u\times_{\mathrm{pr}} D_v$ and call it a \textit{projective box associated with a quadrilateral $\mathcal{Q}$} (see Figure~\ref{FIG:projective_box}).
\label{DFN:projective_box_Q}
\end{dfn}

\begin{figure}
  \includegraphics[height=8cm]{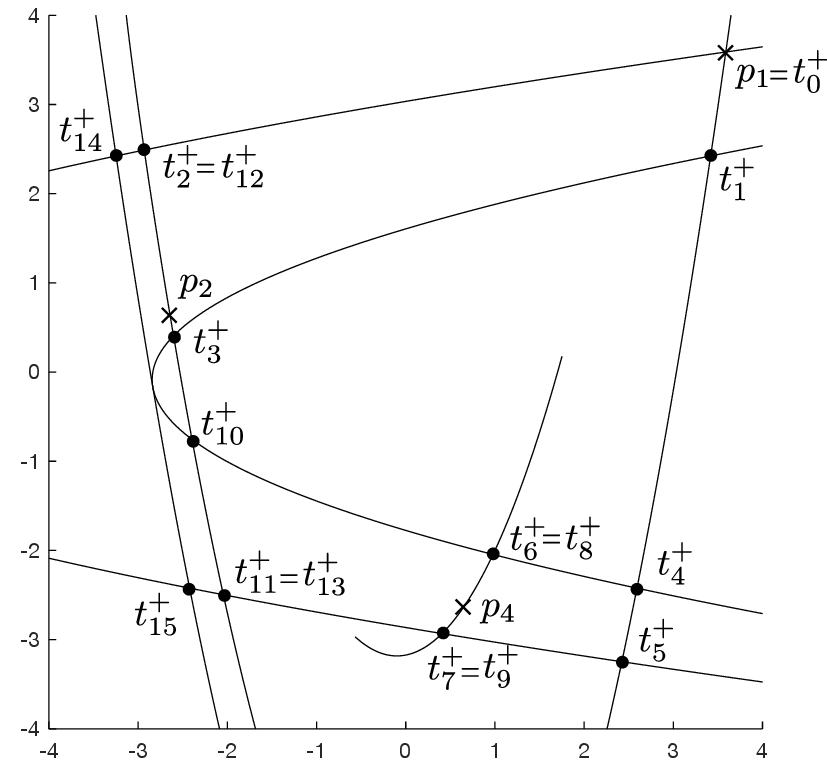} \\
  \vspace{0.5cm}
  \includegraphics[height=9cm]{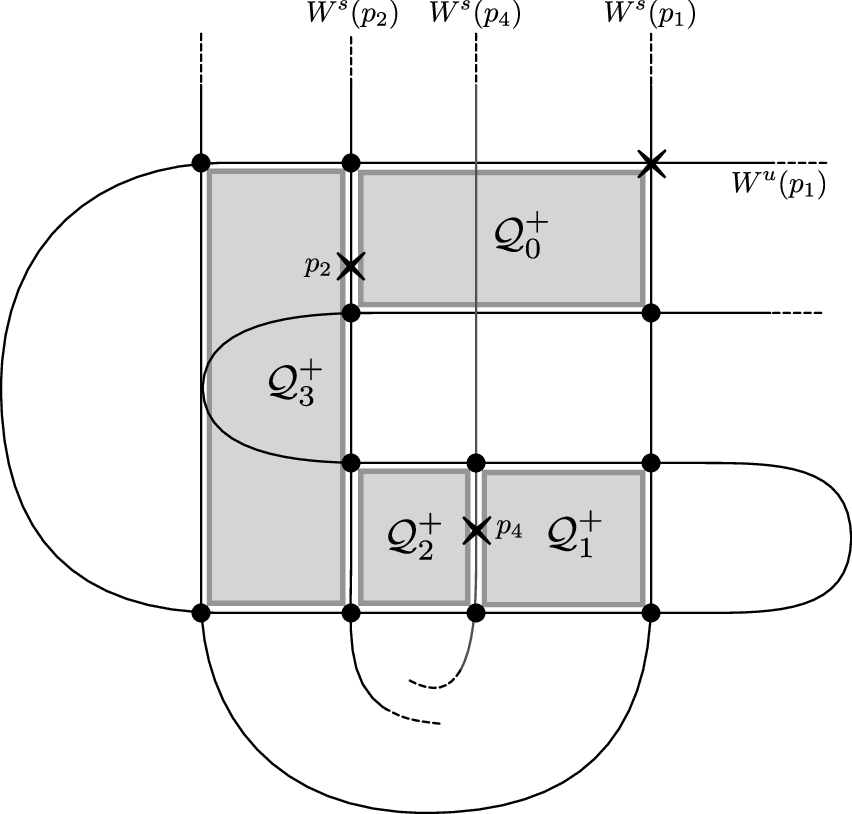}
  \caption{Above: trellis and the quadrilaterals $\{\mathcal{Q}^+_i\}^3_{i=0}$ for $(a, b)=(5.7, 1)$. Below: their cartoon images.}
  \label{FIG:trellis_positive}
\end{figure}

\begin{figure}
  \includegraphics[height=8cm]{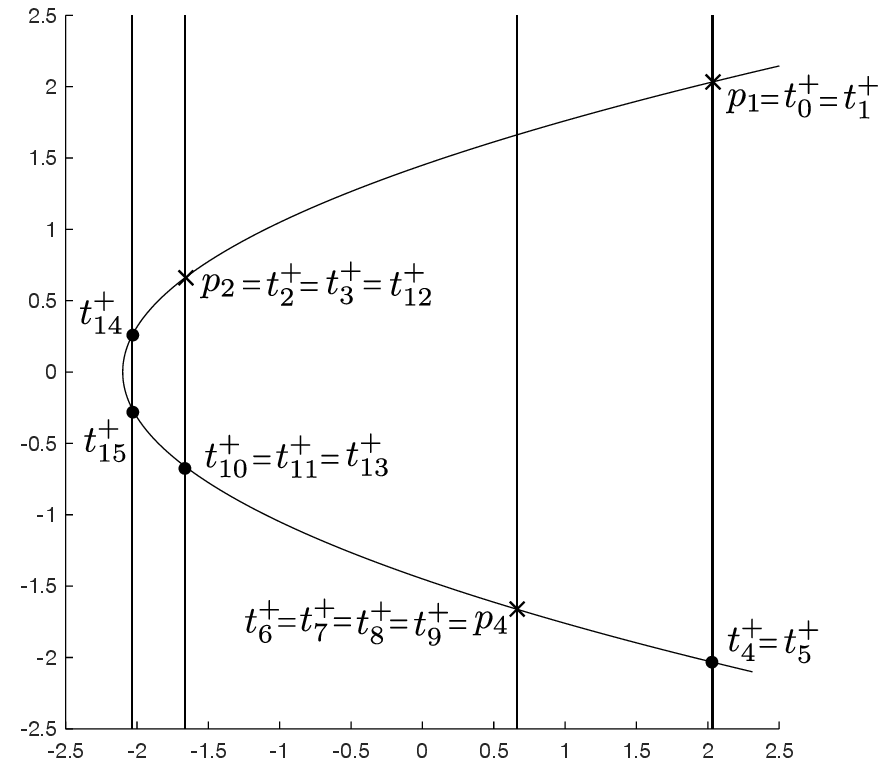} \\
  \vspace{0.5cm}
  \includegraphics[height=9.5cm]{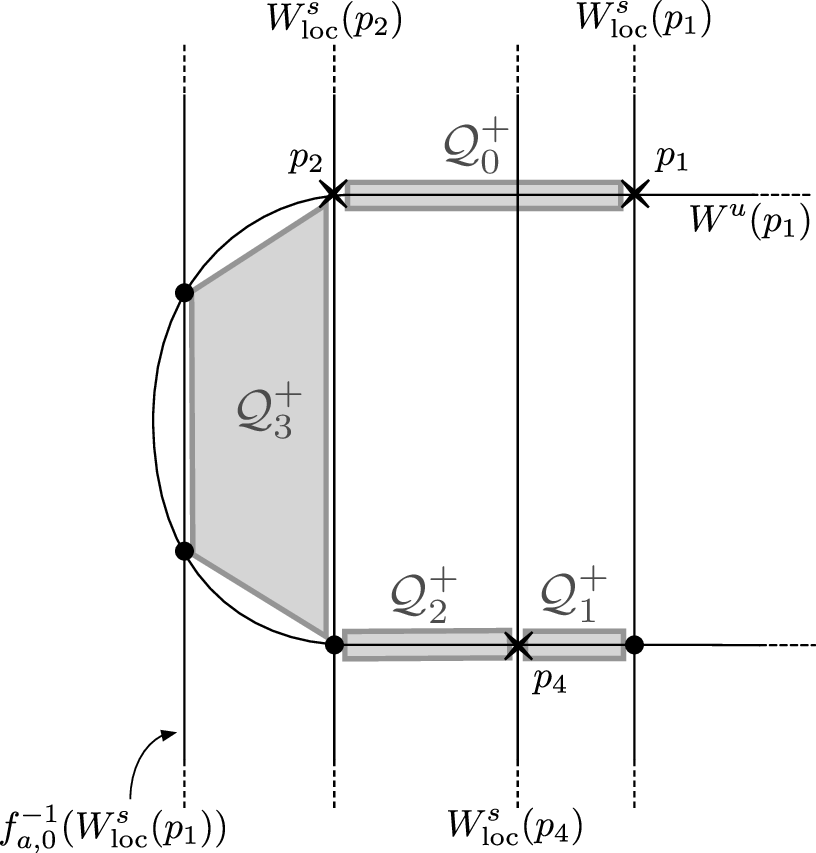}
  \caption{Above: trellis and the quadrilaterals $\{\mathcal{Q}^+_i\}^3_{i=0}$ for $(a, b)=(2.1, 0)$. Below: their cartoon images.}
  \label{FIG:trellis_positive_zero}
\end{figure}

\begin{figure}
  \includegraphics[height=8cm]{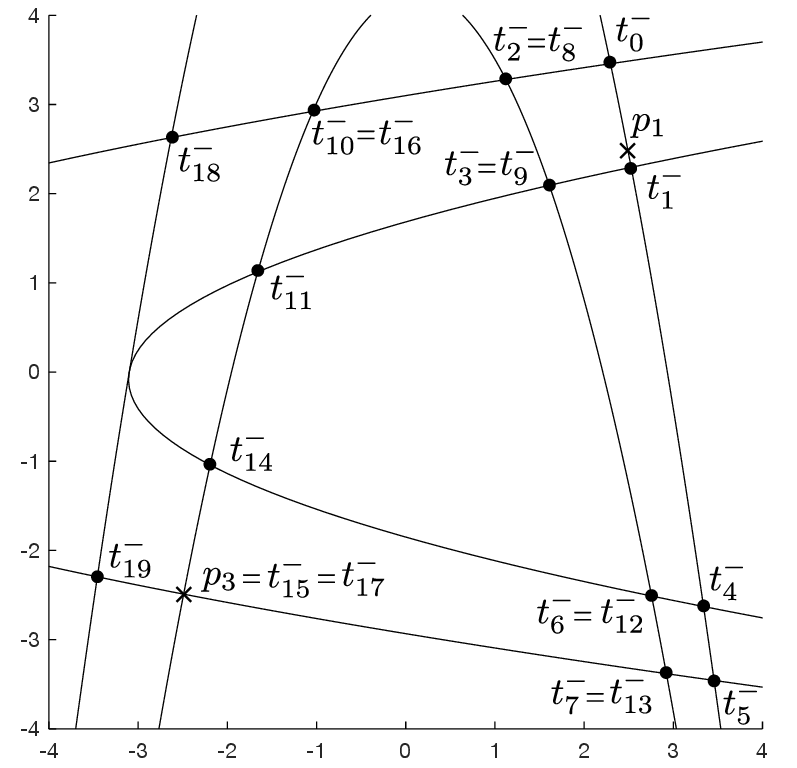} \\
  \vspace{0.5cm}
  \includegraphics[height=9cm]{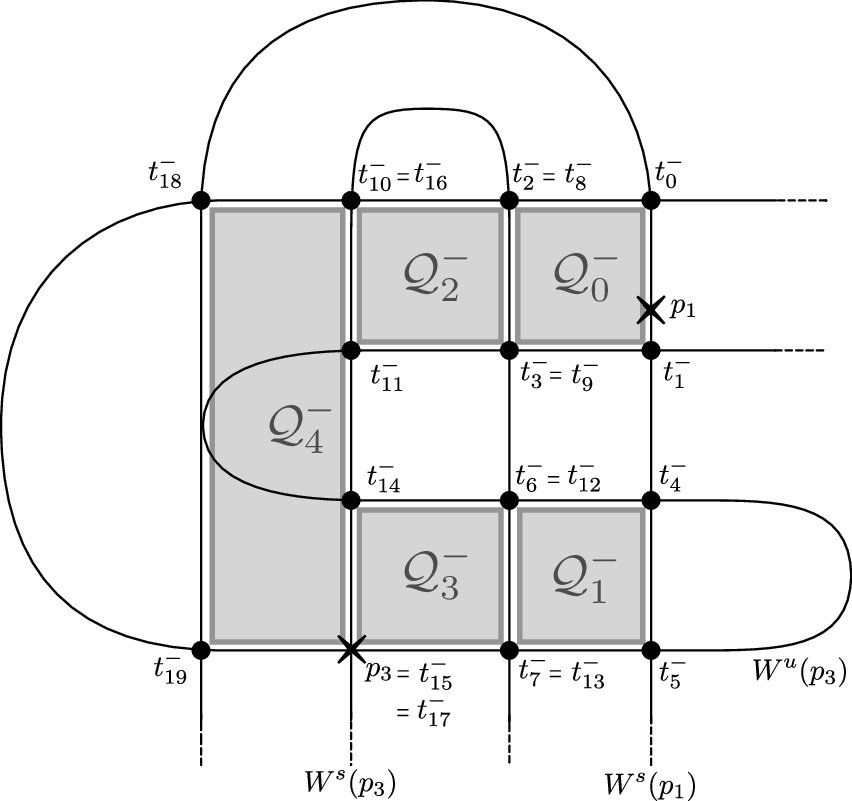}
  \caption{Above: trellis and the quadrilaterals $\{\mathcal{Q}^-_i\}^4_{i=0}$ for $(a, b)=(6.2, -1)$. Below: their cartoon images.}
  \label{FIG:trellis_negative}
\end{figure}

\begin{figure}
  \includegraphics[height=8cm]{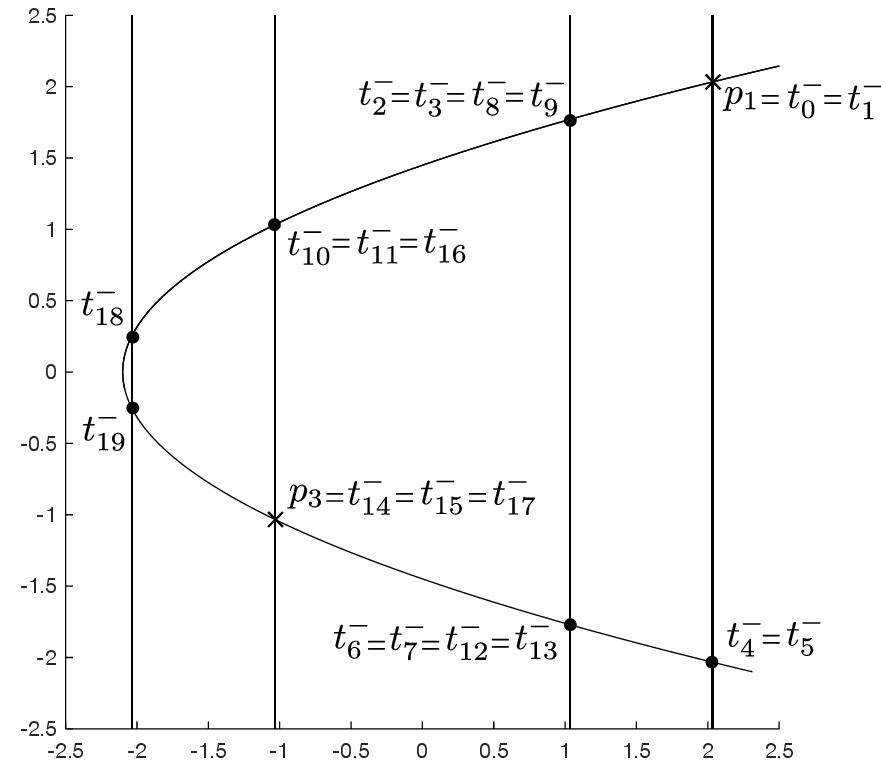} \\
  \vspace{0.5cm}
  \includegraphics[height=9.5cm]{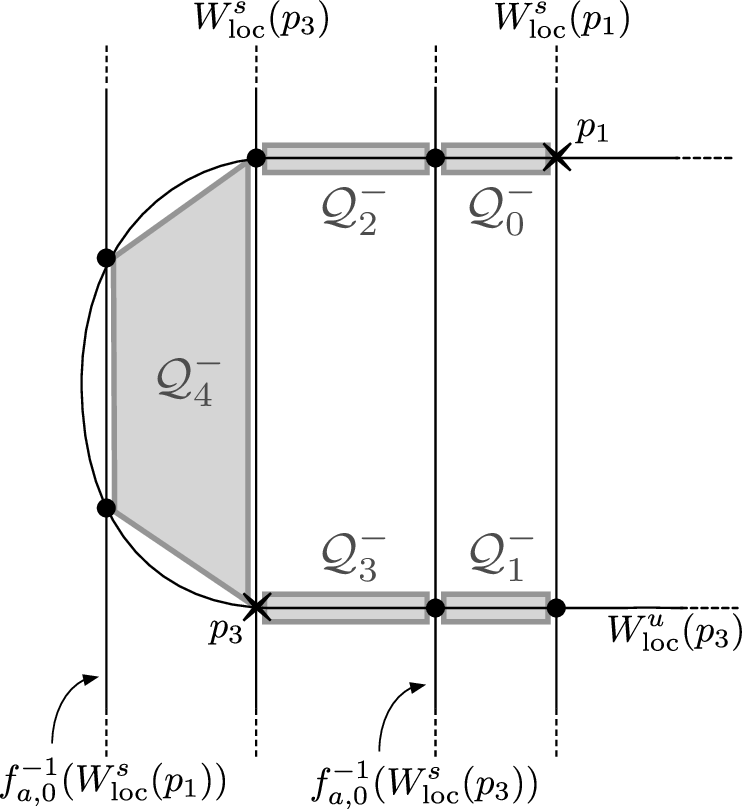}
  \caption{Above: trellis and the quadrilaterals $\{\mathcal{Q}^-_i\}^4_{i=0}$ for $(a, b)=(2.1, 0)$. Below: their cartoon images.}
  \label{FIG:trellis_negative_zero}
\end{figure}

Based on this notion, we construct a family of projective boxes associated with the trellis of $f_{\mathrm{Re}(a), \mathrm{Re}(b)}$ for $(a, b)\in \mathcal{F}^{\pm}$ as follows.

First consider the case $(a, b)\in \mathcal{F}^+$. When $\mathrm{Re}(b)\ne 0$, we compute $12$ intersection points in the trellis generated by $W^u(p_1)$, $W^s(p_1)$, $W^s(p_2)$ and $W^s(p_4)$ of the real map $f_{\mathrm{Re}(a), \mathrm{Re}(b)} : \R^2 \to \R^2$, and name them $t^+_k$ ($0\leq k\leq 15$) as in Figure~\ref{FIG:trellis_positive}. When $\mathrm{Re}(b)=0$, we compute $7$ intersection points in the trellis generated by $W^u(p_1)$, $W^s_{\mathrm{loc}}(p_1)$, $W^s_{\mathrm{loc}}(p_2)$, $W^s_{\mathrm{loc}}(p_4)$ and $f^{-1}_{\mathrm{Re}(a), 0}(W^s_{\mathrm{loc}}(p_1))$ of the real map $f_{\mathrm{Re}(a), 0} : \R^2 \to \R^2$, and name them $t^+_k$ ($0\leq k\leq 15$) as in Figure~\ref{FIG:trellis_positive_zero}. For $(a, b)\in\mathcal{F}^+$, let $\mathcal{Q}^+_i$ ($0\leq i\leq 3$) be the (possibly, degenerate\footnote{When $\mathcal{Q}^+_i$ is degenerate, we fatten it appropriately to obtain a quadrilateral; see Remark~\ref{RMK:degenerate}.}) quadrilateral in $\R^2$ formed by $t^+_{4i}$, $t^+_{4i+1}$, $t^+_{4i+2}$ and $t^+_{4i+3}$ as in Figures~\ref{FIG:trellis_positive} and \ref{FIG:trellis_positive_zero}. We define a projective box $\mathcal{B}^+_i\equiv D^+_{u, i}\times_{\mathrm{pr}} D^+_{v, i}$ associated with $\mathcal{Q}^+_i$ by choosing appropriate topological disks $D^+_{u, i}$ and $D^+_{v, i}$. See Subsection~\ref{subsection6.3} for specific data of the topological disks we will choose in Theorem~\ref{THM:quasitrichotomy} (Quasi-Trichotomy). 

Next consider the case $(a, b)\in \mathcal{F}^-$. When $\mathrm{Re}(b)\ne 0$, we compute $14$ intersection points in the trellis generated by $W^u(p_3)$, $W^s(p_1)$ and $W^s(p_3)$ of the real map $f_{\mathrm{Re}(a), \mathrm{Re}(b)} : \R^2 \to \R^2$, and name them $t^-_k$ ($0\leq k\leq 19$) as in Figure~\ref{FIG:trellis_negative}. When $\mathrm{Re}(b)=0$, we compute $8$ intersection points in the trellis generated by $W^u(p_1)$, $W^s_{\mathrm{loc}}(p_1)$, $W^s_{\mathrm{loc}}(p_3)$, $f^{-1}_{\mathrm{Re}(a), 0}(W^s_{\mathrm{loc}}(p_1))$ and $f^{-1}_{\mathrm{Re}(a), 0}(W^s_{\mathrm{loc}}(p_3))$ of the real map $f_{\mathrm{Re}(a), 0} : \R^2 \to \R^2$, and name them $t^-_k$ ($0\leq k\leq 19$) as in Figure~\ref{FIG:trellis_negative_zero}. For $(a, b)\in\mathcal{F}^-$, let $\mathcal{Q}^-_i$ ($0\leq i\leq 4$) be the (possibly, degenerate) quadrilateral in $\R^2$ formed by $t^-_{4i}$, $t^-_{4i+1}$, $t^-_{4i+2}$ and $t^-_{4i+3}$ as in Figures~\ref{FIG:trellis_negative} and \ref{FIG:trellis_negative_zero}. We define a projective box $\mathcal{B}^-_i\equiv D^-_{u, i}\times_{\mathrm{pr}} D^-_{v, i}$ associated with $\mathcal{Q}^-_i$ by choosing appropriate topological disks $D^-_{u, i}$ and $D^-_{v, i}$. See Subsection~\ref{subsection6.3} for specific data of the topological disks we will choose in Theorem~\ref{THM:quasitrichotomy} (Quasi-Trichotomy). 

\begin{dfn}
We call $\{\mathcal{B}^{\pm}_i\}_i$ a \textit{family of projective boxes associated with the trellis of $f_{\mathrm{Re}(a), \mathrm{Re}(b)}$} for $(a, b)\in \mathcal{F}^{\pm}$. 
\end{dfn}

This kind of a family of boxes has been first used in~\cite{BS2} and also employed to construct the first example of a non-planar hyperbolic H\'enon map~\cite{I1} as well as certain combinatorial objects called the Hubbard trees in~\cite{I2} and the iterated monodromy groups~\cite{I3} for such maps.

\subsection{Crossed mappings}\label{subsection2.3}

The notion of a crossed mapping has been first introduced in~\cite{HO} and will play a crucial role throughout this paper. Here we present the following version of this notion (see Subsection 5.1 in~\cite{ISm}). 

Let $\mathcal{B}=D_u\times_{\mathrm{pr}}D_v$ (resp. $\mathcal{B}'=D'_u\times_{\mathrm{pr}}D'_v$) be a projective box and let $(\pi_u, \pi_v)$ (resp. $(\pi'_u, \pi'_v)$) be the projective coordinates for $\mathcal{B}$ (resp. $\mathcal{B}'$). 

\begin{dfn}[\textbf{Crossed Mapping Condition}]
We say that $f : \mathcal{B}\cap f^{-1}(\mathcal{B}')\to \mathcal{B}'$ satisfies the \textit{crossed mapping condition (CMC)} of degree $d$ if
\[\rho_f\equiv (\pi'_u \circ f, \pi_v \circ \iota) \, : \, \mathcal{B}\cap f^{-1}(\mathcal{B}') \longrightarrow D'_u\times D_v\]
is proper of degree $d$, where $\iota : \mathcal{B}\cap f^{-1}(\mathcal{B}')\to \mathcal{B}$ is the inclusion map. 
\label{DFN:CMC}
\end{dfn}

Let $\mathcal{B}$, $\mathcal{B}'$ and $\mathcal{B}''$ be projective boxes. A proof of the next claim can be found in Proposition 3.7 (b) of~\cite{HO}.

\begin{lmm}
Let $f : \mathcal{B}\cap f^{-1}(\mathcal{B}') \to \mathcal{B}'$ (resp. $g : \mathcal{B}'\cap g^{-1}(\mathcal{B}'') \to \mathcal{B}''$) satisfy the (CMC) of degree $d_f$ (resp. degree $d_g$). Then, the composition $g\circ f : \mathcal{B}\cap f^{-1}(\mathcal{B}'\cap g^{-1}(\mathcal{B}''))\to \mathcal{B}''$ satisfies the (CMC) of degree $d_f\, d_g$.
\label{LMM:composition}
\end{lmm}

Let $\mathcal{B}=D_u\times_{\mathrm{pr}}D_v$ be a projective box. 

\begin{dfn}
A complex one-dimensional (not necessarily connected) submanifold $D$ in $\mathcal{B}$ is called \textit{horizontal} \footnote{We remark that the notion of a horizontal (resp. vertical) submanifold defined here is weaker than a \textit{horizontal-like} (resp. \textit{vertical-like}) submanifold given in~\cite{ISm}. Any tangent vector to a horizontal-like (resp. vertical-like) submanifold is contained in the horizontal (resp. vertical) Poincar\'e cone (see Definition 5.7 in~\cite{ISm}), but a tangent vector to a horizontal (resp. vertical) submanifold can be vertical (resp. horizontal).} of degree $d$ if the projection $\pi_u : D\to D_u$ is a proper map of degree $d$. The notion of a \textit{vertical submanifold} is defined similarly.
\label{DFN:horizontal_submanifold}
\end{dfn}

The next lemma tells that a crossed mapping controls the behavior of horizontal/vertical submanifolds under $f$. A proof can be found in Proposition 3.4 of~\cite{HO}.

\begin{lmm}
If $f : \mathcal{B}\cap f^{-1}(\mathcal{B}')\to \mathcal{B}'$ satisfies the (CMC) of degree $d$ and if $D\subset \mathcal{B}$ is a horizontal submanifold of degree $k$, then $f(D)\cap \mathcal{B}'$ is horizontal of degree $dk$ in $\mathcal{B}'$. If $f^{-1} : \mathcal{B}'\cap f(\mathcal{B})\to \mathcal{B}$ satisfies the (CMC) of degree $d$ and if $D\subset \mathcal{B}'$ is a vertical submanifold of degree $k$, then $f^{-1}(D)\cap \mathcal{B}$ is vertical of degree $dk$ in $\mathcal{B}$.
\label{LMM:mapping_submanifolds}
\end{lmm}

We note that in Lemma~\ref{LMM:mapping_submanifolds} above, the submanifold $f(D)\cap \mathcal{B}'$ may not be connected even when $D$ is connected. 

\begin{figure}
  \includegraphics[height=7cm]{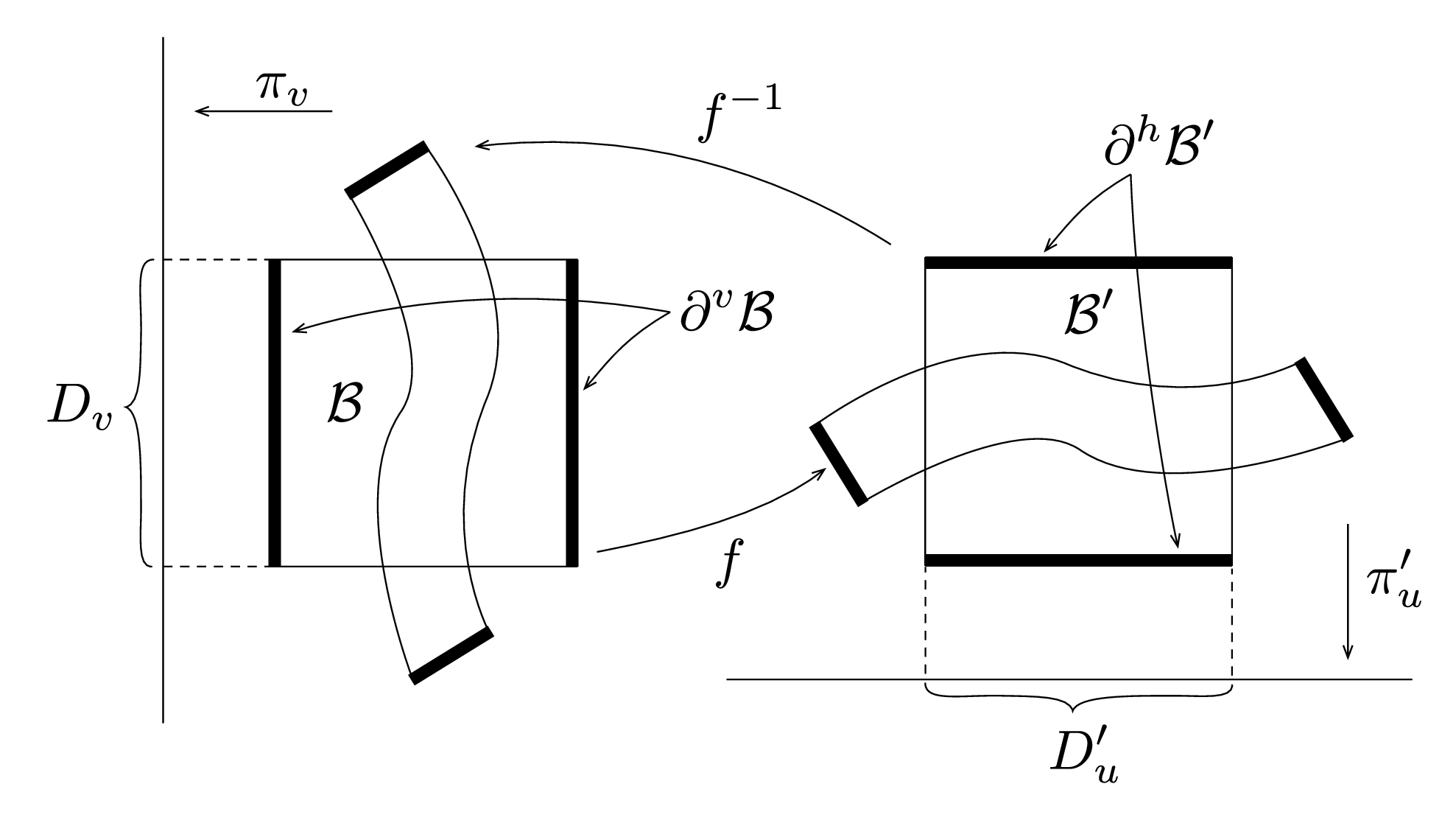}
  \caption{Figure of the boundary compatibility condition.}
  \label{FIG:BCC}
\end{figure}

A more checkable condition for a map to satisfy the (CMC) is given as follows (see Subsection 5.2 in~\cite{I1}). Below we write $\partial^v\mathcal{B}\equiv \partial D_u\times_{\mathrm{pr}} D_v$ and $\partial^h\mathcal{B}\equiv D_u\times_{\mathrm{pr}} \partial D_v$ for $\mathcal{B}=D_u\times_{\mathrm{pr}} D_v$.

\begin{dfn}
We say that $f : \C^2\to \C^2$ satisfies the \textit{boundary compatibility condition (BCC)} with respect to $\mathcal{B}$ and $\mathcal{B}'$ if both $\pi'_u\circ f(\partial^v\mathcal{B})\cap D'_u=\emptyset$ and $\pi_v\circ f^{-1}(\partial^h\mathcal{B}')\cap D_v=\emptyset$ hold (see Figure~\ref{FIG:BCC}). 
\label{DFN:BCC}
\end{dfn}

Note that this last condition $\pi_v\circ f^{-1}(\partial^h\mathcal{B}')\cap D_v=\emptyset$ makes sense even when $f^{-1}$ is not defined; it can be replaced by $f(\mathcal{B})\cap\partial^h\mathcal{B}'=\emptyset$. 

\begin{figure}
  \includegraphics[height=4.5cm]{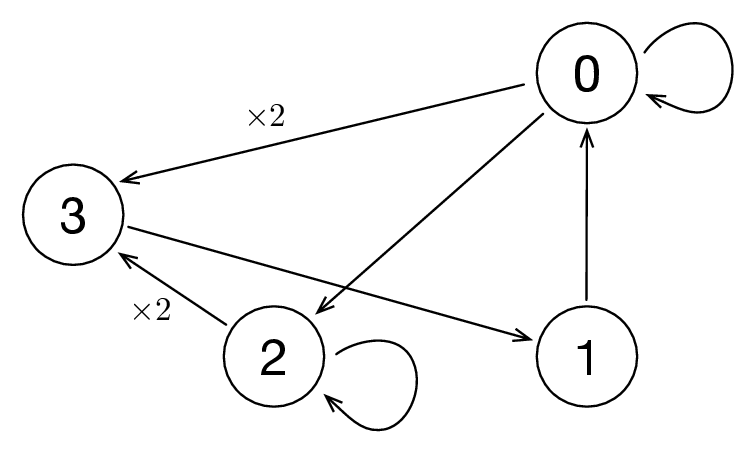}
  \caption{Diagram for admissible transitions $\mathfrak{T}^+$ for $(a, b)\in\mathcal{F}^+$.}
\label{FIG:admissible_positive}
\end{figure}

\begin{figure}
  \includegraphics[height=4.5cm]{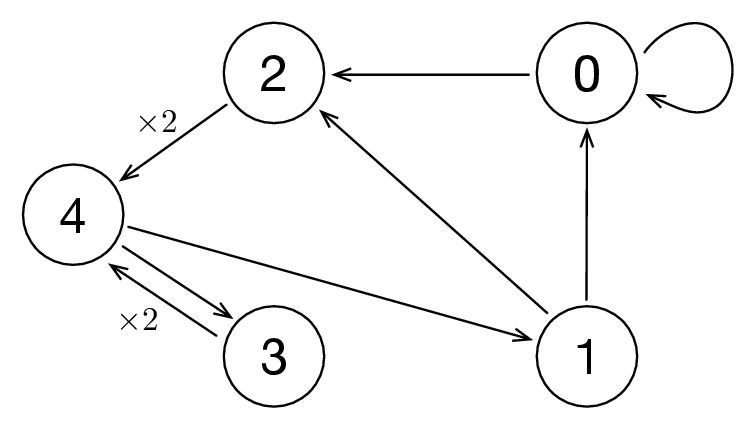}
  \caption{Diagram of admissible transitions $\mathfrak{T}^-$ for $(a, b)\in\mathcal{F}^-$.}
\label{FIG:admissible_negative}
\end{figure}

Below we give an explicit family of four projective boxes $\{\mathcal{B}^+_i\}_{i=0}^3$ for every parameter $(a, b)\in\mathcal{F}^+$ and a family of five projective boxes $\{\mathcal{B}^-_i\}_{i=0}^4$ for every parameter $(a, b)\in\mathcal{F}^-$. We set
\[\mathfrak{T}^+\equiv\bigl\{(0, 0), (0, 2), (0, 3), (1, 0), (2, 2), (2, 3), (3, 1) \bigr\}\]
and 
\[\mathfrak{T}^-\equiv\bigl\{(0, 0), (0, 2), (1, 0), (1, 2), (2, 4), (3, 4), (4, 1), (4, 3) \bigr\}.\]
Elements in $\mathfrak{T}^{\pm}$ are called \textit{admissible transitions}. 

\begin{figure}
  \includegraphics[height=7.5cm]{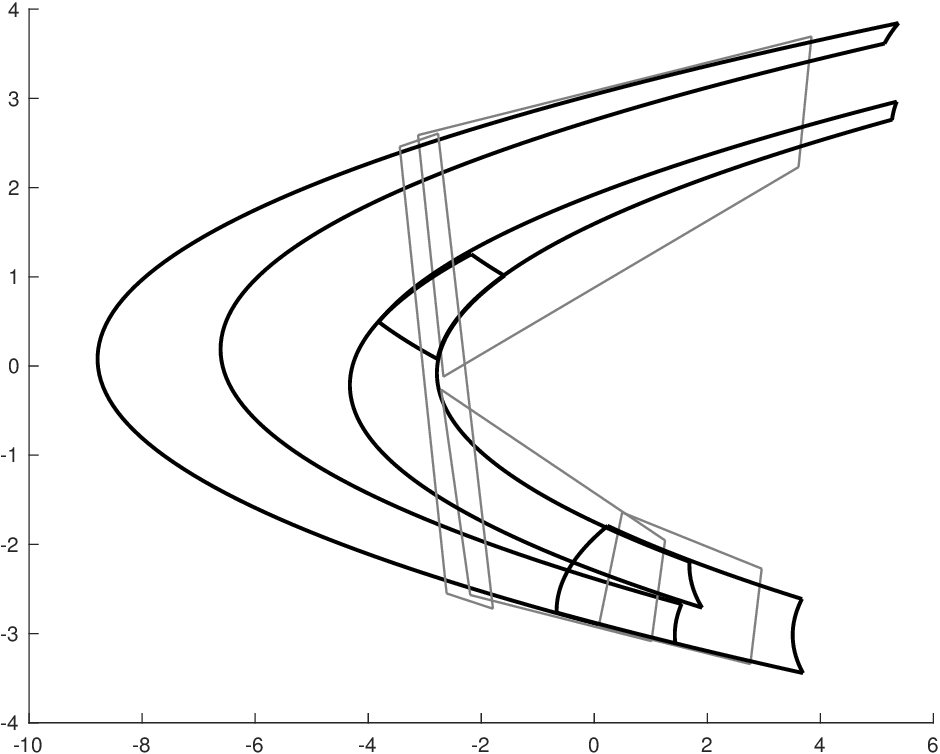} \\
  \vspace{0.5cm}
  \includegraphics[height=7.5cm]{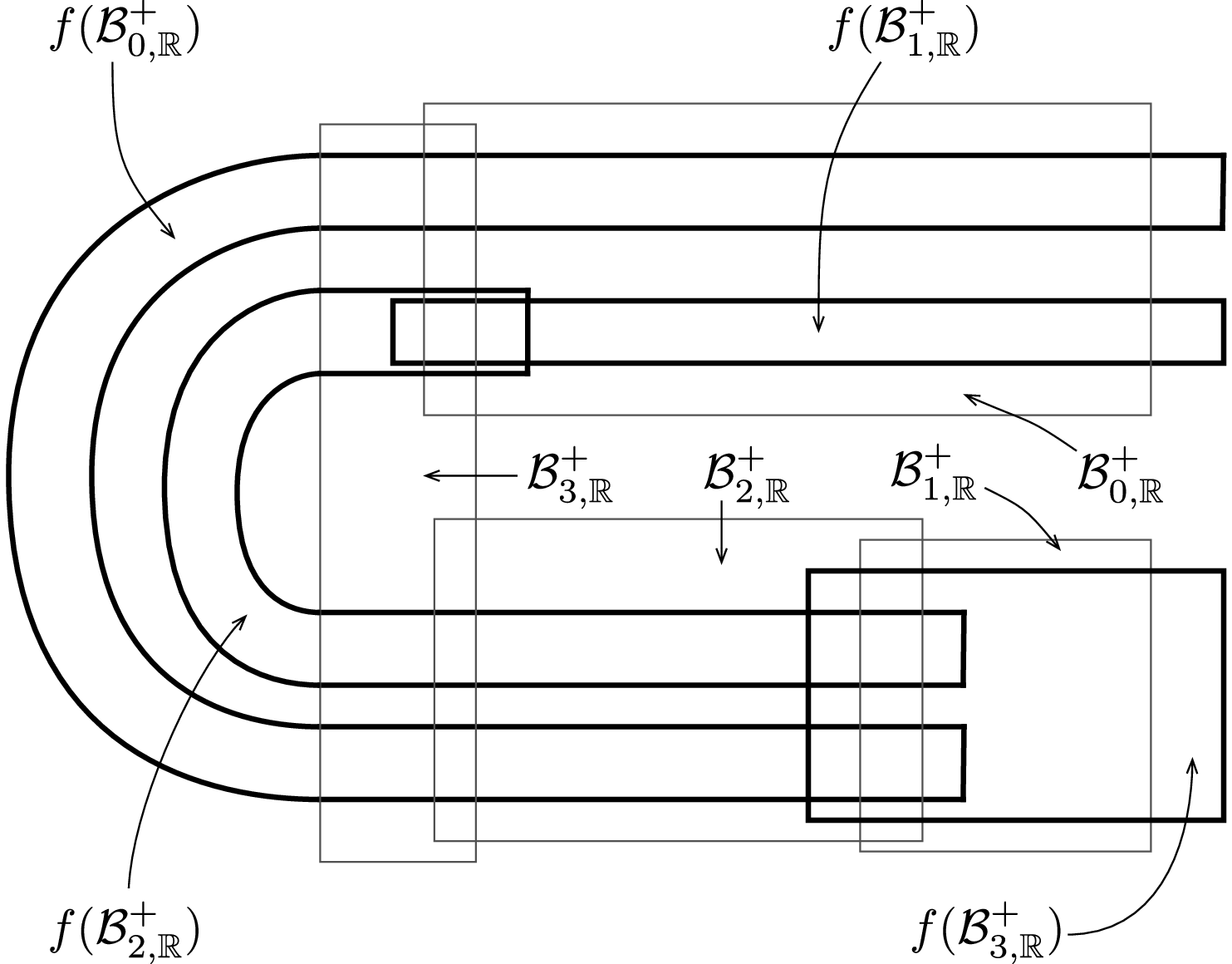}
  \caption{Above: the real slices of the boxes $\mathcal{B}^+_{i, \R}=\mathcal{B}^+_i\cap\R^2$ and their images by $f_{a, b}$ for $(a, b)=(5.7, 1)$. Below: their cartoon images.}
\label{FIG:slice_positive}
\end{figure}

\begin{figure}
  \includegraphics[height=7.5cm]{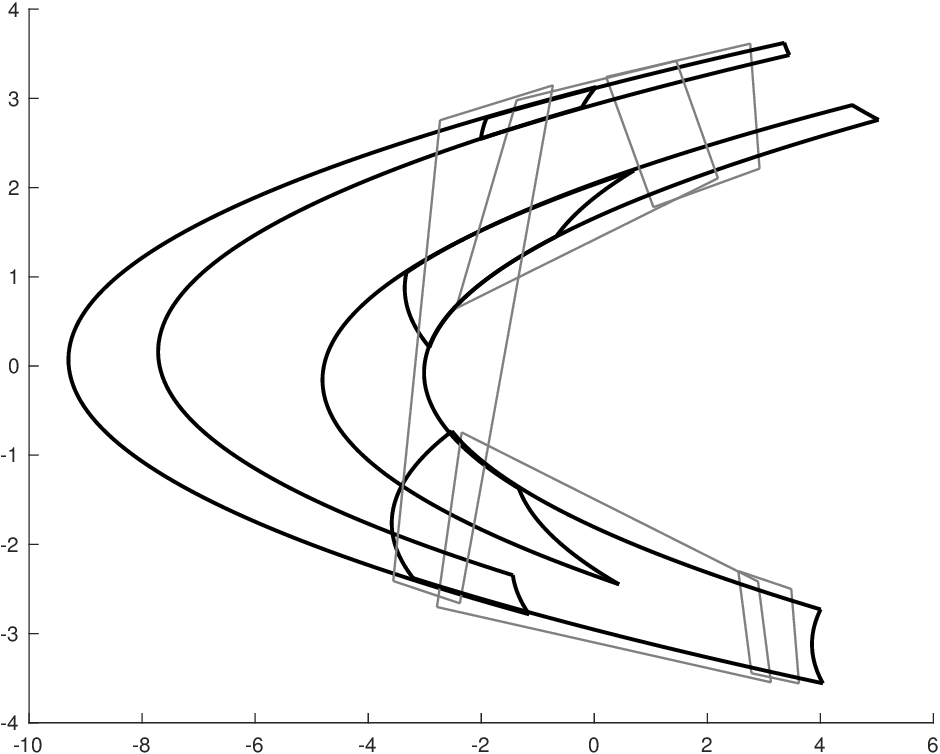} \\
  \vspace{0.5cm}
  \includegraphics[height=7.5cm]{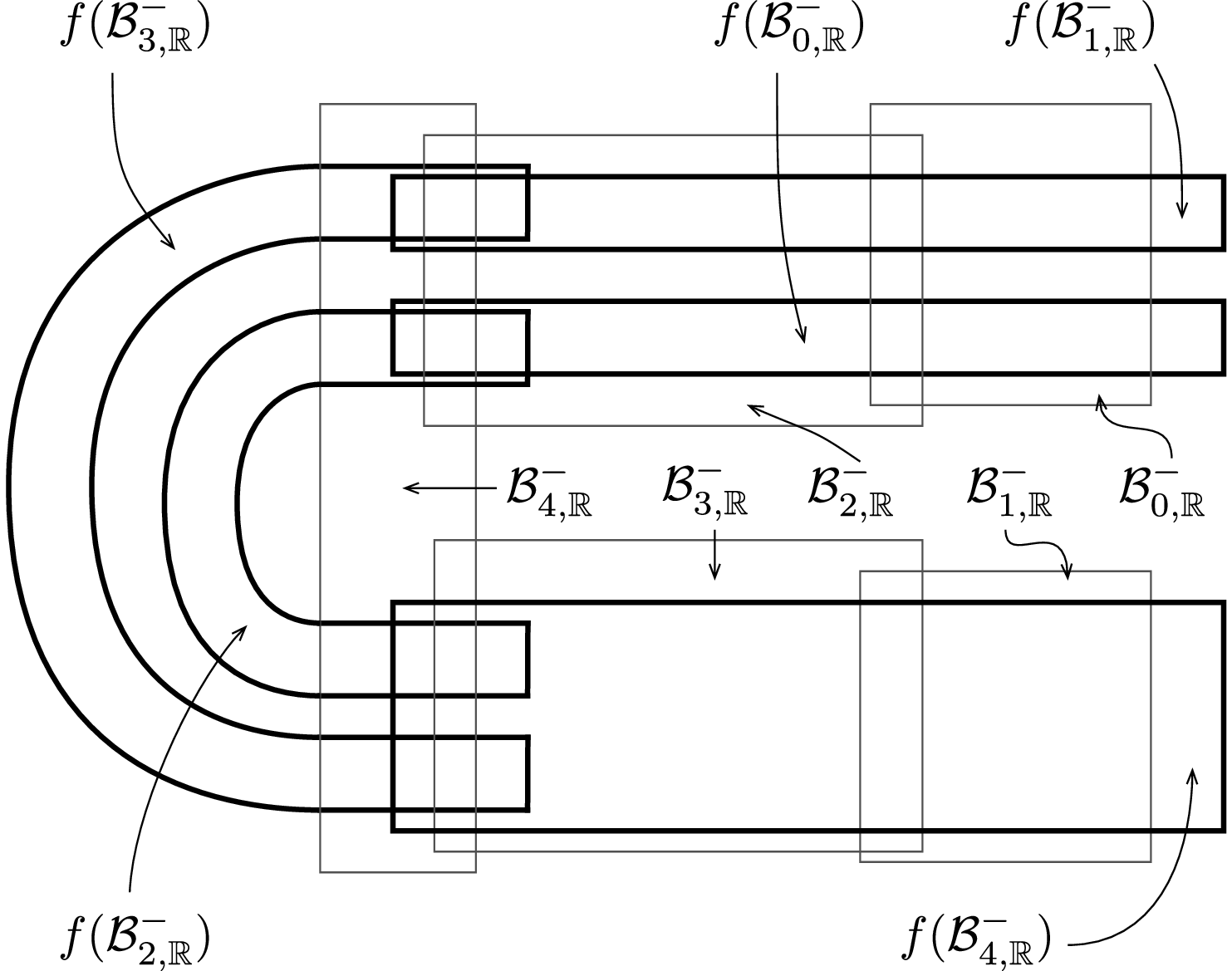}
  \caption{Above: the real slices of the boxes $\mathcal{B}^-_{i, \R}=\mathcal{B}^-_i\cap\R^2$ and their images by $f_{a, b}$ for $(a, b)=(6.2, -1)$. Below: their cartoon images.}
\label{FIG:slice_negative}
\end{figure}

\begin{dfn}
A triple $(f_{a, b}, \{\mathcal{B}^{\pm}_i\}, \mathfrak{T}^{\pm})$ is said to satisfy the (CMC) if $f_{a, b} : \mathcal{B}^{\pm}_i\cap f_{a, b}^{-1}(\mathcal{B}^{\pm}_j)\to \mathcal{B}^{\pm}_j$ satisfies the (CMC) for every $(i, j)\in \mathfrak{T}^{\pm}$. 
\label{DFN:triple}
\end{dfn}

Diagrams in Figures~\ref{FIG:admissible_positive} and \ref{FIG:admissible_negative} describe all the admissible transitions $\mathfrak{T}^+$ and $\mathfrak{T}^-$ respectively, where $\times 2$ indicates that the corresponding transition is a crossed mapping of degree $2$. 

Figures~\ref{FIG:slice_positive} and \ref{FIG:slice_negative} illustrate how the real slices of the boxes $\mathcal{B}^{\pm}_i$ we will choose in Theorem~\ref{THM:quasitrichotomy} are mapped by a H\'enon map $f=f_{a, b}$. There, by comparing with the cartoon figures below, one can see how the boxes in the real figures above are mapped by $f$.

\subsection{Quasi-trichotomy}\label{subsection2.4}

The purpose of this subsection is to classify any H\'enon map into three types; either (i) $f_{a, b}$ does not attain the maximal entropy, (ii) $f_{a, b}$ is a hyperbolic horseshoe on $\R^2$, or (iii) $f_{a, b}$ satisfies the crossed mapping condition. More precisely, we show

\begin{thm}[\textbf{Quasi-Trichotomy}]
We have the following three claims. 
\begin{enumerate}
 \renewcommand{\labelenumi}{(\roman{enumi})}
  \item If $(a, b)\in \R\times (I^{\pm}_{\R}\setminus \{0\})$ and $a\leq a^{\pm}_{\mathrm{aprx}}(b)-\chi^{\pm}(b)$, we have $\htop(f_{a, b}|_{\R^2})<\log 2$.
  \item If $(a, b)\in\R\times (I^{\pm}_{\R}\setminus \{0\})$ and $a\geq a^{\pm}_{\mathrm{aprx}}(b)+\chi^{\pm}(b)$, $f_{a, b}$ is a hyperbolic horseshoe on $\R^2$.
  \item If $(a, b)\in \C\times I^{\pm}$ and $|a-a^{\pm}_{\mathrm{aprx}}(b)|\leq \chi^{\pm}(b)$, one can construct a family of projective boxes $\{\mathcal{B}^{\pm}_i\}_i$ associated with the trellis of $f_{\mathrm{Re}(a), \mathrm{Re}(b)}$ so that $(f_{a, b}, \{\mathcal{B}^{\pm}_i\}_i, \mathfrak{T}^{\pm})$ satisfies the (CMC).
\end{enumerate} 
\label{THM:quasitrichotomy}
\end{thm}

The proof of Theorem~\ref{THM:quasitrichotomy} (Quasi-Trichotomy) requires computer-assistance with rigorous error bounds. Notice that the condition for $(a, b)$ in (iii) is equivalent to $(a, b)\in \mathcal{F}^{\pm}$, hence $a$ and $b$ are allowed to be complex numbers and $b$ can vanish. Note also that the three cases (i), (ii) and (iii) are not exclusive, and this is why we call this theorem ``Quasi-Trichotomy''. 

\begin{rmk}
In our computer-assisted proofs below, the compactness of parameter regions and dynamical regions where we verify numerical criteria is essential since only finitely many statements described in terms of compact intervals can be checked by interval arithmetic. For example, we verify certain numerical condition by computer-assistance for $(a, b) \in \R\times I^{\pm}_{\R}$ with $-(b+1)^2/4 \leq a \leq a^{\pm}_{\mathrm{aprx}} (b) - \chi^{\pm}(b)$ in Lemma~\ref{LMM:period7}, and for $(a, b)\in\R\times I^{\pm}_{\R}$ with $2(1+|b|)^2 \geq a \geq a^{\pm}_{\mathrm{aprx}}(b) + \chi^{\pm}(b)$ in Lemma~\ref{LMM:zero-section} (note that both regions contain the case $b=0$).
\end{rmk}

\begin{proof}[Proof of (i) of Theorem~\ref{THM:quasitrichotomy}]
Recall Theorem 10.1 in~\cite{BLS} which proves that $\htop(f_{a, b}|_{\R^2})=\log 2$ if and only if every periodic point of $f_{a, b}: \C^2 \to \C^2$ is contained in $\R^2$ for $(a, b) \in \R \times \R^{\times}$. Therefore, it suffices to show that there exists a periodic point of $f_{a, b}$ in $\C^2 \setminus \R^2$ for all $(a, b) \in \R\times (I^{\pm}_{\R}\setminus \{0\})$ with $a \leq a^{\pm}_{\mathrm{aprx}}(b)-\chi^{\pm}(b)$.

For $a$ small enough, this can be done by hand; if $a < -(b+1)^2/4$, the two fixed points of $f_{a, b}$ are away from $\R^2$ by solving the quadratic equation defining the fixed point of the map. For the rest of the parameter values, the existence of a non-real periodic point is established by rigorous numerics. In fact, in Subsection~\ref{subsection6.4} we show 

\begin{lmm}
For all $(a, b) \in \R\times I^{\pm}_{\R}$ with $-(b+1)^2/4 \leq a \leq a^{\pm}_{\mathrm{aprx}}(b)-\chi^{\pm}(b)$, there exists a periodic point of period $7$ of $f_{a, b}$ in $\C^2 \setminus \R^2$. 
\label{LMM:period7}
\end{lmm}

The proof first uses Newton's method to find an approximate periodic point in $\C^2 \setminus \R^2$ and then its existence is rigorously proven by the interval Krawczyk method. Remark that the statement of the lemma includes the case $b = 0$, in which $f_{a, b}$ degenerate to the one-dimensional quadratic map. The periodic point continues to the case $b = 0$ and remains in $\C^2 \setminus \R^2$. This completes the proof of the claim (i).
\end{proof}

\begin{proof}[Proof of (ii) of Theorem~\ref{THM:quasitrichotomy}]
We first prove that for $(a, b) \in \R\times (I^{\pm}_{\R}\setminus \{0\})$ with $a > 2(1+|b|)^2$, $f_{a, b}$ is a hyperbolic horseshoe on $\R^2$. Under the assumption $|a|>2(1+|b|)^2$, it has been shown that the restriction of $f_{a, b}$ to its complex non-wandering set $\mathrm{\Omega}(f_{a, b})$ is topologically conjugate to the shift map $\sigma$ on $\{0, 1\}^{\Z}$ (see~\cite{O, U}), and that $f_{a, b}$ is hyperbolic on $\mathrm{\Omega}(f_{a, b})$ (see~\cite{ISm}). Hence our task is to prove that $\mathrm{\Omega}(f_{a, b})$ is contained in $\R^2$ when $a > 2(1+|b|)^2$. 

To do this we first recall the following construction in~\cite{O, U, ISm}. Let us put  
\[R\equiv \frac{1+|b|+\sqrt{(1+|b|)^2+4|a|}}{2}\] 
and define
\[\mathcal{D}\equiv \bigl\{(x, y)\in\C^2 : |x|\leq R, |y|\leq R\bigr\}.\]
Then, we see that $\mathcal{D}\cap f^{-1}(\mathcal{D})$ consists of two connected components, say $\mathcal{D}_0$ and $\mathcal{D}_1$. Given a symbol sequence $\underline{\varepsilon}=\cdots\varepsilon_{-2}\varepsilon_{-1}\cdot\varepsilon_0\varepsilon_1\cdots\in\{0, 1\}^{\Z}$, $\bigcap_{n\geq 0}f^n(\mathcal{D}_{\varepsilon_{-n}})$ is a horizontal submanifold of degree one in $\mathcal{D}$ and $\bigcap_{n\leq 0}f^n(\mathcal{D}_{\varepsilon_{-n}})$ is a vertical submanifold of degree one in $\mathcal{D}$. Therefore, their (complex) intersection $\bigcap_{n\in\Z}f^n(\mathcal{D}_{\varepsilon_{-n}})$ consists of exactly one point which we denote by $\omega(\underline{\varepsilon})\in\mathrm{\Omega}(f_{a, b})$. 

Next we consider their real sections, namely we define $\mathcal{D}_{\R}\equiv\mathcal{D}\cap\R$. Then, $\mathcal{D}_{\R}\cap f(\mathcal{D}_{\R})$ consists of two connected components, say $\mathcal{D}_{\R, 0}$ and $\mathcal{D}_{\R, 1}$, each of which is a strip connecting the left boundary and the right boundary of the square $\mathcal{D}_{\R}$. Now, take a symbol sequence $\underline{\varepsilon}=\cdots\varepsilon_{-2}\varepsilon_{-1}\cdot\varepsilon_0\varepsilon_1\cdots\in\{0, 1\}^{\Z}$. Then, for any $N\geq 0$ one can inductively show that $\bigcap_{0\leq n\leq N}f^n(\mathcal{D}_{\R, \varepsilon_{-n}})$ is a strip connecting the left boundary and the right boundary of the square $\mathcal{D}_{\R}$. A similar argument shows that $\bigcap_{-N\leq n\leq 0}f^n(\mathcal{D}_{\R, \varepsilon_{-n}})$ is a strip connecting the upper boundary and the lower boundary of the square $\mathcal{D}_{\R}$. Therefore, $\bigcap_{-N\leq n\leq N}f^n(\mathcal{D}_{\R, \varepsilon_{-n}})$ is a decreasing sequence in $N$ of non-empty compact sets. It follows from the compactness that $\bigcap_{n\in\Z}f^n(\mathcal{D}_{\R, \varepsilon_{-n}})$ is non-empty. 

Since we have
\[1=\mathrm{Card}\Bigg(\bigcap_{n\in\Z}f^n(\mathcal{D}_{\varepsilon_{-n}})\Bigg) \geq \mathrm{Card}\Bigg(\bigcap_{n\in\Z}f^n(\mathcal{D}_{\R, \varepsilon_{-n}})\Bigg)\geq 1,\]
it follows that the real intersection $\bigcap_{n\in\Z}f^n(\mathcal{D}_{\R, \varepsilon_{-n}})$ consists of exactly one point and hence it coincides with the complex intersection $\omega(\underline{\varepsilon})$. Since $\omega(\{0, 1\}^{\Z})=\mathrm{\Omega}(f_{a, b})$, this yields that $\mathrm{\Omega}(f_{a, b})\subset \R^2$ and $f_{a, b}$ is a hyperbolic horseshoe on $\R^2$.

For the rest of the parameters $(a, b) \in \R\times (I^{\pm}_{\R}\setminus \{0\})$ with $2(1+|b|)^2 \geq a \geq a^{\pm}_{\mathrm{aprx}}(b) + \chi^{\pm}(b)$ we employ the algorithm of~\cite{A1}. The key step is to prove the uniform hyperbolicity of the map. To avoid the difficulty in defining unstable and stable directions, we introduced a weaker notion of hyperbolicity called quasi-hyperbolicity. Let $f : M\to M$ be a smooth map on a differentiable manifold $M$ and $\mathrm{\Lambda}\subset M$ be a compact invariant set of $f$. We denote by $T_{\mathrm{\Lambda}}M$ the restriction of the tangent bundle $TM$ to $\mathrm{\Lambda}$. An orbit of $Df|_{T_{\mathrm{\Lambda}}M}: T_{\mathrm{\Lambda}}M \to T_{\mathrm{\Lambda}}M$ is said to be \textit{trivial} if it is contained in the image of the zero section of $T_{\mathrm{\Lambda}}M$.

\begin{dfn}
We say that $f$ is \textit{quasi-hyperbolic} on $\mathrm{\Lambda}$ if the restriction $Df|_{T_{\mathrm{\Lambda}}M} : T_{\mathrm{\Lambda}}M\to T_{\mathrm{\Lambda}}M$ has no non-trivial bounded orbit, that is, the orbit of every non-zero tangent vector is unbounded with respect to either forward or backward iteration of $Df|_{T_{\mathrm{\Lambda}}M}$.
\label{DFN:quasi-hyperbolic}
\end{dfn}

It is known that quasi-hyperbolicity is strictly weaker than uniform hyperbolicity. However, when the invariant set $\mathrm{\Lambda}$ is the chain recurrent set of the map, these two notions of hyperbolicity coincide~\cite{CFS, SS} (see also Theorem 2.3 of~\cite{A1}). Recall that the \textit{chain recurrent set} $\mathcal{R}(f)$ of $f : M \to M$ is the set of points $x\in M$ such that for any $\varepsilon > 0$ there exists an $\varepsilon$-chain from $x$ to itself. Here, an \textit{$\varepsilon$-chain} from $x$ to $x'$ is a sequence of points $x = x_0, x_1, \ldots, x_n = x'$ satisfying $d(f(x_i), x_{i+1}) < \varepsilon$ for $0\leq i \leq n-1$, where $d$ is the distance function on $M$. Therefore, to show the uniform hyperbolicity of $f$ on $\mathcal{R}(f)$, it suffices to show the quasi-hyperbolicity on $\mathcal{R}(f)$. To do this, it is convenient to rephrase the definition of quasi-hyperbolicity in terms of an isolating neighborhood as follows. Let $N \subset M$ be a compact set. Its \textit{maximal invariant set} $\mathrm{Inv}(f, N)$ is defined as
\[\mathrm{Inv}(f, N) = \bigcap_{n\in\Z}f^n(N).\]
Note that this definition is valid even for non-invertible maps. A compact set $N\subset M$ is called an \textit{isolating neighborhood} with respect to $f$ if $\mathrm{Inv}(f, N)$ is contained in the interior of $N$.

\begin{prp}
Assume that $N\subset T_{\mathrm{\Lambda}}M$ is an isolating neighborhood with respect to $Df|_{T_{\mathrm{\Lambda}}M} : T_{\mathrm{\Lambda}}M \to T_{\mathrm{\Lambda}}M$ and $N$ contains the image of the zero-section of $T_{\mathrm{\Lambda}}M$. Then $f$ is quasi-hyperbolic on $\mathrm{\Lambda}$.
\label{PRP:quasi-hyperbolic}
\end{prp}

See Proposition 2.5 in~\cite{A1} for a proof. With the help of rigorous numerics combined with set-oriented algorithms, we show 

\begin{lmm}
For all $(a, b)\in\R\times I^{\pm}_{\R}$ with $2(1+|b|)^2 \geq a \geq a^{\pm}_{\mathrm{aprx}}(b) + \chi^{\pm}(b)$, one can find an isolating neighborhood $N\subset T _{\mathcal{R}(f)}\R^2$ with respect to $Df|_{T_{\mathcal{R}(f)}\R^2} : T_{\mathcal{R}(f)}\R^2\to T_{\mathcal{R}(f)}\R^2$ containing the image of the zero-section of $T_{\mathcal{R}(f)}\R^2$, where $f=f_{a, b}$.
\label{LMM:zero-section}
\end{lmm}

Remark that the statement of the lemma also includes the case $b = 0$ and hence the set of parameter values to be examined is compact. The details of the proof are given in \cite{A1}. Since the non-wandering set $\mathrm{\Omega} (f_{a, b})$ is always contained in the chain recurrent set $\mathcal{R}(f_{a, b})$ of $f_{a, b}$, the above lemma yields that $f_{a, b}$ is hyperbolic on $\Omega(f_{a, b})$. Since each connected component of the parameter region where we verified hyperbolicity meets $\{a>2(1+|b|)^2\}$, we conclude that $f_{a, b}$ is a hyperbolic horseshoe on $\R^2$. This completes the proof of the claim (ii).
\end{proof}

\begin{figure}
  \includegraphics[width=7.5cm]{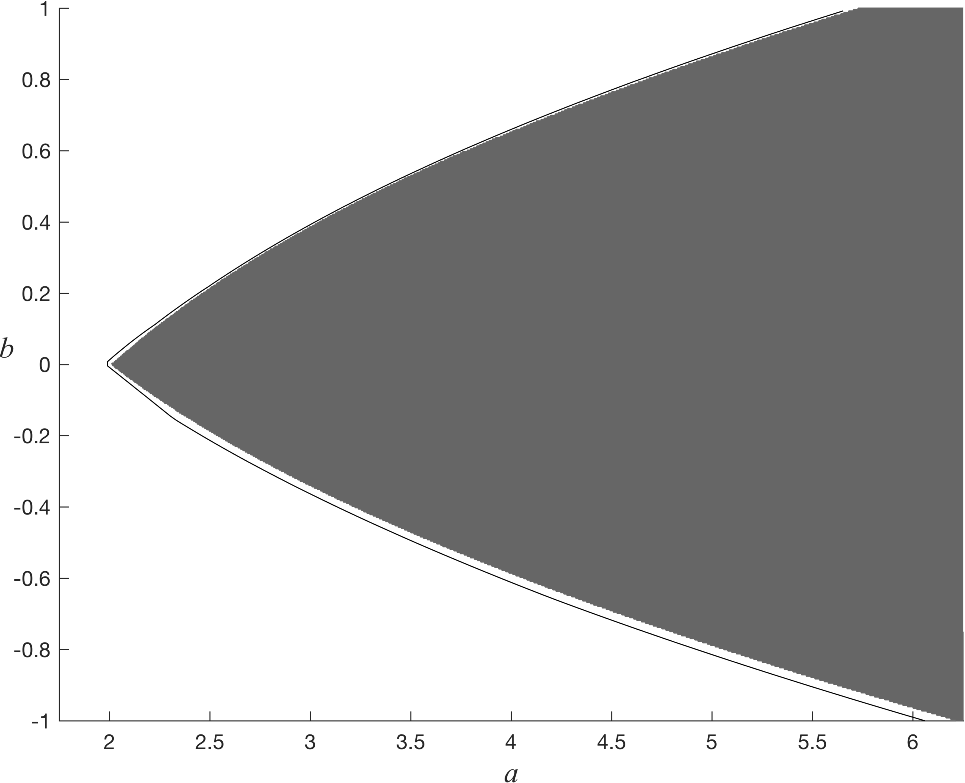}
  \hspace{0.3cm}
  \includegraphics[width=7.5cm]{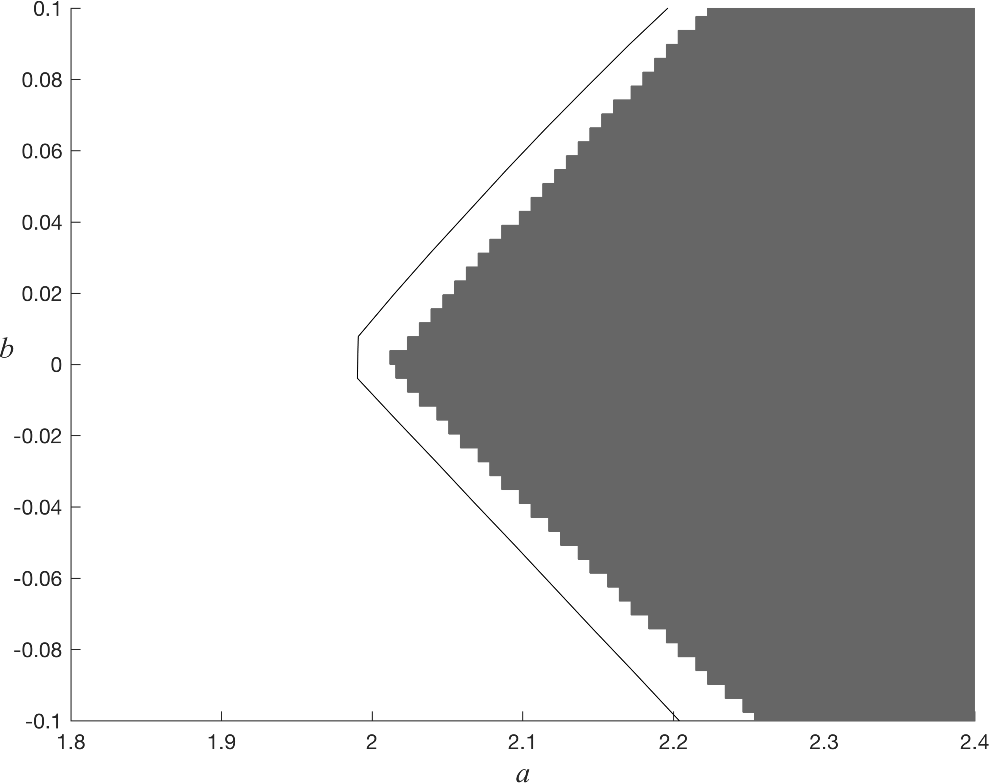}
  \caption{Left: $f_{a, b}$ is rigorously shown to be a hyperbolic horseshoe on $\R^2$ in the shaded region, and $\htop(f_{a, b}|_{\R^2}) < \log 2$ is verified in the parameter region left to the solid line. Right: a closeup view to $(a, b) = (2, 0).$}
  \label{FIG:plateau}
\end{figure}

\begin{proof}[Proof of (iii) of Theorem~\ref{THM:quasitrichotomy}]
For each $(a, b)\in \mathcal{F}^{\pm}_{\R}$ we compute the intersecting points in the trellis of $f_{a, b}$ to obtain the quadrilaterals $\mathcal{Q}^{\pm}_i$ which define projective coordinates $(\pi^{\pm}_{u, i}, \pi^{\pm}_{v, i})$ as explain in the previous subsection. Our main task here is therefore to find appropriate topological disks $D^{\pm}_{u, i}$ and $D^{\pm}_{v, i}$ so that $(f_{a, b}, \{\mathcal{B}^{\pm}_i\}_i, \mathfrak{T}^{\pm})$ satisfies the (CMC), where $\mathcal{B}^{\pm}_i=D^{\pm}_{u, i}\times_{\mathrm{pr}}D^{\pm}_{v, i}$. In Subsection~\ref{subsection6.3} we present a recipe to find appropriate topological disks $D^{\pm}_{u, i}$ and $D^{\pm}_{v, i}$. This construction gives a family of boxes $\{\mathcal{B}^{\pm}_i\}_i$ as well as a family of projective coordinates $\{(\pi^{\pm}_{u, i}, \pi^{\pm}_{v, i})\}_i$. With the help of rigorous numerics, we show 

\begin{lmm}
For every $(a, b)\in \mathcal{F}^{\pm}$, the two conditions $\pi^{\pm}_{u, j}\circ f(\partial^v\mathcal{B}^{\pm}_i)\cap D^{\pm}_{u, j}=\emptyset$ and $\pi^{\pm}_{v, i}\circ f^{-1}(\partial^h\mathcal{B}^{\pm}_j)\cap D^{\pm}_{v, i}=\emptyset$ hold for $(i, j)\in \mathfrak{T}^{\pm}$ where $f=f_{a, b}$.
\label{LMM:CMC}
\end{lmm}

The proof of this lemma is given in Subsection~\ref{subsection6.4}. This completes the proof of (iii). 
\end{proof}

Figure~\ref{FIG:plateau} illustrates the parameter region of our interest. When the parameter $(a, b)$ is in the shaded regions, we can rigorously show that the H\'enon map $f_{a, b}$ is uniformly hyperbolic on its chain recurrent set in $\R^2$ by employing the algorithm of~\cite{A1}. By the structural stability of a hyperbolic horseshoe, it follows that the shaded region is contained in the locus $\mathcal{H}_{\R}$. In Figure~\ref{FIG:plateau} there is also a solid curve close to the shaded region. When the parameter $(a, b)$ is either on the solid curve or on the left side of it, we can rigorously show that the complex H\'enon map $f_{a, b}$ possesses a periodic point in $\C^2 \setminus \R^2$, hence the topological entropy of $f_{a, b}$ on $\R^2$ is strictly less than $\log 2$ by~\cite{BLS}. We will show that the actual tangency curve $a = a_{\mathrm{tgc}}(b)$ is trapped in the narrow gap between the solid curve and the shaded region.

\newpage

\section{Dynamics and Parameter Space over $\C$}\label{section3}

Throughout this section we assume $(a, b)\in \mathcal{F}^{\pm}$ and basically consider the complex dynamics $f=f_{a, b} : \C^2\to\C^2$. 

\subsection{Admissibility}\label{subsection3.1}

Let $K=K_{a, b}$ be the \textit{filled Julia set} of $f_{a, b}$ consisting of points whose forward and backward orbits by $f_{a, b}$ are both bounded. Write $\mathcal{B}^+\equiv \bigcup_{i=0}^3 \mathcal{B}^+_i$ and $\mathcal{B}^-\equiv \bigcup_{i=0}^4 \mathcal{B}^-_i$, where $\mathcal{B}^{\pm}_i$ are the projective boxes constructed in (iii) of Theorem~\ref{THM:quasitrichotomy} (Quasi-Trichotomy). 

\begin{prp}
If $(a, b)\in \mathcal{F}^{\pm}\cap\{b\ne 0\}$, then we have $K_{a, b}=\bigcap_{n\in\Z}f_{a, b}^n(\mathcal{B}^{\pm})$.
\label{PRP:filled_Julia}
\end{prp}

Recall that we have defined 
\[R\equiv \frac{1+|b|+\sqrt{(1+|b|)^2+4|a|}}{2}\] 
and 
\[\mathcal{D}\equiv \bigl\{(x, y)\in\C^2 : |x|\leq R, |y|\leq R\bigr\}.\]
To prove Proposition~\ref{PRP:filled_Julia} we first need 

\begin{lmm} 
For any $(a, b)\in\mathcal{F}^{\pm}\cap \{b\ne 0\}$ there exists $N>0$ so that $\bigcap_{n=-N}^Nf_{a, b}^n(\mathcal{D}) \subset \mathcal{B}^{\pm}$.
\label{LMM:bounding_K}
\end{lmm}

The proof of this lemma requires computer assistance and will be given in Subsection~\ref{subsection6.4}. 

\begin{proof}[Proof of Proposition~\ref{PRP:filled_Julia}.]
One easily sees $K_{a, b}\subset \mathcal{D}$. By the $f_{a, b}$-invariance of $K_{a, b}$, this implies $\bigcap_{n\in\Z}f_{a, b}^n(\mathcal{D})=K_{a, b}$. By Lemma~\ref{LMM:bounding_K} we have $K_{a, b}\subset\mathcal{B}^{\pm}$, which yields the conclusion.
\end{proof}

Let us write $\mathrm{\Sigma}^+\equiv\{0, 1, 2, 3\}$ and $\mathrm{\Sigma}^-\equiv\{0, 1, 2, 3, 4\}$. Define
\[\mathfrak{S}^{\pm}_{\mathrm{fwd}}\equiv\bigl\{(i_n)_{n\geq 0}\in (\mathrm{\Sigma}^{\pm})^{\N} : \mbox{$(i_n, i_{n+1})\in \mathfrak{T}^{\pm}$ for all $n\geq 0$} \bigr\}\]
and call its element a \textit{forward admissible sequence} with respect to $\mathfrak{T}^{\pm}$. Also define 
\[\mathfrak{S}^{\pm}_{\mathrm{bwd}}\equiv\bigl\{(i_n)_{n\leq 0}\in (\mathrm{\Sigma}^{\pm})^{-\N} : \mbox{$(i_{n-1}, i_n)\in \mathfrak{T}^{\pm}$ for all $n\leq 0$} \bigr\}\]
and call its element a \textit{backward admissible sequence} with respect to $\mathfrak{T}^{\pm}$. Finally, we set
\[\mathfrak{S}^{\pm}\equiv\bigl\{(i_n)_{n\in \Z}\in (\mathrm{\Sigma}^{\pm})^{\Z} : \mbox{$(i_n, i_{n+1})\in \mathfrak{T}^{\pm}$ for all $n\in\Z$} \bigr\}\]
and call its element a \textit{bi-infinite admissible sequence} with respect to $\mathfrak{T}^{\pm}$. For $z\in K_{a, b}$ a symbol sequence $(i_n)_{n\geq 0}\in \mathfrak{S}^{\pm}_{\mathrm{fwd}}$ (resp. $(i_n)_{n\leq 0}\in \mathfrak{S}^{\pm}_{\mathrm{bwd}}$) satisfying $f^n(z)\in \mathcal{B}^{\pm}_{i_n}$ for $n\geq 0$ (resp. for $n\leq 0$) is called a \textit{forward itinerary} (resp. \textit{backward itinerary}) of $z$. 

The following Propositions~\ref{PRP:coding_positive} and~\ref{PRP:coding_negative} tell that the orbit of a point in $K_{a, b}$ can be traced by a sequence of appropriate crossed mappings. First consider the case $(a, b)\in\mathcal{F}^+$. 

\begin{prp}
Let $(a, b)\in \mathcal{F}^+\cap\{b\ne 0\}$. Then, for any $z\in K_{a, b}$ there exists a bi-infinite admissible sequence $(i_n)_{n\in\Z}\in \mathfrak{S}^+$ so that $f^n(z)\in \mathcal{B}^+_{i_n}$ holds for all $n\in\Z$.
\label{PRP:coding_positive} 
\end{prp}

\begin{figure}
  \includegraphics[height=4.5cm]{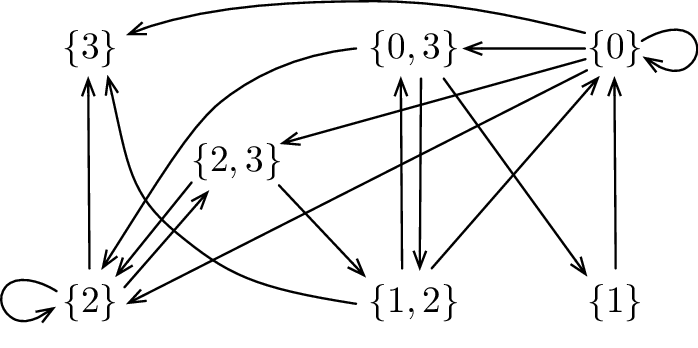}
  \caption{Diagram of allowed transitions for $(a, b)\in\mathcal{F}^+$ in Lemma~\ref{LMM:allowed_positive}.}
\label{FIG:allowed_positive}
\end{figure}

The proof of this proposition goes in the same spirit as (i) of Theorem 4.23 in~\cite{I1}. For 
\[I\in \bigl\{\{0\}, \{1\}, \{2\}, \{3\}, \{0, 1\}, \{0, 2\}, \{0, 3\}, \{1, 2\}, \{1, 3\}, \{2, 3\}\bigr\}\]
we set $\mathcal{B}^+_I\equiv \bigcap_{i\in I}\mathcal{B}^+_i$. A sequence of transitions $\cdots\to I_{n-1}\to I_n\to I_{n+1}\to\cdots$ is said to be \textit{allowed} if there exists a point $z\in\bigcap_{n\in\Z}f^n(\mathcal{B}^+)$ so that $f^n(z)\in\mathcal{B}^+_{I_n}$ holds for all $n\in\Z$. The following claims can be verified by using rigorous computation whose proof will be given in Subsection~\ref{subsection6.4}.

\begin{lmm}
Any allowed transition for $(a, b)\in\mathcal{F}^+$ is a sequence of the following $19$ arrows: $\{0\} \to \{0\}$, $\{0\} \to \{0, 3\}$, $\{0\} \to \{3\}$, $\{0\} \to \{2, 3\}$, $\{0\} \to \{2\}$, $\{0, 3\} \to \{2\}$, $\{0, 3\} \to \{1, 2\}$, $\{0, 3\} \to \{1\}$, $\{3\} \to \{1\}$, $\{2, 3\} \to \{2\}$, $\{2, 3\} \to \{1, 2\}$, $\{2, 3\} \to \{1\}$, $\{2\} \to \{3\}$, $\{2\} \to \{2, 3\}$, $\{2\} \to \{2\}$, $\{1, 2\} \to \{0\}$, $\{1, 2\} \to \{0, 3\}$, $\{1, 2\} \to \{3\}$ and $\{1\} \to \{0\}$. 
\label{LMM:allowed_positive}
\end{lmm}

Figure~\ref{FIG:allowed_positive} describes all the $19$ allowed transitions for $(a, b)\in\mathcal{F}^+$ in Lemma~\ref{LMM:allowed_positive}. The next lemma, which is essential in the proof of Proposition~\ref{PRP:coding_positive}, immediately follows from Lemma~\ref{LMM:allowed_positive}, hence its proof is omitted.

\begin{lmm}
Let $I \to I'$ be one of the $19$ arrows listed in Lemma~\ref{LMM:allowed_positive}. Then, (1) for any $i'\in I'$ there exists $i\in I$ so that $(i, i')\in\mathfrak{T}^+$ holds, and (2) for any $i\in I$ there exists $i'\in I'$ so that $(i, i')\in\mathfrak{T}^+$ holds if $\mathrm{card}(I')=2$.
\label{LMM:rule_positive}
\end{lmm}

\begin{proof}[Proof of Proposition~\ref{PRP:coding_positive}.]
Take a point $z\in K_{a, b}$. Then, there exists a unique $I_n$ so that $f^n(z)\in \mathcal{B}^+_{I_n}$ for any $n\in\Z$. We set $\mathcal{N}\equiv \{n\in\Z : \mathrm{card}(I_n)=1\}$. Assume first that $\mathcal{N}=\emptyset$. Then, the only possible allowed transition is $\cdots\to\{0, 3\}\to\{1, 2\}\to\{0, 3\}\to\{1, 2\}\to\cdots$ (see Figure~\ref{FIG:allowed_positive}). Claims (1) and (2) of Lemma~\ref{LMM:rule_positive} yield that for $n\in\Z$ there exists $i_n\in I_n$ so that $(i_n)_{n\in\Z}\in\mathfrak{S}^+$ holds. Assume next that $\mathcal{N}\ne\emptyset$ and $\sup\mathcal{N}=+\infty$. We may suppose $\inf\mathcal{N}=-\infty$ (the proof for the case $\inf\mathcal{N}>-\infty$ is similar). Let $\cdots <n_{k-1}<n_k<n_{k+1}<\cdots$ ($k\in\Z$) be the elements of $\mathcal{N}$. For any $k\in \Z$ we apply (1) of Lemma~\ref{LMM:rule_positive} to the arrow $I_{n_k-1}\to I_{n_k}$ and next to $I_{n_k-2}\to I_{n_k-1}$ until we arrive at $I_{n_{k-1}}\to I_{n_{k-1}+1}$. This determines $i_{n_{k-1}}\in I_{n_{k-1}}, \dots, i_{n_k}\in I_{n_k}$ for any $k\in\Z$, hence $(i_n)_{n\in\Z}\in\mathfrak{S}^+$. Assume finally that $\mathcal{N}\ne\emptyset$ and $N\equiv \sup\mathcal{N}<+\infty$. We can determine $i_n\in I_n$ for any $n\leq N$ as in the previous case. Note that $\mathrm{card}(I_N)=1$ and $\mathrm{card}(I_n)=2$ hold for all $n>N$. Then, the only possibilities for the transitions $I_N\to I_{N+1}\to I_{N+2}\to\cdots$ are either $\{0\}\to \{0, 3\}\to\{1, 2\}\to\{0, 3\}\to\{1, 2\}\to\cdots$, $\{0\}\to\{2, 3\}\to\{1, 2\}\to\{0, 3\}\to\{1, 2\}\to\{0, 3\}\to\cdots$ or $\{2\}\to\{2, 3\}\to\{1, 2\}\to\{0, 3\}\to\{1, 2\}\to\{0, 3\}\to\cdots$ (see Figure~\ref{FIG:allowed_positive} again). In each of these three cases we can successively apply (2) of Lemma~\ref{LMM:rule_positive} to determine $i_n$ for $n>N$. Hence $(i_n)_{n\in\Z}\in\mathfrak{S}^+$, and this proves Proposition~\ref{PRP:coding_positive}.
\end{proof}

Next consider the case $(a, b)\in\mathcal{F}^-$. 

\begin{prp}
Let $(a, b)\in \mathcal{F}^-\cap\{b\ne 0\}$. Then, for any $z\in K_{a, b}$ there exists a bi-infinite admissible sequence $(i_n)_{n\in\Z}\in \mathfrak{S}^-$ so that $f^n(z)\in \mathcal{B}^-_{i_n}$ holds for all $n\in\Z$. 
\label{PRP:coding_negative}
\end{prp}

\begin{figure}
 \includegraphics[height=4.5cm]{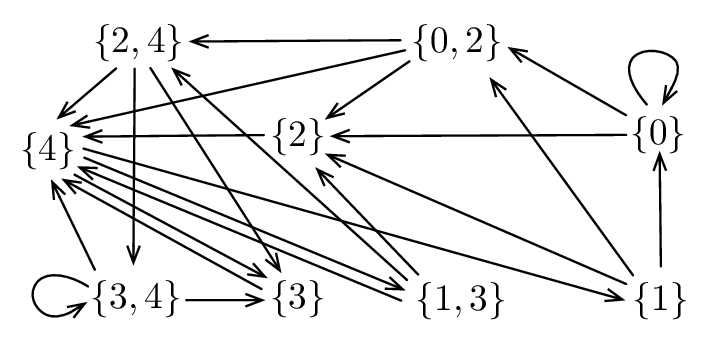}
  \caption{Diagram of allowed transitions for $(a, b)\in\mathcal{F}^-$ in Lemma~\ref{LMM:allowed_negative}.}
\label{FIG:allowed_negative}
\end{figure}

For 
\[I\in \bigl\{\{0\}, \{1\}, \{2\}, \{3\}, \{4\}, \{0, 1\}, \{0, 2\}, \{0, 3\}, \{0, 4\}, \{1, 2\}, \{1, 3\}, \{1, 4\}, \{2, 3\}, \{2, 4\}, \{3, 4\}\bigr\}\]
we set $\mathcal{B}^-_I\equiv \bigcap_{i\in I}\mathcal{B}^-_i$. A sequence of transitions $\cdots\to I_{n-1}\to I_n\to I_{n+1}\to\cdots$ is said to be \textit{allowed} if there exists a point $z\in\bigcap_{n\in\Z}f^n(\mathcal{B}^-)$ so that $f^n(z)\in\mathcal{B}^-_{I_n}$ holds for all $n\in\Z$. The following claims can be verified by using rigorous computation whose proof will be given in Subsection~\ref{subsection6.4}.

\begin{lmm}
Any allowed transition for $(a, b)\in\mathcal{F}^-$ is a sequence of the following $23$ arrows: $\{0\}\to\{0\}$, $\{0\}\to\{0, 2\}$, $\{0\}\to\{2\}$, $\{0, 2\}\to\{2\}$, $\{0, 2\}\to\{2, 4\}$, $\{0, 2\}\to\{4\}$, $\{2\}\to\{4\}$, $\{2, 4\}\to\{3\}$, $\{2, 4\}\to\{3, 4\}$, $\{2, 4\}\to\{4\}$, $\{4\}\to\{1\}$, $\{4\}\to\{1, 3\}$, $\{4\}\to\{3\}$, $\{3, 4\}\to\{3\}$, $\{3, 4\}\to\{3, 4\}$, $\{3, 4\}\to\{4\}$, $\{3\}\to\{4\}$, $\{1, 3\}\to\{2\}$, $\{1, 3\}\to\{2, 4\}$, $\{1, 3\}\to\{4\}$, $\{1\}\to\{0\}$, $\{1\}\to\{0, 2\}$ and $\{1\}\to\{2\}$. 
\label{LMM:allowed_negative}
\end{lmm}

Figure~\ref{FIG:allowed_negative} describes all the $23$ allowed transitions for $(a, b)\in\mathcal{F}^-$ in Lemma~\ref{LMM:allowed_negative}. The next lemma, which is essential in the proof of Proposition~\ref{PRP:coding_negative}, immediately follows from Lemma~\ref{LMM:allowed_negative}, hence its proof is omitted.

\begin{lmm}
Let $I \to I'$ be one of the $23$ arrows listed in Lemma~\ref{LMM:allowed_negative}. Then, (1) for any $i'\in I'$ there exists $i\in I$ so that $(i, i')\in\mathfrak{T}^-$ holds, and (2) for any $i\in I$ there exists $i'\in I'$ so that $(i, i')\in\mathfrak{T}^-$ holds if $\mathrm{card}(I')=2$.
\label{LMM:rule_negative}
\end{lmm}

\begin{proof}[Proof of Proposition~\ref{PRP:coding_negative}.]
Take a point $z\in K_{a, b}$. Then, there exists a unique $I_n$ so that $f^n(z)\in \mathcal{B}^-_{I_n}$ for any $n\in\Z$. We set $\mathcal{N}\equiv \{n\in\Z : \mathrm{card}(I_n)=1\}$. Assume first that $\mathcal{N}=\emptyset$. Then, the only possible allowed transition is $\cdots\to\{3, 4\}\to\{3, 4\}\to\{3, 4\}\to\cdots$ (see Figure~\ref{FIG:allowed_negative}). Claims (1) and (2) of Lemma~\ref{LMM:rule_negative} yield that for $n\in\Z$ there exists $i_n\in I_n$ so that $(i_n)_{n\in\Z}\in\mathfrak{S}^-$ holds. Assume next that $\mathcal{N}\ne\emptyset$ and $\sup\mathcal{N}=+\infty$. We may suppose $\inf\mathcal{N}=-\infty$ (the proof for the case $\inf\mathcal{N}>-\infty$ is similar). Let $\cdots <n_{k-1}<n_k<n_{k+1}<\cdots$ ($k\in\Z$) be the elements of $\mathcal{N}$. For any $k\in \Z$ we apply (1) of Lemma~\ref{LMM:rule_negative} to the arrow $I_{n_k-1}\to I_{n_k}$ and next to $I_{n_k-2}\to I_{n_k-1}$ until we arrive at $I_{n_{k-1}}\to I_{n_{k-1}+1}$. This determines $i_{n_{k-1}}\in I_{n_{k-1}}, \dots, i_{n_k}\in I_{n_k}$ for any $k\in\Z$, hence we have $(i_n)_{n\in\Z}\in\mathfrak{S}^-$. Assume finally that $\mathcal{N}\ne\emptyset$ and $N\equiv \sup\mathcal{N}<+\infty$. We can determine $i_n\in I_n$ for any $n\leq N$ as in the previous case. Note that $\mathrm{card}(I_N)=1$ and $\mathrm{card}(I_n)=2$ hold for all $n>N$. Then, the only possibilities for the transitions $I_N\to I_{N+1}\to I_{N+2}\to\cdots$ are either $\{0\}\to \{0, 2\}\to\{2, 4\}\to\{3, 4\}\to\{3, 4\}\to\{3, 4\}\to\cdots$, $\{1\}\to\{0, 2\}\to\{2, 4\}\to\{3, 4\}\to\{3, 4\}\to\{3, 4\}\to\cdots$ or $\{4\}\to\{1, 3\}\to\{2, 4\}\to\{3, 4\}\to\{3, 4\}\to\{3, 4\}\to\cdots$ (see Figure~\ref{FIG:allowed_negative} again). In each of these three cases we can successively apply (2) of Lemma~\ref{LMM:rule_negative} to determine $i_n$ for $n>N$. Hence $(i_n)_{n\in\Z}\in\mathfrak{S}^-$, and this proves Proposition~\ref{PRP:coding_negative}.
\end{proof}

\begin{figure}
\includegraphics[height=6cm]{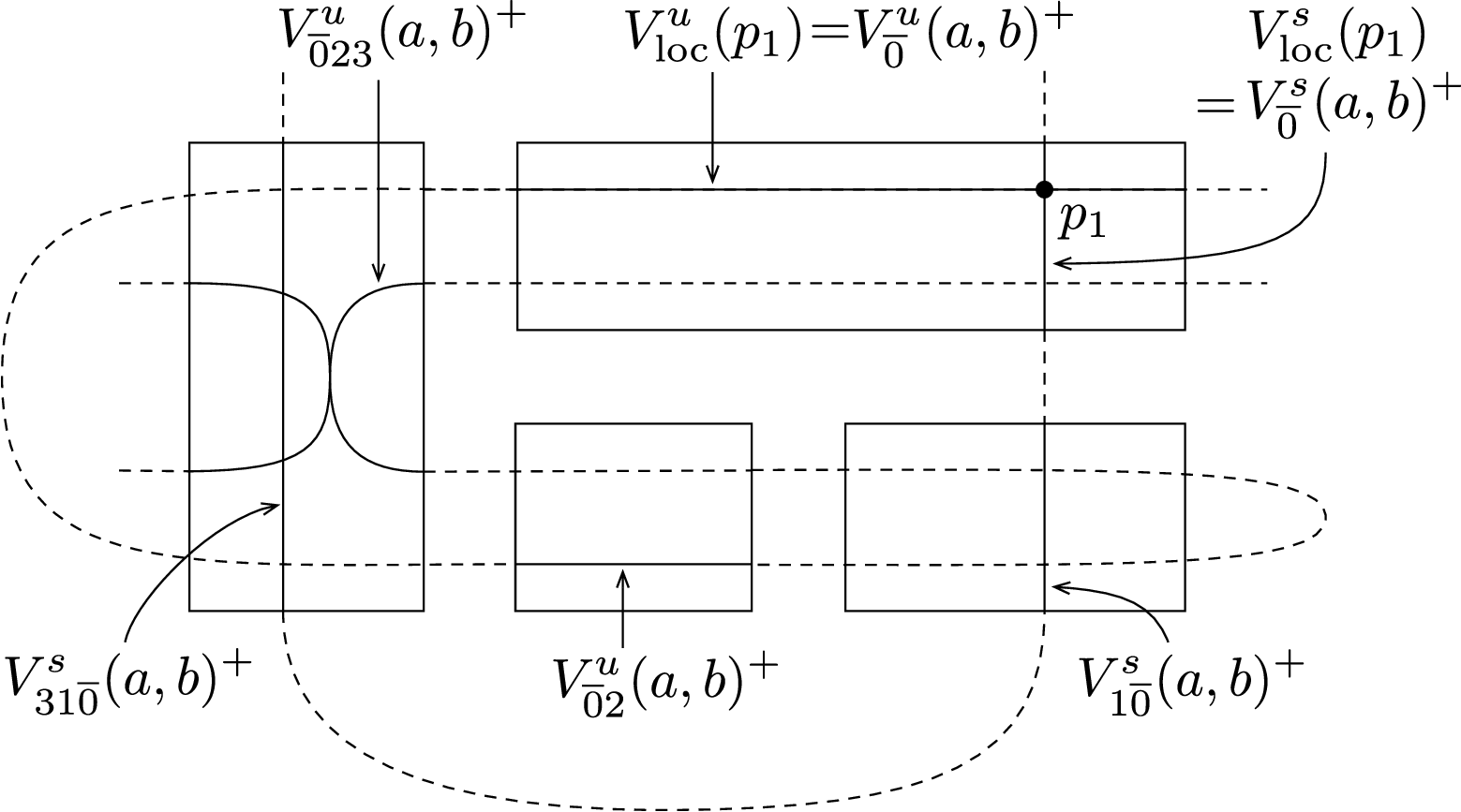}
\caption{Decomposition of invariant manifolds for $(a, b)\in\mathcal{F}^+$.}
\label{FIG:decomposition_positive}
\end{figure}

\begin{figure}
\includegraphics[height=7cm]{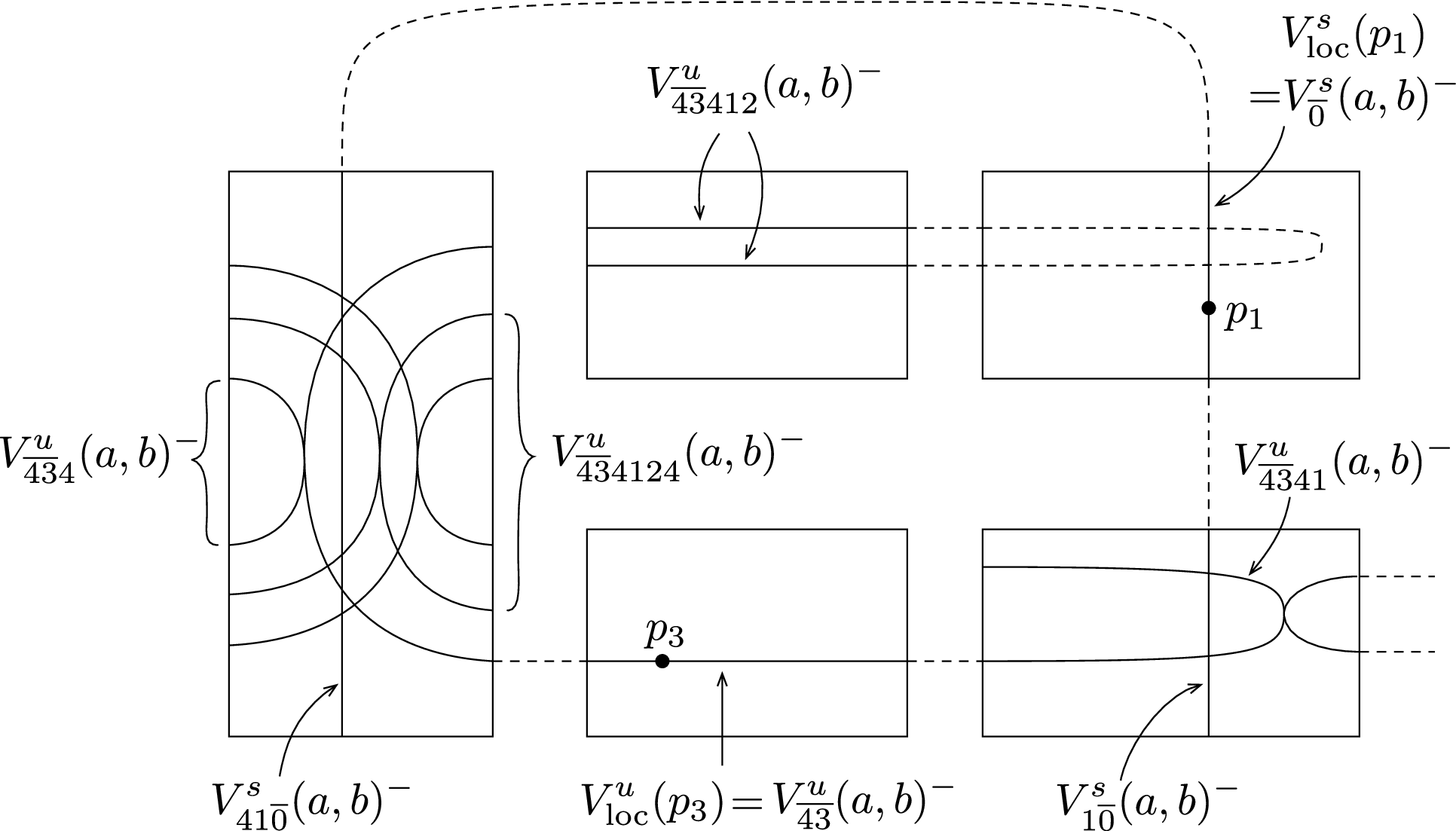}
\caption{Decomposition of invariant manifolds for $(a, b)\in\mathcal{F}^-$.}
\label{FIG:decomposition_negative}
\end{figure}

\subsection{Encoding in $\C^2$}\label{subsection3.2}

In this subsection we decompose the complex stable/unstable manifolds of some saddle points in $\C^2$ according to the projective boxes $\{\mathcal{B}^{\pm}_i\}_i$ found in Theorem~\ref{THM:quasitrichotomy} (Quasi-Trichotomy).

For $(a, b)\in (\mathcal{F}^+\cup\mathcal{F}^-)\cap\{b\ne 0\}$, let $V^{u/s}(p_i)$ be the complex unstable/stable manifolds of $p_i\in\C^2$ for the map $f_{a, b} : \C^2\to \C^2$. For $(a, b)\in (\mathcal{F}^+\cup\mathcal{F}^-)\cap \{b=0\}$, we let $V^u(p_i)\equiv \{(x, y)\in\C^2 : x=y^2-a\}$ and $V^s(p_i)\equiv \{(x, y)\in\C^2 : x=x_i\}$, where $p_i=(x_i, y_i)$. 

For $(a, b)\in\mathcal{F}^+$, let $V_{\mathrm{loc}}^s(p_1)$ be the connected component of $V^s(p_1)\cap\mathcal{B}^+_0$ containing $p_1$ and $V_{\mathrm{loc}}^u(p_1)$ be the connected component of $V^u(p_1)\cap\mathcal{B}^+_0$ containing $p_1$. Since $f : \mathcal{B}^+_0\cap f^{-1}(\mathcal{B}^+_0)\to\mathcal{B}^+_0$ is a crossed mapping of degree one, $V_{\mathrm{loc}}^s(p_1)$ is a vertical submanifold of degree one in $\mathcal{B}^+_0$ and $V_{\mathrm{loc}}^u(p_1)$ is a horizontal submanifold of degree one in $\mathcal{B}^+_0$ (see Figure~\ref{FIG:decomposition_positive}). For $(a, b)\in\mathcal{F}^-$, let $V_{\mathrm{loc}}^s(p_1)$ be the connected component of $V^s(p_1)\cap\mathcal{B}^-_0$ containing $p_1$ and $V_{\mathrm{loc}}^u(p_3)$ be the connected component of $V^u(p_3)\cap\mathcal{B}^-_1$ containing $p_3$. Since $f : \mathcal{B}^-_0\cap f^{-1}(\mathcal{B}^-_0)\to\mathcal{B}^-_0$ is a crossed mapping of degree one, $V_{\mathrm{loc}}^s(p_1)$ is a vertical submanifold of degree one in $\mathcal{B}^-_0$ (see Figure~\ref{FIG:decomposition_negative}). 

\begin{figure}
  \includegraphics[height=5.5cm]{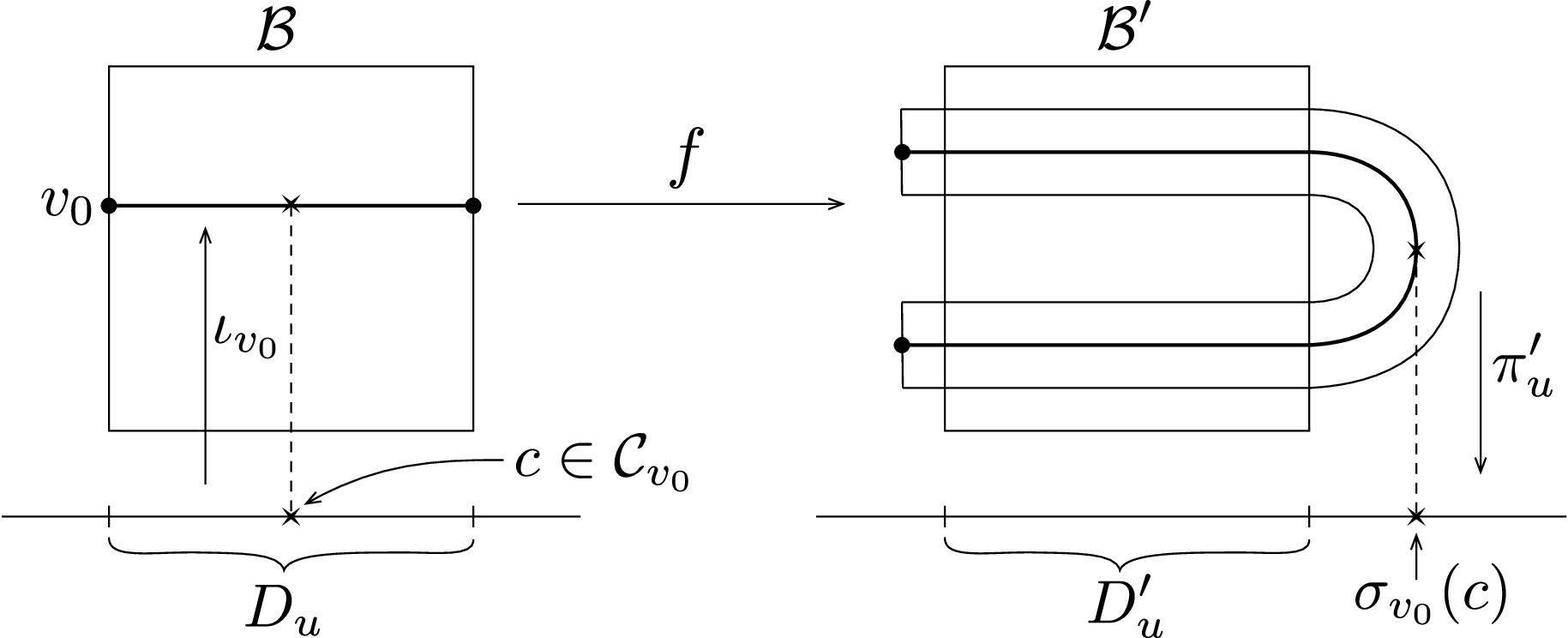}
  \caption{Figure of the off-criticality condition.}
\label{FIG:OCC}
\end{figure}

Characterizing $V_{\mathrm{loc}}^u(p_3)$ for $(a, b)\in\mathcal{F}^-$ in terms of the boxes is problematic. For this, let us recall the following notion from~\cite{ISm}. Let $\mathcal{B}=D_u\times_{\mathrm{pr}}D_v$ and $\mathcal{B}'=D'_u\times_{\mathrm{pr}}D'_v$ be two projective bidisks and let $f : \C^2\to \C^2$ be a complex H\'enon map satisfying the boundary compatibility condition with respect to $\mathcal{B}$ and $\mathcal{B}'$. For each $v_0\in D_v$, define 
\[\sigma_{v_0}\equiv\pi'_u\circ f\circ \iota_{v_0} : D_u\longrightarrow L'_u,\]
where $\iota_{v_0} : D_u\to \mathcal{B}$ is given by $u\mapsto (u, v_0)$ in the projective coordinates of $\mathcal{B}$.

\begin{dfn}
We say that $f : \C^2\to \C^2$ satisfies the \textit{off-criticality condition (OCC)} with respect to $\mathcal{B}$ and $\mathcal{B}'$ if $\sigma_{v_0}(\mathcal{C}_{v_0})\cap D'_u=\emptyset$ holds for every $v_0\in D_v$, where $\mathcal{C}_{v_0}$ denotes the critical points $c$ of $\sigma_{v_0}$ (see Figure~\ref{FIG:OCC}). 
\label{DFN:OCC}
\end{dfn}

\begin{figure}
\includegraphics[width=11cm]{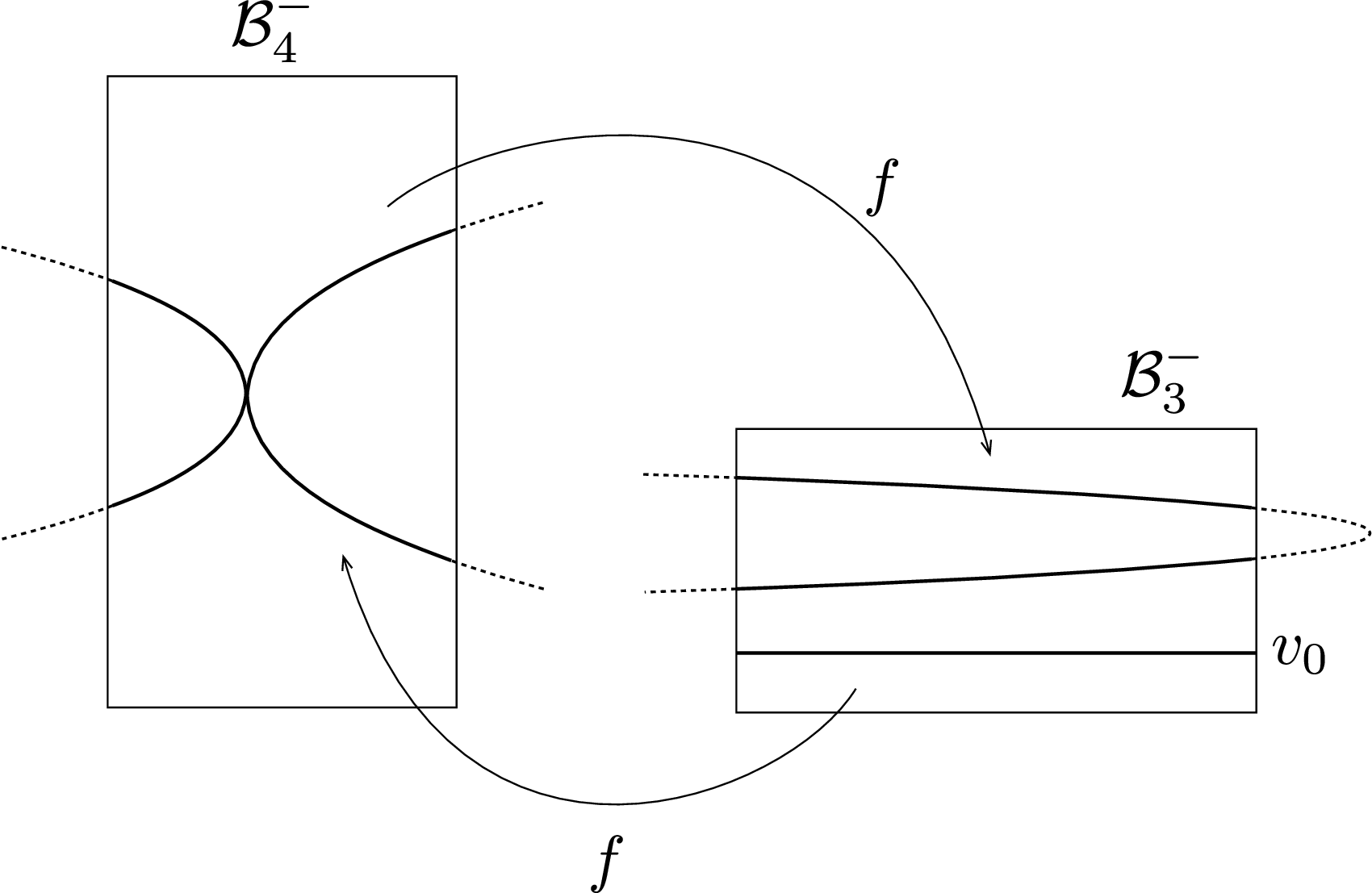}
\caption{Figure of Lemma~\ref{LMM:OCC_negative1}.}
\label{FIG:OCC_negative1}
\end{figure}

With this notion we prove the next claim.

\begin{prp}
For $(a, b)\in \mathcal{F}^-$, $V_{\mathrm{loc}}^u(p_3)$ is a horizontal submanifold of degree one in $\mathcal{B}^-_3$.
\label{PRP:local_unstable}
\end{prp}

To prove this proposition we need

\begin{lmm}
Let $(a, b)\in\mathcal{F}^-$. Then, for every fixed $v_0\in D_{v, 3}^-$ we have 
\[\frac{d}{du}\left\{\pi^-_{u, 3} \circ f^2 \circ\iota_{v_0}(u)\right\}\ne 0\]
for $u\in D_{u, 3}^-$ with $\iota_{v_0}(u)\in \mathcal{B}^-_3\cap f^{-1}(\mathcal{B}^-_4\cap f^{-1}(\mathcal{B}^-_3))$ (see Figure~\ref{FIG:OCC_negative1}).
\label{LMM:OCC_negative1}
\end{lmm}

The proof of this lemma requires computer assistance and will be supplied in Subsection~\ref{subsection6.4}.

\begin{proof}[Proof of Proposition~\ref{PRP:local_unstable}.]
From Lemma~\ref{LMM:OCC_negative1} it follows that $f^2 : \mathcal{B}^-_3\cap f^{-1}(\mathcal{B}^-_4\cap f^{-1}(\mathcal{B}^-_3))\to\mathcal{B}^-_3$ is a crossed mapping of degree two satisfying the (OCC), hence is of \textit{horseshoe type}, that is, $\mathcal{B}^-_3\cap f^{-1}(\mathcal{B}^-_4\cap f^{-1}(\mathcal{B}^-_3))$ has two connected components and the restriction of $f^2$ to each component is of degree one.

Take a horizontal submanifold $D_0$ of degree one in $\mathcal{B}^-_3$ through $p_3$. When $b\ne 0$ (resp. $b=0$), $\mathcal{B}^-_3\cap f(\mathcal{B}^-_4\cap f(D_0))$ consists of two horizontal submanifolds (resp. one horizontal submanifold) of degree one in $\mathcal{B}^-_3$ by the discussion above. Choose the one containing the fixed point $p_3$ and call it $D_1$. We repeat this procedure to obtain a sequence of horizontal submanifolds $D_n$ of degree one in $\mathcal{B}^-_3$. By the Lambda Lemma, $D_n$ converges to $V_{\mathrm{loc}}^u(p_3)$ in the Hausdorff topology, hence $V_{\mathrm{loc}}^u(p_3)$ is a horizontal submanifold of degree one in $\mathcal{B}^-_3$.
\end{proof}

Let us decompose complex stable/unstable manifolds $V^{u/s}(p_i)$ into several pieces according to the family of boxes $\{\mathcal{B}^{\pm}_i\}_i$. Below, $\overline{0}$ means either $\cdots 00$ or $00\cdots$. For a forward admissible sequence of the form $I=i_0 i_1\cdots i_n\overline{0}\in\mathfrak{S}_{\mathrm{fwd}}^+$ we define
\[V^s_I(a, b)^+\equiv \mathcal{B}^+_{i_0}\cap f_{a, b}^{-1}(\mathcal{B}^+_{i_1}\cap\cdots\cap f_{a, b}^{-1}(\mathcal{B}^+_{i_n}\cap f_{a, b}^{-1}(V^s_{\mathrm{loc}}(p_1)))\cdots),\]
and for a backward admissible sequence of the form $J=\overline{0}j_{-n}\cdots j_{-1}j_0\in\mathfrak{S}_{\mathrm{bwd}}^+$ we define
\[V^u_J(a, b)^+\equiv \mathcal{B}^+_{j_0}\cap f_{a, b}(\mathcal{B}^+_{j_{-1}}\cap\cdots\cap f_{a, b}(\mathcal{B}^+_{j_{-n}}\cap f_{a, b}(V^u_{\mathrm{loc}}(p_1)))\cdots).\]
Among these pieces we are particularly interested in 
\[V_{31\overline{0}}^s(a, b)^+\equiv \mathcal{B}^+_3\cap f_{a, b}^{-1}(\mathcal{B}^+_1\cap f_{a, b}^{-1}(V^s_{\mathrm{loc}}(p_1)))\]
which is a degree one vertical submanifold in $\mathcal{B}^+_3$, and
\[V_{\overline{0}23}^u(a, b)^+\equiv \mathcal{B}^+_3\cap f_{a, b}(\mathcal{B}^+_2\cap f_{a, b}(V^u_{\mathrm{loc}}(p_1)))\]
which is a degree two horizontal submanifold in $\mathcal{B}^+_3$. 

Let $(a, b)\in \mathcal{F}^-$. Below, $\overline{43}$ means $\cdots 4343$. For a forward admissible sequence of the form $I=i_0 i_1\cdots i_n\overline{0}\in\mathfrak{S}_{\mathrm{fwd}}^-$ we define
\[V^s_I(a, b)^-\equiv \mathcal{B}^-_{i_0}\cap f_{a, b}^{-1}(\mathcal{B}^-_{i_1}\cap\cdots\cap f_{a, b}^{-1}(\mathcal{B}^-_{i_n}\cap f_{a, b}^{-1}(V^s_{\mathrm{loc}}(p_1)))\cdots),\]
and for a backward admissible sequence of the form $J=\overline{43}j_{-n}\cdots j_{-1}j_0\in\mathfrak{S}_{\mathrm{bwd}}^-$ we define
\[V^u_J(a, b)^-\equiv \mathcal{B}^-_{j_0}\cap f_{a, b}(\mathcal{B}^-_{j_{-1}}\cap\cdots\cap f_{a, b}(\mathcal{B}^-_{j_{-n}}\cap f_{a, b}(V^u_{\mathrm{loc}}(p_3)))\cdots).\]
Among these pieces we are particularly interested in 
\[V_{41\overline{0}}^s(a, b)^-\equiv \mathcal{B}^-_4\cap f_{a, b}^{-1}(\mathcal{B}^-_1\cap f_{a, b}^{-1}(V^s_{\mathrm{loc}}(p_1)))\]
which is a degree one vertical submanifold in $\mathcal{B}^-_2$, and finally we define
\[V_{\overline{43}4124}^u(a, b)^-\equiv \mathcal{B}^-_4\cap f_{a, b}(\mathcal{B}^-_2\cap f_{a, b}(\mathcal{B}^-_1\cap f_{a, b}(\mathcal{B}^-_4\cap f_{a, b}(V^u_{\mathrm{loc}}(p_3))))).\]

The above submanifolds $V_{31\overline{0}}^s(a, b)^+$, $V_{\overline{0}23}^u(a, b)^+$, $V_{41\overline{0}}^s(a, b)^-$ and $V_{\overline{43}4124}^u(a, b)^-$ are called the \textit{special pieces} and will play a important role in what follows. Note that these submanifolds are well-defined even for the case $b=0$. To deal with the last one, it is useful to consider
\[V^u_{\overline{43}412}(a, 0)^-\equiv\mathcal{B}^-_2\cap f_{a, 0}(\mathcal{B}^-_1\cap f_{a, 0}(\mathcal{B}^-_4\cap f_{a, 0}(V_{\mathrm{loc}}^u(p_3)))).\]

\begin{figure}
\includegraphics[width=12cm]{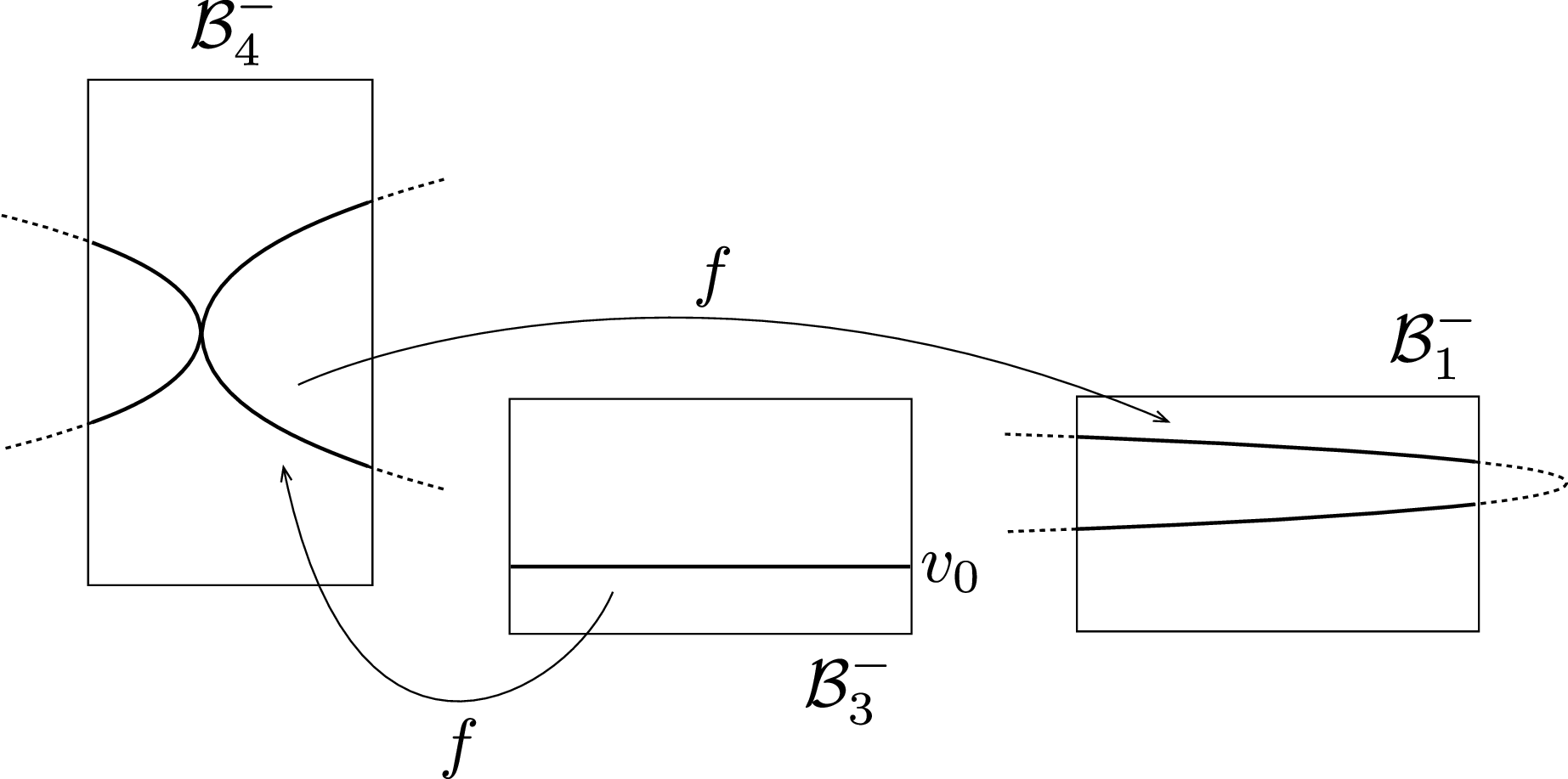}
\caption{Figure of Lemma~\ref{LMM:OCC_negative2} (i).} 
\label{FIG:OCC_negative2}
\end{figure}

\begin{figure}
\includegraphics[width=12cm]{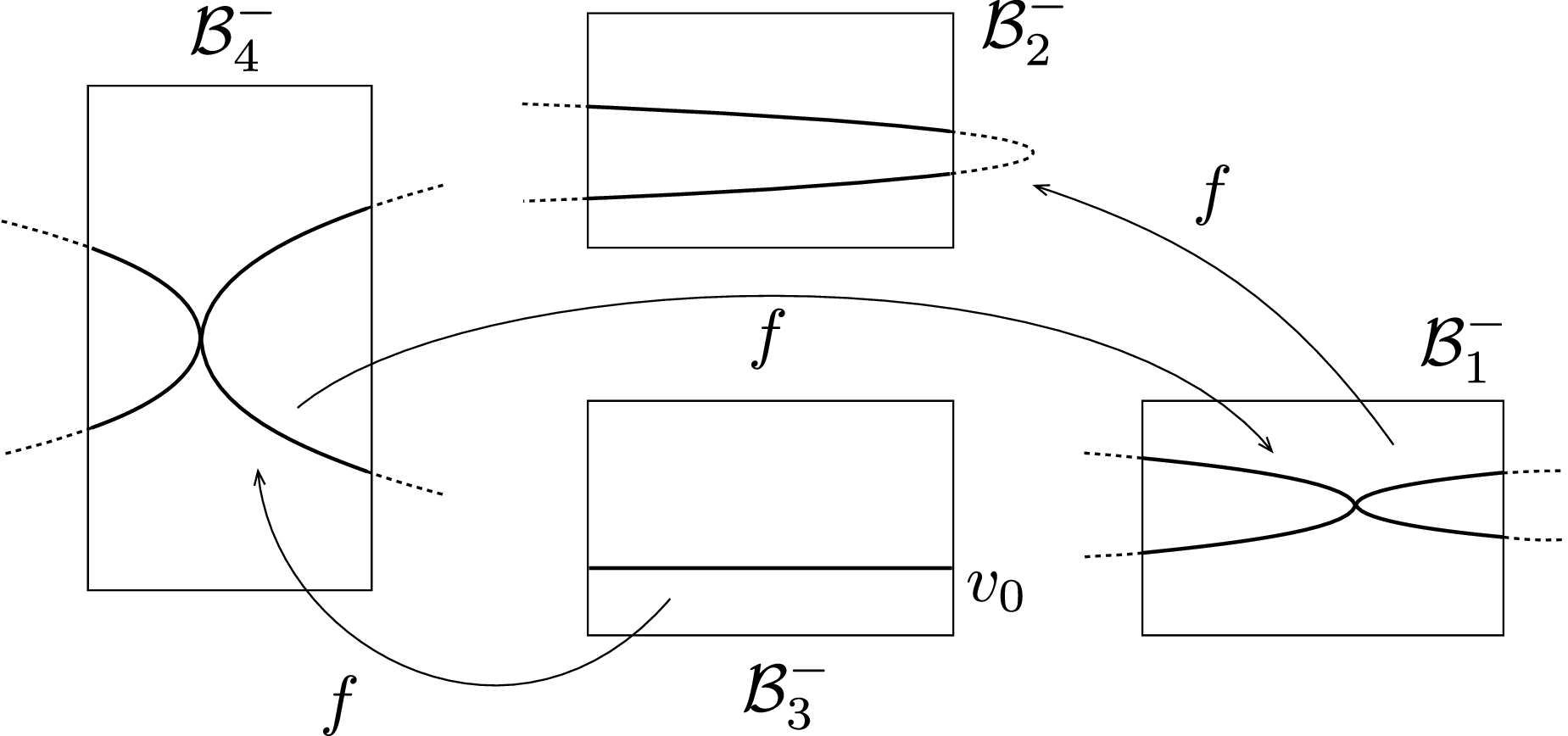}
\caption{Figure of Lemma~\ref{LMM:OCC_negative2} (ii).}
\label{FIG:OCC_negative3}
\end{figure}

\begin{prp}
When $(a, b)\in \mathcal{F}^-\cap\{b\ne 0\}$, $V^u_{\overline{43}412}(a, b)^-$ consists of two horizontal submanifolds of degree one in $\mathcal{B}^-_2$. When $(a, 0)\in \mathcal{F}^-\cap\{b=0\}$, $V^u_{\overline{43}412}(a, 0)^-$ consists of one horizontal submanifold of degree one in $\mathcal{B}^-_2$.
\label{PRP:43412}
\end{prp}

To prove this proposition we need

\begin{lmm}
Let $(a, b)\in\mathcal{F}^-$. Then, one of the following (i) and (ii) holds; 
\begin{enumerate}
\renewcommand{\labelenumi}{(\roman{enumi})}
\item for every fixed $v_0\in D_{v, 3}^-$ we have 
\[\frac{d}{du}\left\{\pi^-_{u, 1} \circ f^2\circ\iota_{v_0}(u)\right\}\ne 0\]
for $u\in D_{u, 3}^-$ with $\iota_{v_0}(u)\in \mathcal{B}^-_3\cap f^{-1}(\mathcal{B}^-_4\cap f^{-1}(\mathcal{B}^-_1))$ (see Figure~\ref{FIG:OCC_negative2}),
\item for every fixed $v_0\in D_{v, 3}^-$ we have 
\[\frac{d}{du}\left\{\pi^-_{u, 2} \circ f^3\circ\iota_{v_0}(u)\right\}\ne 0\]
for $u\in D_{u, 3}^-$ with $\iota_{v_0}(u)\in \mathcal{B}^-_3\cap f^{-1}(\mathcal{B}^-_4\cap f^{-1}(\mathcal{B}^-_1\cap f^{-1}(\mathcal{B}^-_2)))$ (see Figure~\ref{FIG:OCC_negative3}).
\end{enumerate}
\label{LMM:OCC_negative2}
\end{lmm}

The proof of this lemma requires computer assistance and will be supplied in Subsection~\ref{subsection6.4}.

\begin{proof}[Proof of Proposition~\ref{PRP:43412}.]
Since $f : \mathcal{B}^-_1\cap f^{-1}(\mathcal{B}^-_2)\to\mathcal{B}^-_2$ is a crossed mapping of degree one, the case (i) yields that $f^3 : \mathcal{B}^-_3\cap f^{-1}(\mathcal{B}^-_4\cap f^{-1}(\mathcal{B}^-_1\cap f^{-1}(\mathcal{B}^-_2)))\to\mathcal{B}^-_2$ is a crossed mapping of degree two satisfying the (OCC), hence is of horseshoe type. In the case of (ii) we immediately obtain the same conclusion. Hence, in both cases we obtain Proposition~\ref{PRP:43412}. 
\end{proof}

\begin{figure}
\includegraphics[height=12cm]{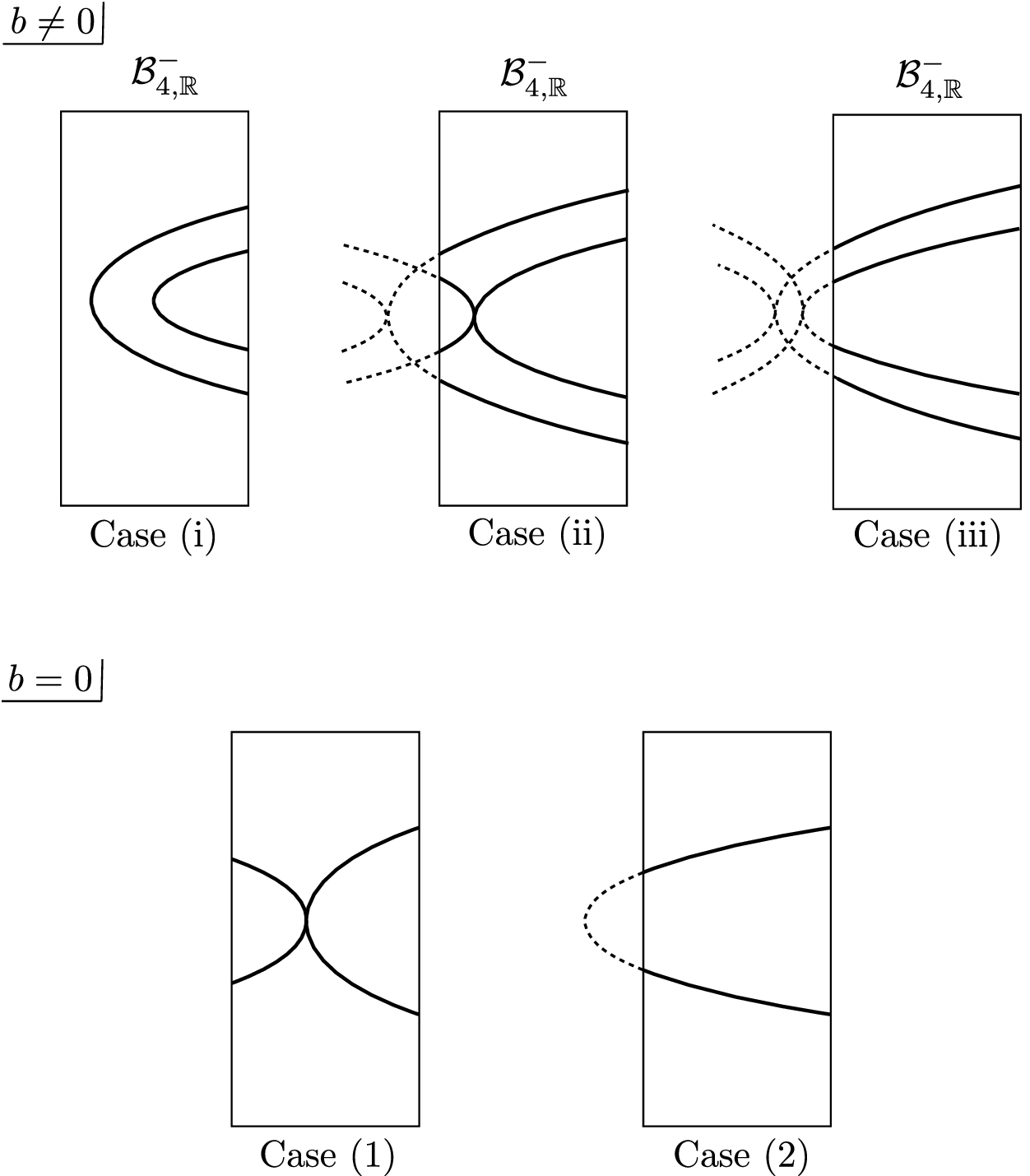}
\caption{Figures of $V^u_{\overline{43}4124}$: (a) case $b\ne 0$, (b) case $b=0$.}
\label{FIG:434124}
\end{figure}

In particular, when $(a, b)\in \mathcal{F}^-\cap\{b\ne 0\}$, the special piece $V_{\overline{43}4124}^u(a, b)^-$ consists of either (i) two mutually disjoint horizontal submanifolds of degree two in $\mathcal{B}^-_4$, (ii) one horizontal submanifold of degree two and two horizontal submanifolds of degree one in $\mathcal{B}^-_4$ all mutually disjoint, or (iii) four mutually disjoint horizontal submanifolds of degree one in $\mathcal{B}^-_4$ (see the top of Figure~\ref{FIG:434124}). When $(a, 0)\in \mathcal{F}^-\cap\{b=0\}$, the special piece $V_{\overline{43}4124}^u(a, 0)^-$ consists of either (1) a single horizontal submanifold of degree two in $\mathcal{B}^-_4$ or (2) two mutually disjoint horizontal submanifolds of degree one in $\mathcal{B}^-_4$ (see the bottom of Figure~\ref{FIG:434124}).

\newpage

\section{Dynamics and Parameter Space over $\R$}\label{section4}

Throughout this section let us assume $(a, b)\in \mathcal{F}^{\pm}_{\R}$ and consider the real dynamics $f_{a, b}|_{\R^2} : \R^2\to\R^2$. Below, we use the notation $f_{\R}\equiv f_{a, b}|_{\R^2}$ and $\mathcal{B}^{\pm}_{i, \R}\equiv\mathcal{B}^{\pm}_i\cap \R^2$. Then, the invariant manifolds of $f_{\R}$ in $\R^2$ are decomposed into several pieces according to the symbolic dynamics given by the family of real boxes $\{\mathcal{B}^{\pm}_{i, \R}\}_i$. The purpose of this section is to investigate the configuration of these pieces in each box by using the crossed mapping condition proved in Theorem~\ref{THM:quasitrichotomy} (Quasi-Trichotomy) and certain plane topology arguments.

\subsection{Encoding in $\R^2$}\label{subsection4.1}

Since each box $\mathcal{B}^{\pm}_i$ moves continuously with respect to the parameters and since $\mathcal{F}^{\pm}_{\R}$ is connected and simply connected, the notions of \textit{upper boundary}, \textit{lower boundary}, \textit{right boundary} and \textit{left boundary} of the real box $\mathcal{B}^{\pm}_{i, \R}$ are well-defined as a continuation from the case $b=0$ where these definitions are obvious. 

\begin{dfn}
A curve in $\mathcal{B}^{\pm}_{i, \R}$ is said to be \textit{horizontal} (resp. \textit{vertical}) if it is a curve between the right and the left (resp. upper and lower) boundaries of $\mathcal{B}^{\pm}_{i, \R}$. We say such a curve is \textit{of degree one} if its horizontal/vertical projection is bijective.
\label{DFN:horizontal_curve}
\end{dfn}

Let $\tau : \C^2\to\C^2$ be the involution in $\C^2$ given by $\tau(x, y)\equiv(\overline{x}, \overline{y})$. A horizontal/vertical disk $D$ in a certain box is said to be \textit{real} if $\tau(D)=D$ holds. 

\begin{lmm}
If $D$ is a horizontal/vertical disk in a box $\mathcal{B}^{\pm}_i$ which is real, then the real section $D\cap\mathcal{B}^{\pm}_{i, \R}$ consists of a nonempty, connected one-dimensional curve.
\label{LMM:real_disk}
\end{lmm}

\begin{proof}
See Proposition 3.1 of~\cite{BS2}. 
\end{proof}

Examples of real disks of degree one are local invariant manifolds $V^u_{\mathrm{loc}}(p_1)$ in $\mathcal{B}_0^+$ and $V^s_{\mathrm{loc}}(p_1)$ in $\mathcal{B}_0^+$ for $(a, b)\in\mathcal{F}^+_{\R}$, and $V^u_{\mathrm{loc}}(p_1)$ in $\mathcal{B}_3^-$ and $V^s_{\mathrm{loc}}(p_1)$ in $\mathcal{B}_0^-$ for $(a, b)\in\mathcal{F}^-_{\R}$. The real sections $W^{u/s}_{\mathrm{loc}}(p_i)\equiv V^{u/s}_{\mathrm{loc}}(p_i)\cap\R^2$ are all corresponding to the local invariant manifolds at $p_i$ for the real dynamics $f_{\R} : \R^2\to\R^2$. It follows that $W^u_{\mathrm{loc}}(p_0)$ is a horizontal curve of degree one in $\mathcal{B}_{0, \R}^+$ and $W^s_{\mathrm{loc}}(p_0)$ is a vertical curve of degree one in $\mathcal{B}_{0, \R}^+$ for $(a, b)\in\mathcal{F}^+_{\R}$,  $W^u_{\mathrm{loc}}(p_3)$ is a horizontal curve of degree one in $\mathcal{B}_{3, \R}^-$ and $W^s_{\mathrm{loc}}(p_0)$ is a vertical curve of degree one in $\mathcal{B}_{0, \R}^-$ for $(a, b)\in\mathcal{F}^-_{\R}$.

Let $(a, b)\in \mathcal{F}^+_{\R}$. For a forward admissible sequence of the form $I=i_0 i_1\cdots i_n\overline{0}\in\mathfrak{S}_{\mathrm{fwd}}^+$ we define
\[W^s_I(a, b)^+\equiv \mathcal{B}^+_{i_0, \R}\cap f_{\R}^{-1}(\mathcal{B}^+_{i_1, \R}\cap\cdots\cap f_{\R}^{-1}(\mathcal{B}^+_{i_n, \R}\cap f_{\R}^{-1}(W^s_{\mathrm{loc}}(p_1)))\cdots),\]
and for a backward admissible sequence of the form $J=\overline{0}j_{-n}\cdots j_{-1}j_0\in\mathfrak{S}_{\mathrm{bwd}}^+$ we define
\[W^u_J(a, b)^+\equiv \mathcal{B}^+_{j_0,\R}\cap f_{\R}(\mathcal{B}^+_{j_{-1},\R}\cap\cdots\cap f_{\R}(\mathcal{B}^+_{j_{-n},\R}\cap f_{\R}(W^u_{\mathrm{loc}}(p_1)))\cdots).\]
Note that these submanifolds are well-defined even for the case $b=0$. Since $f^{-1} : \mathcal{B}^+_0\cap f(\mathcal{B}^+_1)\to\mathcal{B}^+_1$ and $f^{-1} : \mathcal{B}^+_1\cap f(\mathcal{B}^+_3)\to\mathcal{B}^+_3$ are crossed mappings of degree one, $W_{31\overline{0}}^s(a, b)^+$ is a vertical curve of degree one in $\mathcal{B}^+_{3,\R}$. Since $f : \mathcal{B}^+_0\cap f^{-1}(\mathcal{B}^+_2)\to\mathcal{B}^+_2$ is a crossed mappings of degree one and $f : \mathcal{B}^+_2\cap f^{-1}(\mathcal{B}^+_3)\to\mathcal{B}^+_3$ is a crossed mapping of degree two, $W_{\overline{0}23}^u(a, b)^+$ consists of either (i) a single $U$-shaped curve in $\mathcal{B}^+_{3,\R}$ from the right boundary of $\mathcal{B}^+_{3,\R}$ to itself or (ii) two mutually disjoint horizontal curves of degree one in $\mathcal{B}^+_{3,\R}$. This easily follows from an argument in the proof of Proposition 3.4 in~\cite{BS2}.

Let $(a, b)\in \mathcal{F}^-_{\R}$. For a forward admissible sequence of the form $I=i_0 i_1\cdots i_n\overline{0}\in\mathfrak{S}_{\mathrm{fwd}}^-$ we define
\[W^s_I(a, b)^-\equiv \mathcal{B}^-_{i_0, \R}\cap f_{\R}^{-1}(\mathcal{B}^-_{i_1, \R}\cap\cdots\cap f_{\R}^{-1}(\mathcal{B}^-_{i_n, \R}\cap f_{\R}^{-1}(W^s_{\mathrm{loc}}(p_1)))\cdots),\]
and for a backward admissible sequence of the form $J=\overline{43}j_{-n}\cdots j_{-1}j_0\in\mathfrak{S}_{\mathrm{bwd}}^-$ we define
\[W^u_J(a, b)^-\equiv \mathcal{B}^-_{j_0,\R}\cap f_{\R}(\mathcal{B}^-_{j_{-1},\R}\cap\cdots\cap f_{\R}(\mathcal{B}^-_{j_{-n},\R}\cap f_{\R}(W^u_{\mathrm{loc}}(p_3)))\cdots).\]
Note that these submanifolds are well-defined even for the case $b=0$. Since $f^{-1} : \mathcal{B}^-_0\cap f(\mathcal{B}^-_1)\to\mathcal{B}^-_1$ and $f^{-1} : \mathcal{B}^-_1\cap f(\mathcal{B}^-_4)\to\mathcal{B}^-_4$ are crossed mappings of degree one, $W_{41\overline{0}}^s(a, b)$ is a vertical curve of degree one in $\mathcal{B}^-_{4,\R}$. However, we need to be careful for $W^u_{\overline{43}4124}(a, b)^-$.

\begin{lmm}
Let $(a, b)\in \mathcal{F}^-_{\R}\cap\{b\ne 0\}$. Then, $W^u_{\overline{43}412}(a, b)^-$ consists of two mutually disjoint horizontal curves of degree one in $\mathcal{B}^-_{2,\R}$. 
\label{LMM:43412}
\end{lmm}

\begin{proof}
This follows from Proposition~\ref{PRP:43412}. 
\end{proof}

\begin{figure}
  \includegraphics[height=7cm]{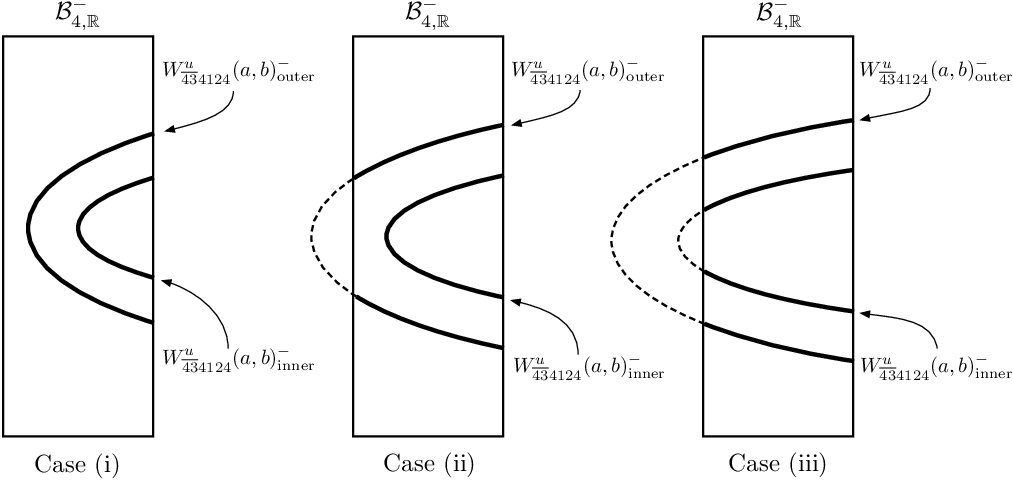}
  \caption{Outer and inner pieces of $W^u_{\overline{43}4124}(a, b)^-$.}
  \label{FIG:inner_and_outer}
\end{figure}

By tracing an argument in the proof of Proposition 3.4 in~\cite{BS2}, we see that $W^u_{\overline{43}4124}(a, b)^-$ consists of either (i) two mutually disjoint $U$-shaped curves in $\mathcal{B}^-_{4,\R}$ from the right boundary of $\mathcal{B}^-_{4,\R}$ to itself, (ii) one $U$-shaped curve as in (i) and two horizontal curves of degree one in $\mathcal{B}^-_{4,\R}$ all mutually disjoint, or (iii) four mutually disjoint horizontal curves of degree one in $\mathcal{B}^-_{4,\R}$ (see Figure~\ref{FIG:inner_and_outer}). 

Thanks to Lemma~\ref{LMM:43412}, we can speak of the upper piece $W^u_{\overline{43}412}(a, b)^-_{\mathrm{upper}}$ of $W^u_{\overline{43}412}(a, b)^-$ and the lower piece $W^u_{\overline{43}412}(a, b)^-_{\mathrm{lower}}$ of $W^u_{\overline{43}412}(a, b)^-$. This enables us to define the ``outer'' and  the ``inner'' pieces of $W^u_{\overline{43}4124}(a, b)^-$. More precisely,

\begin{dfn}
Let $(a, b)\in \mathcal{F}^-_{\R}\cap\{b<0\}$. Then, the \textit{inner piece} of $W^u_{\overline{43}4124}(a, b)^-$ is defined as $W^u_{\overline{43}4124}(a, b)^-_{\mathrm{inner}}\equiv\mathcal{B}^-_{4, \R}\cap f_{a, b}(W^u_{\overline{43}412}(a, b)^-_{\mathrm{upper}})$, and the \textit{outer piece} of $W^u_{\overline{43}4124}(a, b)^-$ is defined as $W^u_{\overline{43}4124}(a, b)^-_{\mathrm{outer}}\equiv\mathcal{B}^-_{4, \R}\cap f_{a, b}(W^u_{\overline{43}412}(a, b)^-_{\mathrm{lower}})$ (see Figure~\ref{FIG:inner_and_outer} again).
\label{DFN:inner_and_outer}
\end{dfn}

\subsection{Sides and signs}\label{subsection4.2}

First we define the notion of \textit{sides} of a real box. 

\begin{figure}
  \includegraphics[width=15cm]{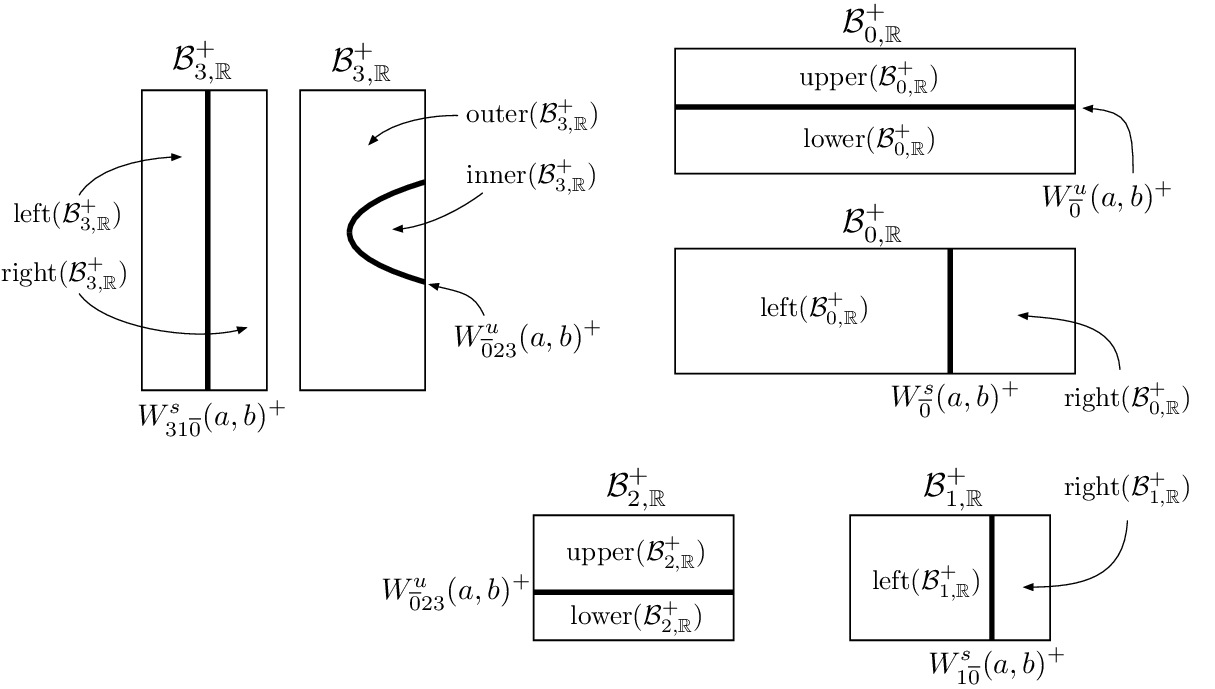}
  \caption{Special pieces and sides of $\mathcal{B}^+_{i, \R}$ for $(a, b)\in \mathcal{F}^+_{\R}\cap\{b>0\}$.}
  \label{FIG:special_varieties_B_pos}
\end{figure}

Let $(a, b)\in \mathcal{F}^+_{\R}\cap\{b>0\}$. By Lemma~\ref{LMM:real_disk} we know that $W^u_{\overline{0}}(a, b)^+$ is a horizontal curve between the right and the left boundaries of $\mathcal{B}^+_{0, \R}$. Hence $\mathcal{B}^+_{0, \R}\setminus W^u_{\overline{0}}(a, b)^+$ consists of two connected components, the one $\mathrm{upper}(\mathcal{B}^+_{0, \R})$ containing the upper boundary of $\mathcal{B}^+_{0, \R}$ and the one $\mathrm{lower}(\mathcal{B}^+_{0, \R})$ containing the lower boundary of $\mathcal{B}^+_{0, \R}$. Since $f : \mathcal{B}^+_0\cap f^{-1}(\mathcal{B}^+_2)\to\mathcal{B}^+_2$ is a crossed mapping of degree one, $W^u_{\overline{0}2}(a, b)^+$ is a horizontal curve between the right and the left boundaries of $\mathcal{B}^+_{2, \R}$. It follows that $\mathcal{B}^+_{2, \R}\setminus W^u_{\overline{0}2}(a, b)^+$ consists of two connected components, the one $\mathrm{upper}(\mathcal{B}^+_{2, \R})$ containing the upper boundary of $\mathcal{B}^+_{2, \R}$ and the one $\mathrm{lower}(\mathcal{B}^+_{2, \R})$ containing the lower boundary of $\mathcal{B}^+_{2, \R}$. Since $f : \mathcal{B}^+_2\cap f^{-1}(\mathcal{B}^+_3)\to\mathcal{B}^+_3$ is a crossed mapping of degree two, $W^u_{\overline{0}23}(a, b)^+$ is either a $U$-shaped curve from the right boundary of $\mathcal{B}^+_{3, \R}$ to itself or two mutually disjoint horizontal curves in $\mathcal{B}^+_{3, \R}$. Let $\mathrm{inner}(\mathcal{B}^+_{3, \R})$ be the connected component of $\mathcal{B}^+_{3, \R}\setminus W^u_{\overline{0}23}(a, b)^+$ which does not contain the upper and the lower boundaries of $\mathcal{B}^+_{3, \R}$ and let $\mathrm{outer}(\mathcal{B}^+_{3, \R})$ be the complement in $\mathcal{B}^+_{3, \R}$ of the union of  $W^u_{\overline{0}23}(a, b)^+$ and $\mathrm{inner}(\mathcal{B}^+_{3, \R})$. Since $W^s_{\overline{0}}(a, b)^+$, $W^s_{1\overline{0}}(a, b)^+$ and $W^s_{31\overline{0}}(a, b)^+$ are vertical curve between the upper and the lower boundaries of $\mathcal{B}^+_{0, \R}$, $\mathcal{B}^+_{1, \R}$ and $\mathcal{B}^+_{3, \R}$ respectively, we can define $\mathrm{right}(\mathcal{B}^+_{i, \R})$ and $\mathrm{left}(\mathcal{B}^+_{i, \R})$ for $i=0, 1, 3$ (see Figure~\ref{FIG:special_varieties_B_pos}). 

\begin{figure}
  \includegraphics[width=15.8cm]{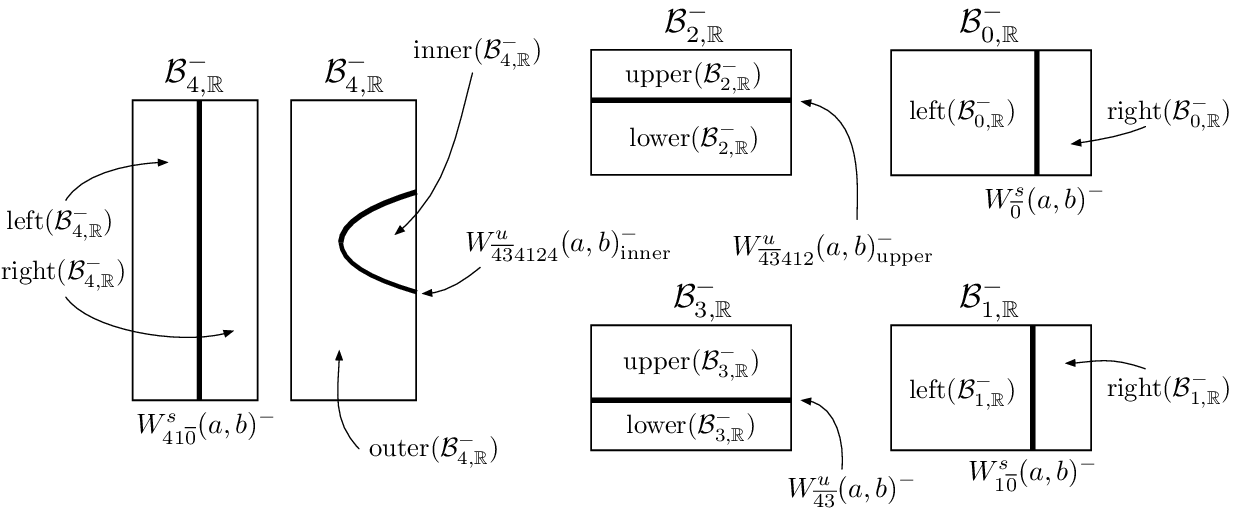}
  \caption{Special pieces and sides of $\mathcal{B}^-_{i, \R}$ for $(a, b)\in \mathcal{F}^-_{\R}\cap\{b<0\}$.}
  \label{FIG:special_varieties_B_neg}
\end{figure}

Let $(a, b)\in \mathcal{F}^-_{\R}\cap\{b<0\}$. We define $\mathrm{right}(\mathcal{B}^-_{0, \R})$ and $\mathrm{left}(\mathcal{B}^-_{0, \R})$ by using $W^s_{\overline{0}}(a, b)^-$, $\mathrm{right}(\mathcal{B}^-_{1, \R})$ and $\mathrm{left}(\mathcal{B}^-_{1, \R})$ by using $W^s_{1\overline{0}}(a, b)^-$, $\mathrm{right}(\mathcal{B}^-_{4, \R})$ and $\mathrm{left}(\mathcal{B}^-_{4, \R})$ by using $W^s_{41\overline{0}}(a, b)^-$, $\mathrm{upper}(\mathcal{B}^-_{3, \R})$ and $\mathrm{lower}(\mathcal{B}^-_{3, \R})$ by using $W^u_{\overline{43}}(a, b)^-$, $\mathrm{upper}(\mathcal{B}^-_{2, \R})$ and $\mathrm{lower}(\mathcal{B}^-_{2, \R})$ by using $W^u_{\overline{43}412}(a, b)^-_{\mathrm{upper}}$, and $\mathrm{outer}(\mathcal{B}^-_{4, \R})$ and $\mathrm{inner}(\mathcal{B}^-_{4, \R})$ by using $W^u_{\overline{43}4124}(a, b)^-_{\mathrm{inner}}$ (see Figure~\ref{FIG:special_varieties_B_neg}).

\begin{dfn}
We call $\mathrm{upper}(\mathcal{B}^{\pm}_{i, \R})$ the \textit{upper side}, $\mathrm{lower}(\mathcal{B}^{\pm}_{i, \R})$ the \textit{lower side}, $\mathrm{right}(\mathcal{B}^{\pm}_{i, \R})$ the \textit{right-hand side}, $\mathrm{left}(\mathcal{B}^{\pm}_{i, \R})$ the \textit{left-hand side}, $\mathrm{outer}(\mathcal{B}^{\pm}_{i, \R})$ the \textit{outer side}, $\mathrm{inner}(\mathcal{B}^{\pm}_{i, \R})$ the \textit{inner side} of a real box $\mathcal{B}^{\pm}_{i, \R}$.
\label{DFN:sides}
\end{dfn}

As in Definition~\ref{DFN:horizontal_curve}, the notion of horizontal and vertical curves can be extended to curves in $\mathrm{right}(\mathcal{B}^{\pm}_{i, \R})$ and in $\mathrm{left}(\mathcal{B}^{\pm}_{i, \R})$ in an obvious way (for appropriate $i$). It can be also extended to curves in the closures of $\mathrm{right}(\mathcal{B}^{\pm}_{i, \R})$ and $\mathrm{left}(\mathcal{B}^{\pm}_{i, \R})$. These notions will be used in Propositions~\ref{PRP:pieces_positive} and~\ref{PRP:pieces_negative} below.

Next we define the notion of \textit{sign pairs} of a crossed mapping. Choose an admissible transition $(i, j)\in\mathfrak{T}^{\pm}$. Assume first that the degree of the crossed mapping $f : \mathcal{B}^{\pm}_i\cap f^{-1}(\mathcal{B}^{\pm}_j)\to\mathcal{B}^{\pm}_j$ is one. In this case $f^{-1} : \mathcal{B}^{\pm}_j\cap f(\mathcal{B}^{\pm}_i)\to\mathcal{B}^{\pm}_i$ is also a crossed mapping of degree one. First, take an \textit{oriented} horizontal curve $C$ of degree one in $\mathcal{B}^{\pm}_{i, \R}$ from the right boundary to the left boundary of $\mathcal{B}^{\pm}_{i, \R}$. Then, $f_{\R}(C)\cap\mathcal{B}^{\pm}_{j, \R}$ is an oriented horizontal curve of degree one in $\mathcal{B}^{\pm}_{j, \R}$. Hence it is a curve either from the right boundary to the left boundary or from the left boundary to the right boundary of $\mathcal{B}^{\pm}_{j, \R}$. In the first case we associate $\varepsilon_u\equiv +$ and in the second case we associate $\varepsilon_u\equiv -$. 

Next, take an oriented vertical curve $C$ of degree one in $\mathcal{B}^{\pm}_{j, \R}$ from the lower boundary to the upper boundary of $\mathcal{B}^{\pm}_{j, \R}$. Then, $f^{-1}_{\R}(C)\cap\mathcal{B}^{\pm}_{i, \R}$ is an oriented vertical curve of degree one in $\mathcal{B}^{\pm}_{i, \R}$. Hence it is a curve either from the lower boundary to the upper boundary or from the upper boundary to the lower boundary of $\mathcal{B}^{\pm}_{i, \R}$. In the first case we associate $\varepsilon_v\equiv +$ and in the second case we associate $\varepsilon_v\equiv -$. When the degree of the crossed mapping $f : \mathcal{B}^{\pm}_i\cap f^{-1}(\mathcal{B}^{\pm}_j)\to\mathcal{B}^{\pm}_j$ is two, we associate $(\varepsilon_u, \varepsilon_v) \equiv (\ast, \ast)$.

\begin{dfn}
We call the pair $(\varepsilon_u, \varepsilon_v)$ defined above the \textit{sign pair} of the admissible transition $(i, j)\in\mathfrak{T}^{\pm}$.
\label{DFN:sign_pair}
\end{dfn}

Using the notion of sign pairs, the following list of transitions of sides is obtained for the case $(a, b)\in\mathcal{F}^+_{\R}\cap\{b>0\}$.

\begin{lmm}
If $(a, b)\in\mathcal{F}^+_{\R}\cap\{b>0\}$, then we have 
\begin{enumerate}
\renewcommand{\labelenumi}{(\roman{enumi})}
\item $f(\mathrm{lower}(\mathcal{B}^+_{0, \R}))\cap \mathcal{B}^+_{0, \R}\subset \mathrm{lower}(\mathcal{B}^+_{0, \R})$ and $f(\mathrm{left}(\mathcal{B}^+_{0, \R}))\cap \mathcal{B}^+_{0, \R}\subset \mathrm{left}(\mathcal{B}^+_{0, \R})$, 
\item $f(\mathrm{lower}(\mathcal{B}^+_{0, \R}))\cap \mathcal{B}^+_{2, \R}\subset \mathrm{upper}(\mathcal{B}^+_{2, \R})$, 
\item $f(\mathrm{lower}(\mathcal{B}^+_{0, \R}))\cap \mathcal{B}^+_{3, \R}\subset \mathrm{outer}(\mathcal{B}^+_{3, \R})$, 
\item $f(\mathcal{B}^+_{1, \R})\cap \mathcal{B}^+_{0, \R}\subset \mathrm{lower}(\mathcal{B}^+_{0, \R})$ and $f(\mathrm{left}(\mathcal{B}^+_{1, \R}))\cap \mathcal{B}^+_{0, \R}\subset \mathrm{left}(\mathcal{B}^+_{0, \R})$,
\item $f(\mathrm{upper}(\mathcal{B}^+_{2, \R}))\cap \mathcal{B}^+_{2, \R}\subset \mathrm{upper}(\mathcal{B}^+_{2, \R})$,
\item $f(\mathrm{upper}(\mathcal{B}^+_{2, \R}))\cap \mathcal{B}^+_{3, \R}\subset \mathrm{outer}(\mathcal{B}^+_{3, \R})$, 
\item $f(\mathrm{right}(\mathcal{B}^+_{3, \R}))\cap \mathcal{B}^+_{1, \R}\subset \mathrm{left}(\mathcal{B}^+_{1, \R})$.
\end{enumerate}
\label{LMM:sides_positive}
\end{lmm}

\begin{proof}
When $(a, b)\in\mathcal{F}^+_{\R}\cap\{b>0\}$, we first examine the sign pair for every admissible transition $(i, j)\in\mathfrak{T}^+$. By referring to Figure~\ref{FIG:special_varieties_B_pos}, the sign pairs are given by $(\varepsilon_u, \varepsilon_v)=(+, +)$ for $(i, j)=(0, 0)$, $(\varepsilon_u, \varepsilon_v)=(-, -)$ for $(i, j)=(0, 2)$, $(\varepsilon_u, \varepsilon_v)=(\ast, \ast)$ for $(i, j)=(0, 3)$, $(\varepsilon_u, \varepsilon_v)=(+, +)$ for $(i, j)=(1, 0)$, $(\varepsilon_u, \varepsilon_v)=(-, -)$ for $(i, j)=(2, 2)$, $(\varepsilon_u, \varepsilon_v)=(\ast, \ast)$ for $(i, j)=(2, 3)$ and $(\varepsilon_u, \varepsilon_v)=(-, -)$ for $(i, j)=(3, 1)$. These claims obviously hold when $b>0$ is close to zero. Since the boxes vary continuously with respect to $(a, b)\in\mathcal{F}^+_{\R}\cap\{b>0\}$, they hold for any $(a, b)\in\mathcal{F}^+_{\R}\cap\{b>0\}$. By using this list, it is easy to show that the claims (i), (v), (vi) and (vii) hold. 

To prove the rest of the claims we first consider the case $b>0$ close to zero and then use the continuity argument. When $b>0$ close to zero, the $y$-coordinate of any point in $\mathcal{B}^+_{0, \R}$ is larger than the $y$-coordinate of any point in $\mathcal{B}^+_{2, \R}$, hence (vi) implies (ii). When $b>0$ close to zero, the $y$-coordinate of any point in $\mathcal{B}^+_{0, \R}$ is larger than the $y$-coordinate of any point in $\mathcal{B}^+_{2, \R}$, hence (vi) implies (iii). When $b>0$ close to zero, the $y$-coordinate of any point in $\mathcal{B}^+_{0, \R}$ is larger than the $y$-coordinate of any point in $\mathcal{B}^+_{1, \R}$, hence (vi) implies (iv). 
\end{proof}

In the case $(a, b)\in\mathcal{F}^-_{\R}\cap\{b<0\}$, let $\widetilde{\mathcal{B}}^-_{4, \R}$ be the closure of the subregion of $\mathcal{B}^-_{4, \R}$ surrounded by $W^u_{\overline{43}4}(a, b)^-$, the right boundary and the left boundary of $\mathcal{B}^-_{4, \R}$ (the left boundary of $\mathcal{B}^-_{4, \R}$ is not necessary when $W^u_{\overline{43}4}(a, b)^-$ consists of a single curve from the right boundary of $\mathcal{B}^-_{4, \R}$ to itself), and let $\widetilde{\mathcal{B}}^-_{1, \R}\equiv f(\widetilde{\mathcal{B}}^-_{4, \R})\cap \mathcal{B}^-_{1, \R}$. Then, the following list of transitions of sides is obtained for the case $(a, b)\in\mathcal{F}^-_{\R}\cap\{b<0\}$.

\begin{lmm}
If $(a, b)\in\mathcal{F}^-_{\R}\cap\{b<0\}$, then we have 
\begin{enumerate}
\renewcommand{\labelenumi}{(\roman{enumi})}
\item $f(\mathrm{left}(\mathcal{B}^-_{0, \R}))\cap \mathcal{B}^-_{0, \R}\subset \mathrm{left}(\mathcal{B}^-_{0, \R})$, 
\item $f(\mathcal{B}^-_{0, \R})\cap \mathcal{B}^-_{2, \R}\subset \mathrm{lower}(\mathcal{B}^-_{2, \R})$,
\item $f(\mathrm{left}(\mathcal{B}^-_{1, \R}))\cap \mathcal{B}^-_{0, \R}\subset \mathrm{left}(\mathcal{B}^-_{0, \R})$, 
\item $f(\widetilde{\mathcal{B}}^-_{1, \R})\cap \mathcal{B}^-_{2, \R}\subset \mathrm{lower}(\mathcal{B}^-_{2, \R})$, 
\item $f(\mathrm{lower}(\mathcal{B}^-_{2, \R}))\cap \mathcal{B}^-_{4, \R}\subset \mathrm{outer}(\mathcal{B}^-_{4, \R})$, 
\item $f(\mathrm{upper}(\mathcal{B}^-_{3, \R}))\cap \mathcal{B}^-_{4, \R}\subset \mathrm{outer}(\mathcal{B}^-_{4, \R})$, 
\item $f(\mathrm{right}(\mathcal{B}^-_{4, \R}))\cap \mathcal{B}^-_{1, \R}\subset \mathrm{left}(\mathcal{B}^-_{1, \R})$, 
\item $f(\widetilde{\mathcal{B}}^-_{4, \R}) \cap \mathcal{B}^-_{3, \R}\subset \mathrm{upper}(\mathcal{B}^-_{3, \R})$.
\end{enumerate}
\label{LMM:sides_negative}
\end{lmm}

\begin{proof}
When $(a, b)\in\mathcal{F}^-_{\R}\cap\{b<0\}$, we first examine the sign pair for every admissible transition $(i, j)\in\mathfrak{T}^-$. By referring to Figure~\ref{FIG:special_varieties_B_neg}, the sign pairs are given by $(\varepsilon_u, \varepsilon_v)=(+, -)$ for $(i, j)=(0, 0)$, $(\varepsilon_u, \varepsilon_v)=(+, -)$ for $(i, j)=(0, 2)$, $(\varepsilon_u, \varepsilon_v)=(+, -)$ for $(i, j)=(1, 0)$, $(\varepsilon_u, \varepsilon_v)=(+, -)$ for $(i, j)=(1, 2)$, $(\varepsilon_u, \varepsilon_v)=(\ast, \ast)$ for $(i, j)=(2, 4)$, $(\varepsilon_u, \varepsilon_v)=(\ast, \ast)$ for $(i, j)=(3, 4)$, $(\varepsilon_u, \varepsilon_v)=(-, +)$ for $(i, j)=(4, 1)$ and $(\varepsilon_u, \varepsilon_v)=(-, +)$ for $(i, j)=(4, 3)$. Using this list, it is easy to show the claims (i), (iii), (v), (vii) and (viii). The claim (iv) immediately follows from the definition of $W^u_{\overline{43}412}(a, b)^-_{\mathrm{upper}}$.

To prove the rest of the claims we argue as in Lemma~\ref{LMM:sides_positive}. When $b<0$ close to zero, the $y$-coordinate of any point in $\mathcal{B}^-_{0, \R}$ is larger than the $y$-coordinate of any point in $\mathcal{B}^-_{1, \R}$, hence (iv) implies (ii). Similarly, when $b<0$ close to zero, the $y$-coordinate of any point in $\mathcal{B}^-_{2, \R}$ is larger than the $y$-coordinate of any point in $\mathcal{B}^-_{3, \R}$, hence  (v) implies (vi).
\end{proof}

\subsection{Special pieces}\label{subsection4.3}

In this subsection we show that a condition on the intersection between special pieces controls a certain global dynamical behavior. Below $\mathrm{card}(X)$ means the cardinality of a set $X$ counted \textit{without} multiplicity.

\begin{figure}
  \includegraphics[height=7cm]{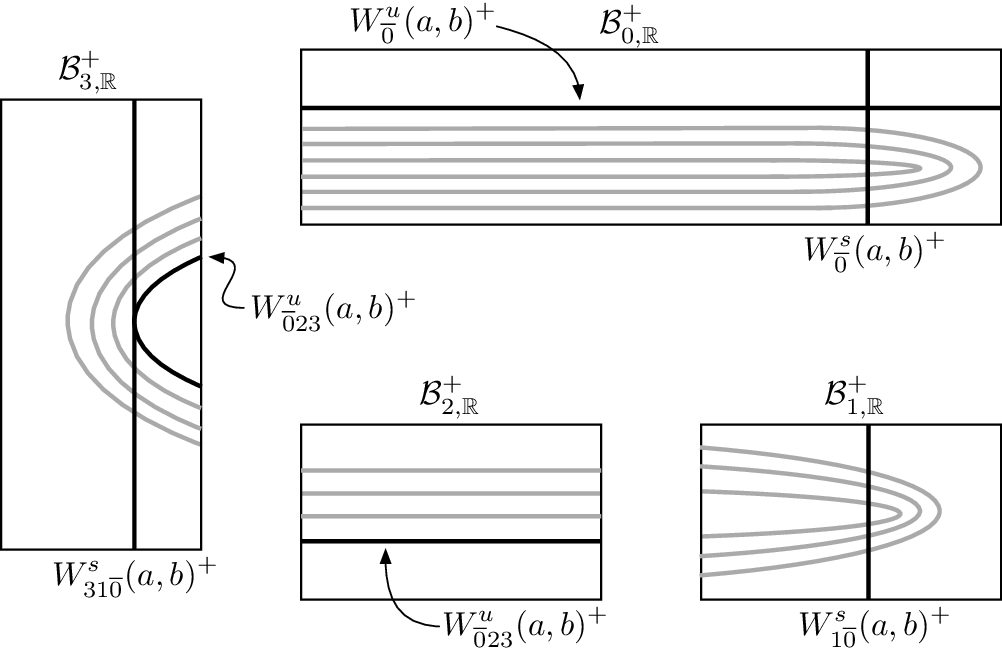}
  \caption{Special pieces (black) and $W^u_J(a, b)^+$ (gray) for $(a, b)\in \mathcal{F}^+_{\R}\cap\{b>0\}$.}
  \label{FIG:special_varieties_Wu_pos}
\end{figure}

First let us consider the case $(a, b)\in \mathcal{F}^+_{\R}\cap\{b>0\}$. 

\begin{prp}
Let $(a, b)\in \mathcal{F}^+_{\R}\cap\{b>0\}$. If $\mathrm{card}(W_{31\overline{0}}^s(a, b)^+\cap W_{\overline{0}23}^u(a, b)^+)\geq 1$, then
\begin{enumerate}
\renewcommand{\labelenumi}{(\roman{enumi})}
\item each connected component of $W^u_J(a, b)^+\cap \mathrm{left}(\mathcal{B}^+_{0, \R})$ is a horizontal curve of degree one in $\mathrm{left}(\mathcal{B}^+_{0, \R})$ and is contained in $\overline{\mathrm{lower}(\mathcal{B}^+_{0, \R})}$ for any backward admissible sequence of the form $J=\overline{0}j_{-n}\cdots j_{-1}0\in\mathfrak{S}_{\mathrm{bwd}}^+$,
\item  each connected component of $W^u_J(a, b)^+\cap \mathrm{left}(\mathcal{B}^+_{1, \R})$ is a horizontal curve of degree one in $\mathrm{left}(\mathcal{B}^+_{1, \R})$ for any backward admissible sequence of the form $J=\overline{0}j_{-n}\cdots j_{-1}1\in\mathfrak{S}_{\mathrm{bwd}}^+$,
\item each connected component of $W^u_J(a, b)^+\cap \mathcal{B}^+_{2, \R}$ is a horizontal curve of degree one in $\mathcal{B}^+_{2, \R}$ and is contained in $\overline{\mathrm{upper}(\mathcal{B}^+_{2, \R})}$ for any backward admissible sequence of the form $J=\overline{0}j_{-n}\cdots j_{-1}2\in\mathfrak{S}_{\mathrm{bwd}}^+$,
\item each connected component of $W^u_J(a, b)^+\cap \mathrm{right}(\mathcal{B}^+_{3, \R})$ is a horizontal curve of degree one in $\mathrm{right}(\mathcal{B}^+_{3, \R})$ for any backward admissible sequence of the form $J=\overline{0}j_{-n}\cdots j_{-1}3\in\mathfrak{S}_{\mathrm{bwd}}^+$ (see Figure~\ref{FIG:special_varieties_Wu_pos}).
\end{enumerate}
If moreover $\mathrm{card}(W_{31\overline{0}}^s(a, b)^+\cap W_{\overline{0}23}^u(a, b)^+)=2$ holds, then $\mathrm{left}(\mathcal{B}^+_{0, \R})$, $\mathrm{left}(\mathcal{B}^+_{1, \R})$ and $\mathrm{right}(\mathcal{B}^+_{3, \R})$ in the above statements can be replaced by $\overline{\mathrm{left}(\mathcal{B}^+_{0, \R})}$, $\overline{\mathrm{left}(\mathcal{B}^+_{1, \R})}$ and $\overline{\mathrm{right}(\mathcal{B}^+_{3, \R})}$ respectively (see Figure~\ref{FIG:special_varieties_Wu_pos} again).
\label{PRP:pieces_positive}
\end{prp}

\begin{proof}
We prove the claim for $\mathrm{card}(W_{31\overline{0}}^s(a, b)^+\cap W_{\overline{0}23}^u(a, b)^+)\geq 1$ by induction on $n$. 

When $n=0$, the claim (i) holds since $W^u_{\mathrm{loc}}(p_0)$ is a horizontal curve of degree one in $\mathcal{B}^+_{0, \R}$, the claim (ii) holds since $W^u_J(a, b)^+\cap \mathrm{left}(\mathcal{B}^+_{1, \R})$ is empty when $j_0=1$, the claim (iii) holds since $f : \mathcal{B}^+_0\cap f^{-1}(\mathcal{B}^+_2)\to \mathcal{B}^+_2$ is a crossed mapping of degree one, the claim (iv) holds by the assumption $\mathrm{card}(W_{31\overline{0}}^s(a, b)^+\cap W_{\overline{0}23}^u(a, b)^+)\geq 1$.

Assume that the claims hold for $k=n-1$ and consider the case $k=n$. Choose a backward admissible sequence $J=\overline{0}j_{-k}\cdots j_{-1}j_0\in\mathfrak{S}_{\mathrm{bwd}}^+$ and write $J'=\overline{0}j_{-k}\cdots j_{-1}\in\mathfrak{S}_{\mathrm{bwd}}^+$. 

If $j_0=0$, then either $j_{-1}=0$ or $j_{-1}=1$ holds. Suppose first the case $j_{-1}=0$. Since $f : \mathcal{B}^+_0\cap f^{-1}(\mathcal{B}^+_0)\to \mathcal{B}^+_0$ is a crossed mapping of degree one and since each connected component of $W^u_{J'}(a, b)^+\cap \mathrm{left}(\mathcal{B}^+_{0, \R})$ is a horizontal curve of degree one in $\mathrm{left}(\mathcal{B}^+_{0, \R})$ by induction assumption, each connected component of $W^u_J(a, b)^+\cap \mathrm{left}(\mathcal{B}^+_{0, \R})$ is a horizontal curve of degree one in $\mathrm{left}(\mathcal{B}^+_{0, \R})$. It is contained in $\overline{\mathrm{lower}(\mathcal{B}^+_{0, \R})}$ thanks to (i) of Lemma~\ref{LMM:sides_positive}. The proof for the case $j_{-1}=1$ is identical, and this proves the claim (i) for $k=n$. 

If $j_0=1$, then $j_{-1}=3$ holds. Since $f : \mathcal{B}^+_3\cap f^{-1}(\mathcal{B}^+_1)\to \mathcal{B}^+_1$ is a crossed mapping of degree one and since each connected component of $W^u_{J'}(a, b)^+\cap \mathrm{right}(\mathcal{B}^+_{3, \R})$ is a horizontal curve of degree one in $\mathrm{right}(\mathcal{B}^+_{3, \R})$ by induction assumption, each component of $W^u_J(a, b)^+\cap \mathrm{left}(\mathcal{B}^+_{1, \R})$ is a horizontal curve of degree one in $\mathrm{left}(\mathcal{B}^+_{1, \R})$. This proves (ii) for $k=n$. 

If $j_0=2$, then either $j_{-1}=0$ or $j_{-1}=2$ holds. Suppose first the case $j_{-1}=0$. Since $f : \mathcal{B}^+_0\cap f^{-1}(\mathcal{B}^+_2)\to \mathcal{B}^+_2$ is a crossed mapping of degree one and since each connected component of $W^u_{J'}(a, b)^+\cap \mathrm{left}(\mathcal{B}^+_{0, \R})$ is a horizontal curve of degree one in $\mathrm{left}(\mathcal{B}^+_{0, \R})$ by induction assumption, that each connected component of $W^u_J(a, b)^+\cap \mathcal{B}^+_{2, \R}$ is a horizontal curve of degree one in $\mathcal{B}^+_{2, \R}$. It is contained in $\overline{\mathrm{upper}(\mathcal{B}^+_{2, \R})}$ thanks to (ii) of Lemma~\ref{LMM:sides_positive}. The proof for the case $j_{-1}=2$ is identical, and this proves the claim (iii) for $k=n$. 

If $j_0=3$, then either $j_{-1}=0$ or $j_{-1}=2$ holds. Suppose first the case $j_{-1}=2$. Since $f : \mathcal{B}^+_2\cap f^{-1}(\mathcal{B}^+_3)\to \mathcal{B}^+_3$ is a crossed mapping of degree two and since each connected component of $W^u_{J'}(a, b)^+\cap \mathcal{B}^+_{2, \R}$ is a horizontal curve of degree one in $\overline{\mathrm{upper}(\mathcal{B}^+_{2, \R})}$ by induction assumption, each connected component of $W^u_J(a, b)^+\cap \mathrm{right}(\mathcal{B}^+_{3, \R})$ is a horizontal curve of degree one in $\mathrm{right}(\mathcal{B}^+_{3, \R})$ by the assumption $\mathrm{card}(W_{31\overline{0}}^s(a, b)^+\cap W_{\overline{0}23}^u(a, b)^+)\geq 1$. The proof for the case $j_{-1}=0$ is identical, and this proves the claim (iv) for $k=n$. 

The proof for the case $\mathrm{card}(W_{31\overline{0}}^s(a, b)^+\cap W_{\overline{0}23}^u(a, b)^+)=2$ is similar, hence omitted. 
\end{proof}

Let us write $K_{\R}\equiv K_{a, b}\cap\R^2$. To globalize this statement, we need

\begin{lmm}
We have
\[\bigcup_I W^s_I(a, b)^+\supset W^s(p_1)\cap K_{\R},\]
where $I$ runs over all forward admissible sequences of the form $I=i_0 i_1\cdots i_n\overline{0}\in\mathfrak{S}_{\mathrm{fwd}}^+$, and
\[\bigcup_J W^u_J(a, b)^+\supset W^u(p_1)\cap K_{\R},\]
where $J$ runs over all backward admissible sequences of the form $J=\overline{0}j_{-n}\cdots j_{-1}j_0\in\mathfrak{S}_{\mathrm{bwd}}^+$. 
\label{LMM:globalize_positive}
\end{lmm}

\begin{proof}
This is an easy consequence of Proposition~\ref{PRP:coding_positive}. 
\end{proof}

As a consequence of this lemma we show that the special intersection determines the non-existence of tangencies between $W^u(p_1)$ and $W^s(p_1)$ when $(a, b)\in \mathcal{F}^+_{\R}\cap\{b>0\}$. 

\begin{cor}
Let $(a, b)\in \mathcal{F}^+_{\R}\cap\{b>0\}$. If $\mathrm{card}(W_{31\overline{0}}^s(a, b)^+\cap W_{\overline{0}23}^u(a, b)^+)=2$, then there is no tangency between $W^u(p_1)$ and $W^s(p_1)$.
\label{COR:no_tangency_positive}
\end{cor}

\begin{proof}
From (iii) of Proposition~\ref{PRP:pieces_positive} we see $\mathrm{card}(W_{31\overline{0}}^s(a, b)^+\cap W_J^u(a, b)^+)=2$ and hence there is no tangency between $W_J^u(a, b)^+$ and $W_{31\overline{0}}^s(a, b)^+$ for any backward admissible sequence of the form $J=\overline{0}j_{-n}\cdots j_{-1}3\in\mathfrak{S}_{\mathrm{bwd}}^+$.

It is enough to show that if there exists no tangency between $W_J^u(a, b)^+$ and $W_{31\overline{0}}^s(a, b)^+$ then there exists no tangency between $W^u(p_1)$ and $W^s(p_1)$. Assume that there is a tangency $q\in W^u(p_1)\cap W^s(p_1)$. Then, $f^n(q)\in W_{31\overline{0}}^s(a, b)^+$ for $n\geq 0$ sufficiently large. Since $f^n(q)\in W^u(p_1)\cap W^s(p_1)\subset W^u(p_1)\cap K_{\R}$, we can find a backward admissible sequence of the form $J=\overline{0}j_{-n}\cdots j_{-1}j_0\in\mathfrak{S}_{\mathrm{bwd}}^+$ so that $f^n(q)\in W^u_J(a, b)^+$ by Lemma~\ref{LMM:globalize_positive}. Since $q\in W^u(p_1)\cap W^s(p_1)$ is a tangency, $f^n(q)$ is a tangency between $W_J^u(a, b)^+$ and $W_{31\overline{0}}^s(a, b)^+$, a contradiction. 
\end{proof}

\begin{figure}
  \includegraphics[height=7cm]{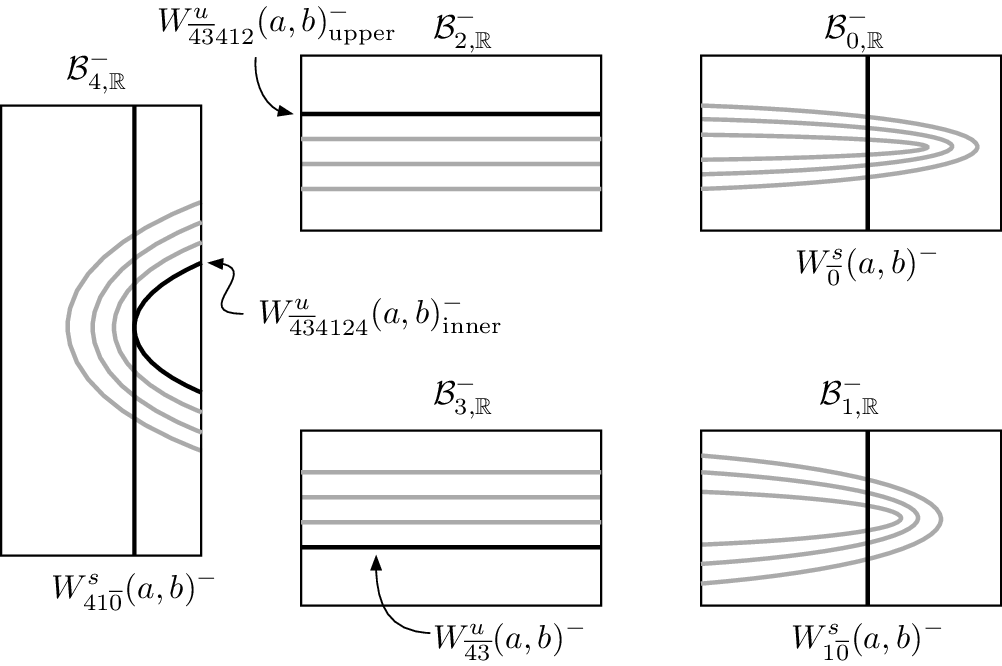}
  \caption{Special pieces (black) and $W^u_J(a, b)^-$ (gray) for $(a, b)\in \mathcal{F}^-_{\R}\cap\{b<0\}$.}
  \label{FIG:special_varieties_Wu_neg}
\end{figure}

Next let us show the corresponding claims for $(a, b)\in \mathcal{F}^-_{\R}\cap\{b<0\}$. 

\begin{prp}
Let $(a, b)\in \mathcal{F}^-_{\R}\cap\{b<0\}$. If $\mathrm{card}(W_{41\overline{0}}^s(a, b)^-\cap W_{\overline{43}4124}^u(a, b)^-_{\mathrm{inner}})\geq 1$, then
\begin{enumerate}
\renewcommand{\labelenumi}{(\roman{enumi})}
\item each connected component of $W^u_J(a, b)^-\cap \mathrm{left}(\mathcal{B}^-_{0, \R})$ is a horizontal curve of degree one in $\mathrm{left}(\mathcal{B}^-_{0, \R})$ for any backward admissible sequence of the form $J=\overline{43}j_{-n}\cdots j_{-1}0\in\mathfrak{S}_{\mathrm{bwd}}^-$,
\item  each connected component of $W^u_J(a, b)^-\cap \mathrm{left}(\mathcal{B}^-_{1, \R})$ is a horizontal curve of degree one in $\mathrm{left}(\mathcal{B}^-_{1, \R})$ and is contained in $\widetilde{\mathcal{B}}^-_{1, \R}$ for any backward admissible sequence of the form $J=\overline{43}j_{-n}\cdots j_{-1}1\in\mathfrak{S}_{\mathrm{bwd}}^-$,
\item each connected component of $W^u_J(a, b)^-\cap \mathcal{B}^-_{2, \R}$ is a horizontal curve of degree one in $\mathcal{B}^-_{2, \R}$ and is contained in $\overline{\mathrm{lower}(\mathcal{B}^-_{2, \R})}$ for any backward admissible sequence of the form $J=\overline{43}j_{-n}\cdots j_{-1}2\in\mathfrak{S}_{\mathrm{bwd}}^-$,
\item each connected component of $W^u_J(a, b)^-\cap \mathcal{B}^-_{3, \R}$ is a horizontal curve of degree one in $\mathcal{B}^-_{3, \R}$ and is contained in $\overline{\mathrm{upper}(\mathcal{B}^-_{3, \R})}$ for any backward admissible sequence of the form $J=\overline{43}j_{-n}\cdots j_{-1}3\in\mathfrak{S}_{\mathrm{bwd}}^-$,
\item each connected component of $W^u_J(a, b)^-\cap \mathrm{right}(\mathcal{B}^-_{4, \R})$ is a horizontal curve of degree one in $\mathrm{right}(\mathcal{B}^-_{4, \R})$ and is contained in $\widetilde{\mathcal{B}}^-_{4, \R}$ for any backward admissible sequence of the form $J=\overline{43}j_{-n}\cdots j_{-1}4\in\mathfrak{S}_{\mathrm{bwd}}^-$ (see Figure~\ref{FIG:special_varieties_Wu_neg}).
\end{enumerate}
If moreover $\mathrm{card}(W_{41\overline{0}}^s(a, b)^-\cap W_{\overline{43}4124}^u(a, b)^-_{\mathrm{inner}})=2$, then $\mathrm{left}(\mathcal{B}^-_{0, \R})$, $\mathrm{left}(\mathcal{B}^-_{1, \R})$ and $\mathrm{right}(\mathcal{B}^-_{4, \R})$ in the above statements can be replaced by $\overline{\mathrm{left}(\mathcal{B}^-_{0, \R})}$, $\overline{\mathrm{left}(\mathcal{B}^-_{1, \R})}$ and $\overline{\mathrm{right}(\mathcal{B}^-_{4, \R})}$ respectively (see Figure~\ref{FIG:special_varieties_Wu_neg} again).
\label{PRP:pieces_negative}
\end{prp}

\begin{proof}
Together with the definition of $\widetilde{\mathcal{B}}^-_{1, \R}$ and $\widetilde{\mathcal{B}}^-_{4, \R}$, the proof is similar to Proposition~\ref{PRP:pieces_positive}, hence omitted. 
\end{proof}

The proof of the following lemma is identical to the case $(a, b)\in \mathcal{F}^+_{\R}\cap\{b>0\}$. 

\begin{lmm}
We have
\[\bigcup_I W^s_I(a, b)^-\supset W^s(p_1)\cap K_{\R},\]
where $I$ runs over all forward admissible sequences of the form $I=i_0 i_1\cdots i_n\overline{0}\in\mathfrak{S}_{\mathrm{fwd}}^-$, and
\[\bigcup_J W^u_J(a, b)^-\supset W^u(p_3)\cap K_{\R},\]
where $J$ runs over all backward admissible sequences of the form $J=\overline{43}j_{-n}\cdots j_{-1}j_0\in\mathfrak{S}_{\mathrm{bwd}}^-$.
\label{LMM:globalize_negative}
\end{lmm}

As a consequence of this lemma we show that the special intersection determines the non-existence of tangencies between $W^u(p_3)$ and $W^s(p_1)$ when $(a, b)\in \mathcal{F}^-_{\R}\cap\{b<0\}$.

\begin{cor}
Let $(a, b)\in \mathcal{F}^-_{\R}\cap\{b<0\}$. If $\mathrm{card}(W_{41\overline{0}}^s(a, b)^-\cap W^u_{\overline{43}4124}(a, b)^-_{\mathrm{inner}})=2$, then there is no tangency between $W^u(p_3)$ and $W^s(p_1)$.
\label{COR:no_tangency_negative}
\end{cor}

\newpage

\section{Synthesis: Proof of the Main Theorem}\label{section5}

In this section we integrate the ideas developed in the previous sections to finish the proof of the Main Theorem. To achieve this we analyze carefully the complex tangency loci $\mathcal{T}^{\pm}$ (see Definition~\ref{DFN:tangency_loci}) and their real sections.

\subsection{Maximal entropy}\label{subsection5.1}

The purpose of this subsection is to show that the intersections of certain special pieces of $W^{u/s}(p_i)$ characterize the H\'enon maps with maximal entropy. Namely, we prove

\begin{thm}[\textbf{Maximal Entropy}]
When $(a, b)\in \mathcal{F}^+_{\R}\cap\{b>0\}$, we have $\htop(f_{a, b}|_{\R^2})=\log 2$ iff $\mathrm{card}(W_{31\overline{0}}^s(a, b)^+\cap W_{\overline{0}23}^u(a, b)^+)\geq 1$. When $(a, b)\in \mathcal{F}^-_{\R}\cap\{b<0\}$, we have $\htop(f_{a, b}|_{\R^2})=\log 2$ iff ${\mathrm{card}}(W_{41\overline{0}}^s(a, b)^-\cap W^u_{\overline{43}4124}(a, b)^-_{\mathrm{inner}})\geq 1$.
\label{THM:maximalentropy}
\end{thm}

Before proving this theorem, let us recall the following facts. For $f=f_{a, b} : \C^2\to\C^2$ with $(a, b)\in\R\times\R^{\times}$, it has been shown in Theorem 10.1 of~\cite{BLS} that the condition: 
\begin{quote}
(1) $\htop(f_{\R})=\log 2$
\end{quote}
is equivalent to
\begin{quote}
(2) for any saddle periodic point $p\in \C^2$, we have $V^u(p)\cap V^s(p)\subset\R^2$. 
\end{quote}
Let us consider a stronger condition: 
\begin{quote}
(2$'$) for any saddle periodic points $p, q\in \C^2$, we have $V^u(p)\cap V^s(q)\subset\R^2$.
\end{quote}

\begin{lmm}
The condition (2$'$) is equivalent to (2), hence to (1).
\label{LMM:heteroclinic}
\end{lmm}

\begin{proof}
Since we know that (2) implies (1) and (2$'$) implies (2), it is enough to show that (1) implies (2$'$). Suppose that (1) holds. By Theorem 10.1 of~\cite{BLS} we see that the filled Julia set of $f$ is contained in $\R^2$. Since every point in $V^u(p)\cap V^s(q)$ has forward and backward bounded orbits, the condition (2$'$) follows. 
\end{proof}

\begin{proof}[Proof of Theorem~\ref{THM:maximalentropy}]
Consider first the case $(a, b)\in \mathcal{F}^+_{\R}\cap\{b>0\}$. Choose any point $q\in V^u(p_1)\cap V^s(p_1)$ with $q\ne p_1$ and assume that $\mathrm{card}(W_{31\overline{0}}^s(a, b)^+\cap W_{\overline{0}23}^u(a, b)^+)\geq 1$ holds. Replacing $q$ by $f^m(q)$, if necessary, we may assume $q\in V^u_{\mathrm{loc}}(p_1)$. Since $q\in K_{a, b}$ and $q\ne p_1$, there exists $0i_1i_2\cdots\in\mathfrak{S}^+_{\mathrm{fwd}}$ different from $\overline{0}$ so that $f^n(q)\in \mathcal{B}^+_{i_n}$ holds for $n\geq 0$ by Proposition~\ref{PRP:coding_positive}. By taking $m\in\Z$ as large as possible, we may assume $i_1\ne 0$. Then, there exists $N\geq 0$ so that $i_1\cdots i_{N}=2\cdots 2$ (when $N=0$ this term disappears) and $i_{N+1}=3$. Since $f : \mathcal{B}^+_0\cap f^{-1}(\mathcal{B}^+_2)\to\mathcal{B}^+_2$ and $f : \mathcal{B}^+_2\cap f^{-1}(\mathcal{B}^+_2)\to\mathcal{B}^+_2$ are both crossed mappings of degree one, 
\[V^u_{\overline{0}2\cdots 2}(a, b)^+ \equiv \mathcal{B}^+_2\cap f(\mathcal{B}^+_2\cap\cdots\cap f(\mathcal{B}^+_2\cap f(V^u_{\mathrm{loc}}(p_1)))\cdots)\]
is a horizontal submanifold of degree one in $\mathcal{B}^+_2$ containing $f^N(q)$. Since $f : \mathcal{B}^+_2\cap f^{-1}(\mathcal{B}^+_3)\to \mathcal{B}^+_3$ is a crossed mapping of degree two, $V^s_{31\overline{0}}(a, b)^+\cap V^u_{\overline{0}2\cdots 23}(a, b)^+$ contains exactly two points in $\mathcal{B}^+_3$ counted with multiplicity, one of which is $f^{N+1}(q)$. By (iii) of Proposition~\ref{PRP:pieces_positive} together with the assumption $\mathrm{card}(W_{31\overline{0}}^s(a, b)^+\cap W_{\overline{0}23}^u(a, b)^+)\geq 1$, we see that $W^s_{31\overline{0}}(a, b)^+\cap W^u_{\overline{0}2\cdots 23}(a, b)^+=V^s_{31\overline{0}}(a, b)^+\cap V^u_{\overline{0}2\cdots 23}(a, b)^+$ holds. Hence $f^{N+1}(q)\in\R^2$ and this implies $q\in\R^2$. It follows that $V^u(p_1)\cap V^s(p_1)\subset \R^2$, and so $\htop(f_{\R})=\log 2$ thanks to Theorem 10.1 of~\cite{BLS}.

Next consider the case $(a, b)\in \mathcal{F}^-_{\R}\cap\{b<0\}$. Choose any point $q\in V^u(p_3)\cap V^s(p_1)$ with $q\ne p_1$ and assume that $\mathrm{card}(W_{41\overline{0}}^s(a, b)^-\cap W^u_{\overline{43}4124}(a, b)^-_{\mathrm{inner}})\geq 1$ holds. As before, we may assume $q\in V^u_{\mathrm{loc}}(p_3)$. Recall that $V^u_{\mathrm{loc}}(p_3)$ is a degree one horizontal submanifold in $\mathcal{B}^-_3$ by Proposition~\ref{PRP:local_unstable}. Since $f : \mathcal{B}^-_3\cap f^{-1}(\mathcal{B}^-_4)\to \mathcal{B}^-_4$ is a crossed mapping of degree two, $f(V^u_{\mathrm{loc}}(p_3))\cap V^s_{41\overline{0}}(a, b)^-$ contains exactly two points, one of which is $f(q)$. Since the submanifolds $V^u_{\mathrm{loc}}(p_3)$ and $V^s_{41\overline{0}}(a, b)^-$ are real, we see that these two points belong to $\R^2$ by (iii) of Proposition~\ref{PRP:pieces_negative}. The rest of the argument stays the same as in the case $(a, b)\in \mathcal{F}^+_{\R}\cap\{b>0\}$, where Theorem 10.1 of~\cite{BLS} is replaced by Lemma~\ref{LMM:heteroclinic}.

To prove the converse, consider first the case $(a, b)\in \mathcal{F}^+_{\R}\cap\{b>0\}$ and assume that $W_{31\overline{0}}^s(a, b)^+\cap W_{\overline{0}23}^u(a, b)^+=\emptyset$ holds. Since $V_{31\overline{0}}^s(a, b)^+$ is a vertical submanifold of degree one in $\mathcal{B}^+_3$ and $V_{\overline{0}23}^u(a, b)^+$ is a horizontal submanifold of degree two in $\mathcal{B}^+_3$, the intersection $V_{31\overline{0}}^s(a, b)^+\cap V_{\overline{0}23}^u(a, b)^+$ consists of two points in $\mathcal{B}^+_3$ counted with multiplicity. By the assumption we see that the two points do not belong to $\R^2$, hence $V^u(p_1)\cap V^s(p_1)$ has elements outside $\R^2$. It follows from Theorem 10.1 of~\cite{BLS} that $\htop(f_{a, b}|_{\R^2})<\log 2$ holds. 

When $(a, b)\in \mathcal{F}^-_{\R}\cap\{b<0\}$, we must analyze the heteroclinic intersection $V^u(p_3)\cap V^s(p_1)$. However, thanks to Lemma~\ref{LMM:heteroclinic}, the above argument works in this case as well. This finishes the proof of Theorem~\ref{THM:maximalentropy} (Maximal Entropy). 
\end{proof}

A similar characterization for the H\'enon maps which are hyperbolic horseshoes on $\R^2$ in terms of the intersections of special pieces will be given in Theorem~\ref{THM:hyperbolichorseshoe} (Hyperbolic Horseshoes).

\subsection{Tin can argument}\label{subsection5.2}

As we have seen in Theorem~\ref{THM:maximalentropy} (and we will see in Theorem~\ref{THM:hyperbolichorseshoe}), the intersections of certain special pieces of $W^{u/s}(p_i)$ is responsible for a H\'enon map to be a hyperbolic horseshoe on $\R^2$ (and for a H\'enon map to attain the maximal entropy on $\R^2$). We are thus led to introduce the following complex tangency loci in the complexified parameter space $\mathcal{F}^{\pm}$. 

\begin{dfn}[\textbf{Complex Tangency Loci}]
We define 
\[\mathcal{T}^+\equiv\bigl\{(a, b)\in\mathcal{F}^+ : \mbox{$V_{31\overline{0}}^s(a, b)^+$ and $V_{\overline{0}23}^u(a, b)^+$ intersect tangentially}\bigr\}\]
and
\[\mathcal{T}^-\equiv\bigl\{(a, b)\in\mathcal{F}^- : \mbox{$V_{41\overline{0}}^s(a, b)^-$ and $V_{\overline{43}4124}^u(a, b)^-$ intersect tangentially}\bigr\},\]
and call them the \textit{complex tangency loci}.
\label{DFN:tangency_loci}
\end{dfn}

Let us write 
\[\partial^v \mathcal{F}^{\pm} \equiv \bigl\{(a, b) \in\C\times I^{\pm} : |a-a^{\pm}_{\mathrm{aprx}}(b)|= \chi^{\pm}(b) \bigr\}.\]
The purpose of this section is to show the following theorem.

\begin{thm}[\textbf{Tin Can}\footnote{A similar condition has been first introduced in~\cite{BS2} where $\partial^v\mathcal{F}^{\pm}$ is replaced by the vertical boundary of a bidisk which looks like a tin can.}]
We have (i) $\overline{\mathcal{T}^+}\cap \partial^v \mathcal{F}^+=\emptyset$ and (ii) $\overline{\mathcal{T}^-}\cap \partial^v \mathcal{F}^-=\emptyset$.
\label{THM:tincan}
\end{thm}

\begin{figure}
  \includegraphics[height=8cm]{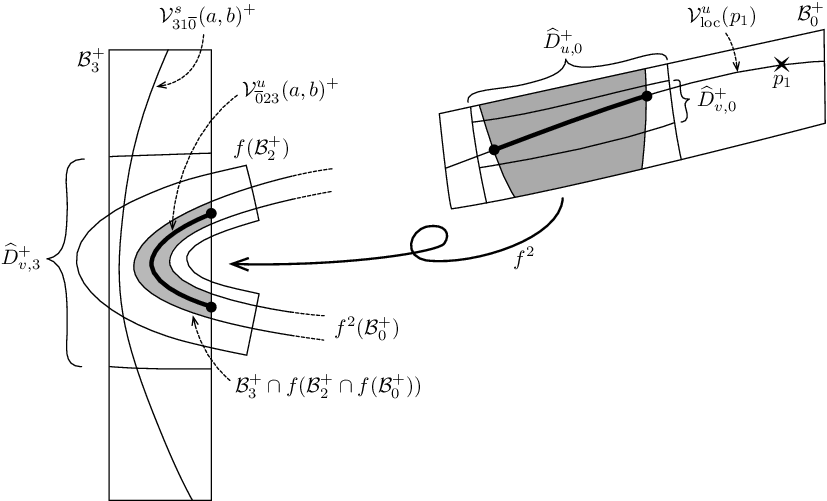}
  \caption{Figure of $\widehat{D}^+_{u, 0}$, $\widehat{D}^+_{v, 0}$ and $\widehat{D}^+_{v, 3}$.}
  \label{FIG:tincan_positive}
\end{figure}

\noindent
\textit{Proof of (i).} 
When we write $\mathcal{B}^+_3=D_{u, 3}^+\times_{\mathrm{pr}} D_{v, 3}^+$, one can choose\footnote{As seen in Figure~\ref{FIG:tincan_positive}, the piece $V^s_{31\overline{0}}(a, b)^+$ of the stable manifold $V^s(p_1)$ is ``curvy'' when $b$ is close to $1$. Hence, we choose a smaller $\widehat{D}_{v, 3}^+\subset D_{v, 3}^+$ so that $\pi^+_{u, 3}(\mathcal{V}^s_{31\overline{0}}(a, b)^+\cap\widehat{\mathcal{B}}^+_3)$ becomes smaller.} a smaller $\widehat{D}_{v, 3}^+\subset D_{v, 3}^+$ so that $\widehat{\mathcal{B}}^+_3\equiv D_{u, 3}^+\times_{\mathrm{pr}} \widehat{D}_{v, 3}^+$ contains $\mathcal{B}^+_3\cap f(\mathcal{B}^+_2)$. 

Let $\varphi : \C\to \C^2$ be a uniformization of $V^u(p_1)$ and let $\pi^+_{u, 3} : \widehat{\mathcal{B}}^+_3\to D_{u, 3}^+$ be the $u$-projection in $\mathcal{B}^+_3$. Denote by $\mathcal{C}(a, b)$ the set of critical points of $\pi^+_{u, 3}\circ \varphi : \varphi^{-1}(V_{\overline{0}23}^u(a, b)^+)\to D_{u, 3}^+$. To prove Theorem~\ref{THM:tincan} (Tin Can), it is sufficient to show 
\begin{equation}
\pi^+_{u, 3}\circ \varphi (\mathcal{C}(a, b))\cap \pi^+_{u, 3}(V_{31\overline{0}}^s(a, b)^+\cap\widehat{\mathcal{B}}^+_3)=\emptyset
\label{tincan1}
\end{equation}
for all $(a, b)\in \partial^v \mathcal{F}^+$. Note that the boxes $\mathcal{B}^+_i$ as well as the maps $\pi^+_{u, 3}$ and $\varphi$ depend continuously on $(a, b)\in \partial^v \mathcal{F}^+$. 

To achieve this, we introduce certain ``neighborhoods'' of the special pieces $V^s_{31\overline{0}}(a, b)^+$ and $V^u_{\overline{0}23}(a, b)^+$ as follows. Choose a large $N\geq 1$ and write 
\[\mathcal{V}^s_{\mathrm{loc}}(p_1)\equiv \mathcal{B}^+_0\cap f^{-1}(\mathcal{B}^+_0)\cap \cdots \cap f^{-N+1}(\mathcal{B}^+_0)\cap f^{-N}(\mathcal{B}^+_0).\]
Define
\[\mathcal{V}_{31\overline{0}}^s(a, b)^+\equiv \widehat{\mathcal{B}}^+_3\cap f^{-1}(\mathcal{B}^+_1\cap f^{-1}(\mathcal{V}^s_{\mathrm{loc}}(p_1))).\]
Similarly, choose a large $M\geq 1$ and write 
\[\mathcal{V}^u_{\mathrm{loc}}(p_1)\equiv \mathcal{B}^+_0\cap f(\mathcal{B}^+_0)\cap \cdots \cap f^{M-1}(\mathcal{B}^+_0)\cap f^M(\mathcal{B}^+_0).\]
Take smaller $\widehat{D}_{u, 0}^+\subset D_{u, 0}^+$ and $\widehat{D}_{v, 0}^+\subset D_{v, 0}^+$ so that\footnote{First take smaller $\widehat{D}_{u, 0}^+\subset D_{u, 0}^+$ so that $\mathcal{B}^+_3\cap f(\mathcal{B}^+_2\cap f(\widehat{D}_{u, 0}^+\times_{\mathrm{pr}} D_{v, 0}^+))$ contains $\mathcal{B}^+_3\cap f(\mathcal{B}^+_2\cap f(\mathcal{B}^+_0))$, and second take a smaller $\widehat{D}_{v, 0}^+\subset D_{v, 0}^+$ so that $\mathcal{V}_{\overline{0}23}^u(a, b)^+$ contains $\mathcal{B}^+_3\cap f(\mathcal{B}^+_2\cap f(\mathcal{V}^u_{\mathrm{loc}}(p_1)))$ (see Figure~\ref{FIG:tincan_positive} again).}
\[\mathcal{V}_{\overline{0}23}^u(a, b)^+\equiv \mathcal{B}^+_3\cap f(\mathcal{B}^+_2\cap f(\widehat{D}_{u, 0}^+\times_{\mathrm{pr}} \widehat{D}_{v, 0}^+))\]
contains $\mathcal{B}^+_3\cap f(\mathcal{B}^+_2\cap f(\mathcal{V}^u_{\mathrm{loc}}(p_1)))$. 

The above construction immediately implies

\begin{lmm}
We have $V_{31\overline{0}}^s(a, b)^+ \cap\widehat{\mathcal{B}}^+_3\subset \mathcal{V}_{31\overline{0}}^s(a, b)^+$ and $V_{\overline{0}23}^u(a, b)^+ \subset \mathcal{V}_{\overline{0}23}^u(a, b)^+$.
\label{LMM:neighborhood_positive}
\end{lmm}

The following claim can be verified by using rigorous numerics and its proof will be supplied in Subsection~\ref{subsection6.4}.

\begin{lmm}
Let $(a, b)\in \partial^v \mathcal{F}^+$. Then, for every fixed $v_0\in \widehat{D}_{v, 0}^+$ we have 
\[\frac{d}{d u}\left\{\pi^+_{u, 3} \circ f^2\circ \iota_{v_0}(u)\right\}\ne 0\]
for $u\in D_{u, 0}^+$ with $\iota_{v_0}(u)\in \mathcal{B}^+_0\cap f^{-1}(\mathcal{B}^+_2\cap f^{-1}(\pi^+_{u, 3}(\mathcal{V}^s_{31\overline{0}}(a, b)^+)\times_{\mathrm{pr}}D_{v, 3}^+))$.
\label{LMM:tangency_positive}
\end{lmm}

Lemmas~\ref{LMM:neighborhood_positive} and \ref{LMM:tangency_positive} yield Eqn.~(\ref{tincan1}), which finishes the proof of (i).

\begin{figure}
  \includegraphics[height=8.5cm]{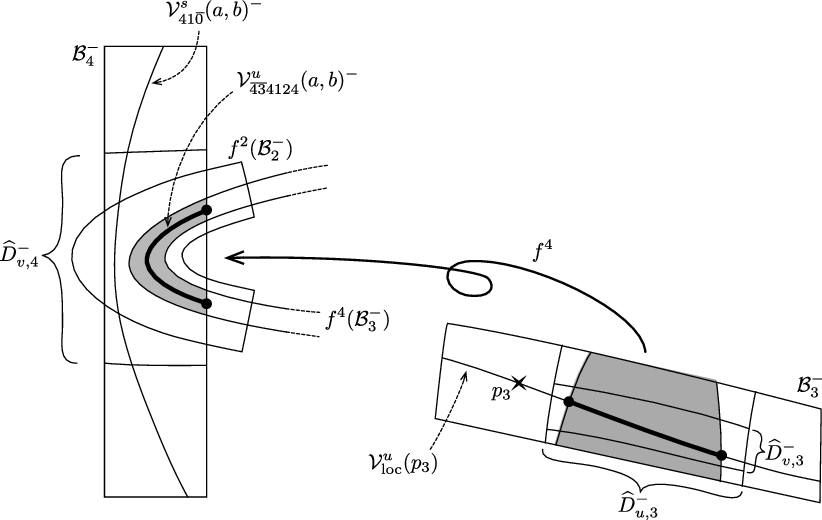}
  \caption{Figure of $\widehat{D}^-_{u, 3}$, $\widehat{D}^-_{v, 3}$ and $\widehat{D}^-_{v, 4}$.}
\label{FIG:tincan_negative}
\end{figure}

\begin{proof}[Proof of (ii)]
As in the previous case, one can choose a smaller $\widehat{D}_{v, 4}^-\subset D_{v, 4}^-$ so that $\widehat{\mathcal{B}}^-_4\equiv D_{u, 4}^-\times_{\mathrm{pr}} \widehat{D}_{v, 4}^-$ contains $\mathcal{B}^-_4\cap f(\mathcal{B}^-_2)$ (see Figure~\ref{FIG:tincan_negative}).

Let $\varphi : \C\to \C^2$ be a uniformization of $V^u(p_3)$ and let $\pi^-_{u, 4} : \widehat{\mathcal{B}}^-_4\to D_{u, 4}^-$ be the vertical projection in $\mathcal{B}^-_4$. Denote by $\mathcal{C}(a, b)$ the set of critical points of $\pi^-_{u, 4}\circ \varphi : \varphi^{-1}(V_{\overline{43}4124}^u(a, b)^-)\to D_{u, 4}^-$. To prove the theorem, it is sufficient to show 
\begin{equation}
\pi^-_{u, 4}\circ \varphi (\mathcal{C}(a, b))\cap \pi^-_{u, 4}(V_{41\overline{0}}^s(a, b)^-)=\emptyset
\label{tincan2}
\end{equation}
for all $(a, b)\in \partial^v \mathcal{F}^-$. 

Choose a large $N\geq 1$ and write 
\[\mathcal{V}^s_{\mathrm{loc}}(p_1)\equiv \mathcal{B}^-_0\cap f^{-1}(\mathcal{B}^-_0)\cap \cdots \cap f^{-N+1}(\mathcal{B}^-_0)\cap f^{-N}(\mathcal{B}^-_0).\]
Define
\[\mathcal{V}_{41\overline{0}}^s(a, b)^-\equiv \widehat{\mathcal{B}}^-_4\cap f^{-1}(\mathcal{B}^-_1\cap f^{-1}(\mathcal{V}^s_{\mathrm{loc}}(p_1))).\]
Recall that $f^2 : \mathcal{B}^-_3\cap f^{-1}(\mathcal{B}^-_4\cap f^{-1}(\mathcal{B}^-_3))\to\mathcal{B}^-_3$ is a crossed mapping of degree two of horseshoe type by Lemma~\ref{LMM:OCC_negative1}. Let $\mathcal{V}^0\equiv\mathcal{B}^-_3$ and define inductively 
\[\mathcal{V}^n\equiv\mathcal{B}^-_3\cap_{p_3} f(\mathcal{B}^-_4\cap f(\mathcal{V}^{n-1})),\]
where $\mathcal{B}^-_3\cap_{p_3} f(\mathcal{B}^-_4\cap f(\mathcal{V}^{n-1}))$ means the connected component of $\mathcal{B}^-_3\cap f(\mathcal{B}^-_4\cap f(\mathcal{V}^{n-1}))$ containing the fixed point $p_3$. Let us choose a large $M\geq 1$ and write $\mathcal{V}^u_{\mathrm{loc}}(p_3)\equiv\mathcal{V}^M$. We take smaller $\widehat{D}_{u, 3}^-\subset D_{u, 3}^-$ and $\widehat{D}_{v, 3}^-\subset D_{v, 3}^-$ so that
\[\mathcal{V}_{\overline{43}4124}^u(a, b)^-\equiv \mathcal{B}^-_4\cap f(\mathcal{B}^-_2\cap f(\mathcal{B}^-_1\cap f(\mathcal{B}^-_4\cap f(\widehat{D}_{u, 3}^-\times_{\mathrm{pr}} \widehat{D}_{v, 3}^-))))\]
contains $\mathcal{B}^-_4\cap f(\mathcal{B}^-_2\cap f(\mathcal{B}^-_1\cap f(\mathcal{B}^-_4\cap f(\mathcal{V}^u_{\mathrm{loc}}(p_3)))))$ (see Figure~\ref{FIG:tincan_negative} again). Then, as in the previous case, 

\begin{lmm}
We have $V_{41\overline{0}}^s(a, b)^- \cap\widehat{\mathcal{B}}^-_4\subset \mathcal{V}_{41\overline{0}}^s(a, b)^-$ and $V_{\overline{43}4124}^u(a, b)^- \subset \mathcal{V}_{\overline{43}4124}^u(a, b)^-$.
\label{LMM:neighborhood_negative}
\end{lmm}

\begin{proof}
Recall the proof of Proposition~\ref{PRP:local_unstable}. It is easy to see that the horizontal submanifold $D_n$ in the proof is contained in $\mathcal{V}^n$ above, so the conclusion follows. 
\end{proof}

The following claim can be verified by using rigorous numerics and its proof will be supplied in Subsection~\ref{subsection6.4}.

\begin{lmm}
Let $(a, b)\in \partial^v \mathcal{F}^-$. Then, for every fixed $v_0\in \widehat{D}_{v, 3}^-$ we have 
\[\frac{d}{d u}\left\{\pi^-_{u, 4} \circ f^4\circ \iota_{v_0}(u)\right\}\ne 0\]
for $u\in D_{u, 3}^-$ with $\iota_{v_0}(u)\in\mathcal{B}^-_3\cap f^{-1}(\mathcal{B}^-_4\cap f^{-1}(\mathcal{B}^-_1\cap f^{-1}(\mathcal{B}^-_2\cap f^{-1}(\pi^-_{u, 4}(\mathcal{V}^s_{41\overline{0}}(a, b)^-)\times_{\mathrm{pr}}D_{v, 4}^-))))$.
\label{LMM:tangency_negative}
\end{lmm}

Lemmas~\ref{LMM:neighborhood_negative} and \ref{LMM:tangency_negative} yield Eqn.~(\ref{tincan2}), which finishes the proof of (ii). 
\end{proof}

\subsection{Tangency loci}\label{subsection5.3}

In this subsection another definition of the special pieces is given to analyze the local complex analytic property of the tangency loci $\mathcal{T}^{\pm}$. Below we let $p_3\equiv (z_3, z_3)$ be the unique fixed point of $f_{a, b}$ in $\mathcal{B}^-_3\cap \mathcal{B}^-_4$ for $(a, b)\in\mathcal{F}^-$. The following construction can be adapted to the other fixed point $p_1\equiv (z_1, z_1)\in\mathcal{B}^{\pm}_0$ of $f_{a, b}$ for $(a, b)\in\mathcal{F}^{\pm}$ as well.

We first examine the case $b\ne 0$. Let $(a, b)\in\mathcal{F}^-\cap\{b\ne 0\}$. Let $\Psi_{a, b} : \C\to \C^2$ be the uniformization of $V^u(p_3)$ with $\Psi_{a, b}(0)=p_3$ and $(\pi_x\circ \Psi_{a, b})'(0)=1$. By the functional equation $\Psi_{a, b}(\lambda z)=f_{a, b}(\Psi_{a, b}(z))$ we see that $\Psi_{a, b}$ is of the form $\Psi_{a, b}(z)=(\varphi_{a, b}(z), \varphi_{a, b}(z/\lambda))$, where $\lambda$ is the unstable eigenvalue of $Df_{a, b}$ at $p_3$. Let $V_{\mathrm{loc}}^u(p_3)$ be the connected component of $V^u(p_3)\cap\mathcal{B}^-_3$ containing $p_3$ and set $\mathrm{\Omega}_{\mathrm{loc}}(a, b)\equiv\Psi_{a, b}^{-1}(V_{\mathrm{loc}}^u(p_3))\subset\C$. We generalize this definition to any backward admissible sequence of the form $J=\overline{43}j_{-n}\cdots j_{-1}j_0\in\mathfrak{S}_{\mathrm{bwd}}^-$ as
\[\mathrm{\Omega}_J(a, b)\equiv \lambda^{n+1}\mathrm{\Omega}_{\mathrm{loc}}(a, b)\cap \Psi_{a, b}^{-1}(\mathcal{B}^-_{j_0}\cap f_{a, b}(\mathcal{B}^-_{j_{-1}}\cap \cdots \cap f_{a, b}(\mathcal{B}^-_{j_{-n+1}}\cap f_{a, b}(\mathcal{B}^-_{j_{-n}}))\cdots )).\]

\begin{lmm}
For $(a, b)\in\mathcal{F}^-\cap\{b\ne 0\}$, $\mathrm{\Omega}_{\overline{43}412}(a, b)=\lambda^3\mathrm{\Omega}_{\mathrm{loc}}(a, b)\cap \Psi_{a, b}^{-1}(\mathcal{B}^-_2 \cap f_{a, b}^2(\mathcal{B}^-_1\cap f_{a, b}^3(\mathcal{B}^-_4)))$ consists of two connected components with disjoint closures. 
\label{LMM:degenerate_components}
\end{lmm}

\begin{proof}
Since one can verify 
\begin{align*}
\Psi_{a, b}(\mathrm{\Omega}_{\overline{43}412}(a, b)) = & \Psi_{a, b}(\lambda^3\mathrm{\Omega}_{\mathrm {loc}}(a, b)\cap \Psi_{a, b}^{-1}(\mathcal{B}^-_2\cap f_{a, b}(\mathcal{B}^-_1\cap f_{a, b}(\mathcal{B}^-_4)))) \\
= & f^3_{a, b}(V_{\mathrm{loc}}^u(p_3))\cap \mathcal{B}^-_2\cap f_{a, b}(\mathcal{B}^-_1\cap f_{a, b}(\mathcal{B}^-_4)) \\
= & V^u_{\overline{43}412}(a, b)^-
\end{align*}
and since $\Psi_{a, b}$ is injective, Proposition~\ref{PRP:43412} yields that $\mathrm{\Omega}_{\overline{43}412}(a, b)$ has two connected components with disjoint closures. 
\end{proof}

We next examine the case $b=0$. Let $(a, 0)\in\mathcal{F}^-\cap\{b=0\}$. Let $\varphi_a : \C\to\C$ be the linearization of $p_a(z)=z^2-a$ with $\varphi_a(0)=z_3$ and $\varphi'_a(0)=1$. Since it satisfies $p_a(\varphi_a(z))=\varphi_a(\lambda z)$ where $\lambda\equiv p'_a(z_3)$, the map $\Psi_{a, 0} : \C\to\mathrm{\Gamma}_a\equiv\{(x, y)\in\C^2 : x=y^2-a\}$ defined by $\Psi_{a, 0}(z)\equiv (\varphi_a(z), \varphi_a(z/\lambda))$ satisfies the functional equation $\Psi_{a, 0}(\lambda z)=f_{a, 0}(\Psi_{a, 0}(z))$. Let $V_{\mathrm{loc}}^u(p_3)$ be the connected component of $\mathrm{\Gamma}_a\cap\mathcal{B}^-_3$ containing $p_3$ and let $\mathrm{\Omega}_{\mathrm{loc}}(a, 0)\subset\C$ be the connected component of $\Psi_{a, 0}^{-1}(V_{\mathrm{loc}}^u(p_3))$ containing the origin. Note that $\Psi_{a, 0}(\mathrm{\Omega}_{\mathrm{loc}}(a, 0))=V_{\mathrm{loc}}^u(p_3)$ holds. We generalize this definition to any backward admissible sequence of the form $J=\overline{43}j_{-n}\cdots j_{-1}j_0\in\mathfrak{S}_{\mathrm{bwd}}^-$ as 
\begin{align*}
\mathrm{\Omega}_J(a, 0) \equiv & \lambda^{n+1}\mathrm{\Omega}_{\mathrm{loc}}(a, 0)\cap \varphi_a^{-1}(D^-_{x, j_0}\cap p_a(D^-_{x, j_{-1}}\cap \cdots \cap p_a(D^-_{x, j_{-n+1}}\cap p_a(D^-_{x, j_{-n}}))\cdots )) \\
= & 
\lambda^{n+1}\mathrm{\Omega}_{\mathrm{loc}}(a, 0)\cap \Psi_{a, 0}^{-1}(\mathcal{B}^-_{j_0}\cap f_{a, 0}(\mathcal{B}^-_{j_{-1}}\cap\cdots \cap f_{a, 0}(\mathcal{B}^-_{j_{-n+1}}\cap f_{a, 0}(\mathcal{B}^-_{j_{-n}}))\cdots)),
\end{align*}
where we write $\mathcal{B}^-_i=D_{x, i}^-\times D_{y, i}^-$ with respect to the standard Euclidean coordinates. As before, one can verify $\Psi_{a, 0}(\mathrm{\Omega}_{\overline{43}412}(a, 0))=V^u_{\overline{43}412}(a, 0)^-$, but $\Psi_{a, 0}$ is not injective anymore. Hence we have to show

\begin{lmm}
For $(a, 0)\in\mathcal{F}^-\cap\{b=0\}$, $\mathrm{\Omega}_{\overline{43}412}(a, 0)=\lambda^3\mathrm{\Omega}_{\mathrm{loc}}(a, 0)\cap \Psi_{a, 0}^{-1}(\mathcal{B}^-_2 \cap f_{a, 0}(\mathcal{B}^-_1\cap f_{a, 0}(\mathcal{B}^-_4)))$ consists of two connected components with disjoint closures. 
\label{LMM:two_components}
\end{lmm}

\begin{proof}
Below, we essentially follow the proof of Lemma 4.4 in~\cite{BS2}. First recall that the crossed mapping $f_{a, 0}^2 : \mathcal{B}^-_3\cap f_{a, 0}^{-1}(\mathcal{B}^-_4\cap f_{a, 0}^{-1}(\mathcal{B}^-_3))\to\mathcal{B}^-_3$ of degree two satisfies the (OCC) by Lemma~\ref{LMM:OCC_negative1}. This means that $p^2_a : D_{x, 3}^-\cap p_a^{-1}(D_{x, 4}^-\cap p_a^{-1}(D_{x, 3}^-))\to D_{x, 3}^-$ is a covering of degree two, so $D_{x, 3}^-\cap p_a^{-1}(D_{x, 4}^-\cap p_a^{-1}(D_{x, 3}^-))$ consists of two disjoint submanifolds. Let $D$ be the one containing the fixed point $z_3$. Then, $D$ is compactly contained in $D_{x, 3}^-$ and $p_a^2 : D\to D_{x, 3}^-$ is a conformal equivalence. So we may define $\lim_{n\to\infty}\lambda^{2n}(p_a^2|_D)^{-n} : D_{x, 3}^-\to\C$, which is the inverse of $\varphi_a$. It follows that $\varphi_a : \mathrm{\Omega}_{\mathrm{loc}}(a, 0)\to D_{x, 3}^-$ is a univalent function. Secondly we compute as
\begin{align}
\varphi_a'(\lambda^n z)\cdot \lambda^n=(p_a^n\circ\varphi_a)'(z)
& =p_a'(p_a^{n-1}\circ\varphi_a(z))\cdots p_a'(p_a\circ\varphi_a(z))p_a'(\varphi_a(z))\varphi_a'(z) \label{computation_1} \\
& =p_a'(\varphi_a(z/\lambda^{n-1}))\cdots p_a'(\varphi_a(z/\lambda))p_a'(\varphi_a(z))\varphi_a'(z). \label{computation_2}
\end{align}
This result will be useful in the discussion below.

Let $c\in \mathrm{\Omega}_{\mathrm{loc}}(a, 0)$ be the unique point so that $\varphi_a(c)=0$ holds. Then, by Eqn.~(\ref{computation_1}) we have $\varphi_a'(\lambda c)\cdot \lambda=p_a'(0)\varphi_a'(c)=0$, hence $\varphi_a'(\lambda c)=0$. Conversely, if $z\in\lambda\mathrm{\Omega}_{\mathrm{loc}}(a, 0)$ and $\varphi_a'(z)=0$, then again by Eqn.~(\ref{computation_2}) we have $0=\varphi_a'(z)\cdot \lambda=p_a'(\varphi_a(z/\lambda))\varphi_a'(z/\lambda)$. Since $\varphi_a$ is univalent on $\mathrm{\Omega}_{\mathrm{loc}}(a, 0)$, one sees $p_a'(\varphi_a(z/\lambda))=0$, hence $\varphi_a(z/\lambda)=0$ and $z=\lambda c$. It follows that $z=\lambda c$ is the unique critical point of $\varphi_a$ in $\lambda\mathrm{\Omega}_{\mathrm{loc}}(a, 0)$. This implies that $\Psi_{a, 0}(z)=(\varphi_a(z), \varphi_a(z/\lambda))$ has no critical point in $\lambda\mathrm{\Omega}_{\mathrm{loc}}(a, 0)$. Since $\Psi_{a, 0}(\lambda\mathrm{\Omega}_{\mathrm{loc}}(a, 0)\cap \Psi_{a, 0}^{-1}(\mathcal{B}^-_4))=f_{a, 0}(\Psi_{a, 0}(\mathrm{\Omega}_{\mathrm{loc}}(a, 0)))\cap \mathcal{B}^-_4=f_{a, 0}(V_{\mathrm{loc}}^u(p_3)) \cap \mathcal{B}^-_4=\mathrm{\Gamma}_a\cap\mathcal{B}^-_4$ is simply connected, it follows that $\Psi_{a, 0} : \lambda\mathrm{\Omega}_{\mathrm{loc}}(a, 0)\cap \Psi_{a, 0}^{-1}(\mathcal{B}^-_4)\to \mathrm{\Gamma}_a\cap\mathcal{B}^-_4$ is univalent. In particular, we see that $\Psi_{a, 0} : \lambda\mathrm{\Omega}_{\mathrm{loc}}(a, 0)\cap \Psi_{a, 0}^{-1}(\mathcal{B}^-_4\cap f^{-1}_{a, 0}(\mathcal{B}^-_1\cap f^{-1}_{a, 0}(\mathcal{B}^-_2)))\to \mathrm{\Gamma}_a\cap(\mathcal{B}^-_4\cap f^{-1}_{a, 0}(\mathcal{B}^-_1\cap f^{-1}_{a, 0}(\mathcal{B}^-_2)))$ is univalent.

The above calculation Eqn.~(\ref{computation_1}) also shows $\varphi_a'(\lambda^2 c)\cdot \lambda^2=p_a'(p_a(0))p_a'(0)\varphi_a'(c)=0$ and $\varphi_a'(\lambda^3 c)\cdot \lambda^3=p_a'(p_a^2(0))p_a'(p_a(0))p_a'(0)\varphi_a'(c)=0$, hence one has $\varphi_a'(\lambda^2 c)=0$ and $\varphi_a'(\lambda^3 c)=0$. Conversely, if we assume $z\in \lambda^3\mathrm{\Omega}_{\mathrm{loc}}(a, 0)$ and $\varphi_a'(z)=0$, then once again by the above computation Eqn.~(\ref{computation_2}) one sees $0 = \varphi_a'(z) \cdot \lambda^3 = p_a'(\varphi_a(z/\lambda)) p_a' (\varphi_a(z/\lambda^2)) p_a'(\varphi_a(z/\lambda^3)) \varphi_a'(z/\lambda^3)$. This implies $z=\lambda^2c, \lambda^3 c$, and hence $z=\lambda^2c, \lambda^3c$ are the only critical points of $\Psi_{a, 0}$ in $\lambda^3\mathrm{\Omega}_{\mathrm{loc}}(a, 0)$. Now, $\Psi_{a, 0}(\lambda^2c)=(\varphi_a(\lambda^2c), \varphi_a(\lambda c))=(p_a^2(0), p_a(0))$ does not belong to $\mathrm{\Gamma}_a\cap \mathcal{B}^-_3$ by Lemma~\ref{LMM:OCC_negative1} and $\Psi_{a, 0}(\lambda^3c)=(\varphi_a(\lambda^3c), \varphi_a(\lambda^2 c))=(p_a^3(0), p_a^2(0))$ does not belong to $\mathrm{\Gamma}_a\cap \mathcal{B}^-_2$ by Lemma~\ref{LMM:OCC_negative2}. It then follows that $\Psi_{a, 0}$ does not have critical points in $\lambda^3\mathrm{\Omega}_{\mathrm{loc}}(a, 0)\cap \Psi_{a, 0}^{-1}(\mathcal{B}^-_2)$ and hence not in the closure of $\mathrm{\Omega}_{\overline{43}412}(a, 0)$. 

By Proposition~\ref{PRP:43412}, $V^u_{\overline{43}412}(a, 0)^-\equiv\mathcal{B}^-_2\cap f_{a, 0}(\mathcal{B}^-_1\cap f_{a, 0}(\mathcal{B}^-_4\cap f_{a, 0}(V_{\mathrm{loc}}^u(p_3))))$ is a horizontal submanifold of degree one in $\mathcal{B}^-_2$. Recall that $f^2_{a, 0} : \mathrm{\Gamma}_a\cap (\mathcal{B}^-_4\cap f^{-1}_{a, 0}(\mathcal{B}^-_1\cap f^{-1}_{a, 0}(\mathcal{B}^-_2)))\to V^u_{\overline{43}412}(a, 0)^-$ is a covering map of degree two thanks to Lemma~\ref{LMM:OCC_negative2}. Since one can check that $\lambda^2(\lambda\mathrm{\Omega}_{\mathrm{loc}}(a, 0)\cap \Psi_{a, 0}^{-1}(\mathcal{B}^-_4\cap f^{-1}_{a, 0}(\mathcal{B}^-_1\cap f^{-1}_{a, 0}(\mathcal{B}^-_2))))=\mathrm{\Omega}_{\overline{43}412}(a, 0)$, it follows that $\Psi_{a, 0} : \mathrm{\Omega}_{\overline{43}412}(a, 0)\to V^u_{\overline{43}412}(a, 0)^-$ is a covering of degree two. In particular, $\mathrm{\Omega}_{\overline{43}412}(a, 0)$ consists of two submanifolds with disjoint closures and each of them is conformally equivalent to $V^u_{\overline{43}412}(a, 0)^-$ by $\Psi_{a, 0}$. Thus we are done. 
\end{proof}

Since $\Psi_{a, b}$ converges to $\Psi_{a, 0}$ as $b\to 0$ uniformly on compact sets, we see that $V^u_{\overline{43}412}(a, b)^-$ converges to $V^u_{\overline{43}412}(a, 0)^-$ as $b\to 0$ with respect to the Hausdorff topology. 

\begin{prp}
We have the following properties of $\mathcal{T}^{\pm}$.
\begin{enumerate}
\renewcommand{\labelenumi}{(\roman{enumi})}
\item $\mathcal{T}^{\pm}$ is a complex subvariety of $\mathcal{F}^{\pm}$.
\item $\mathcal{T}^-$ is reducible, i.e. one can write $\mathcal{T}^-=\mathcal{T}^-_1 \cup \mathcal{T}^-_2$ where $\mathcal{T}^-_i$ is a complex subvariety of $\mathcal{F}^-$ for $i=1, 2$.
\item The projection to the $b$-axis: 
\[\mathrm{pr}^+ : \mathcal{T}^+ \longrightarrow I^+\]
is a proper map of degree one. Similarly, the projection to the $b$-axis: 
\[\mathrm{pr}^- : \mathcal{T}^-_i \longrightarrow I^-\]
is a proper map of degree one for $i=1, 2$.
\item $\mathcal{T}^+$ (resp. $\mathcal{T}^-_i$) is a complex submanifold of $\mathcal{F}^+$ (resp. $\mathcal{F}^-$).
\end{enumerate}
\label{PRP:complex_loci}
\end{prp}

Note that for the complex locus $\mathcal{T}^-$, we can not a priori ``distinguish'' $\mathcal{T}^-_1$ and $\mathcal{T}^-_2$.

\begin{proof}
Below we first show (i), (ii) and (iii) for $b\ne 0$, and then prove all the claims for the general case.

(i) Proposition~\ref{PRP:T_is_subvariety} yields that $\mathcal{T}^{\pm}\cap\{b\ne 0\}$ is a subvariety in $\mathcal{F}^{\pm}\cap\{b\ne 0\}$. 

(ii) For $(a, b)\in \mathcal{F}^-$, let $\mathrm{\Omega}(a, b)'$ and $\mathrm{\Omega}(a, b)''$ be the two connected components of $\mathrm{\Omega}_{\overline{43}412}(a, b)$ as in Lemmas~\ref{LMM:degenerate_components} and \ref{LMM:two_components}. These define a splitting of $V^u_{\overline{43}412}(a, b)^-$ into two parts $\Psi_{a, b}(\mathrm{\Omega}(a, b)')$ and $\Psi_{a, b}(\mathrm{\Omega}(a, b)'')$ (they coincide when $b=0$). Hence by letting $\mathcal{T}^-_1$ to be the parameter locus where $\mathcal{B}^-_4\cap f_{a, b}(\Psi_{a, b}(\mathrm{\Omega}(a, b)'))$ intersects $V^s_{41\overline{0}}(a, b)^-$ tangentially and $\mathcal{T}^-_2$ the parameter locus where $\mathcal{B}^-_4\cap f_{a, b}(\Psi_{a, b}(\mathrm{\Omega}(a, b)''))$ intersects $V^s_{41\overline{0}}(a, b)^-$ tangentially, the locus $\mathcal{T}^-$ can be written as $\mathcal{T}^-=\mathcal{T}^-_1\cup \mathcal{T}^-_2$. Moreover, Proposition~\ref{PRP:T_is_subvariety} yields that $\mathcal{T}^-_i\cap \{b\ne 0\}$ is a complex subvariety in $\mathcal{F}^-\cap \{b\ne 0\}$ for $i=1, 2$.

(iii) Thanks to Theorem~\ref{THM:tincan} (Tin Can), the condition $\overline{A}\cap (\partial D\times E)=\emptyset$ in Lemma~\ref{LMM:proper} is satisfied. Hence it follows that $\mathrm{pr}^+ : \mathcal{T}^+\cap \{b\ne 0\}\rightarrow I^+\cap \{b\ne 0\}$ is a proper map. Since $\mathcal{T}^+$ is non-empty, its degree is at least one. Below we prove that the degree is at most one.  

For this, we consider the quadratic family in one variable $p_a(x)=x^2-a$. Its critical value is $c(a)=-a$. One of the fixed points of $p_a$ is $q(a)=(1+\sqrt{1+4a})/2$. Let $\tilde{q}(a)=-(1+\sqrt{1+4a})/2$, which satisfies $\tilde{q}(a)\ne q(a)$ and $p_a(\tilde{q}(a))=q(a)$. For all $a_0>0$, an easy computation shows
\[\frac{d}{da}(\tilde{q}-c)(a_0)<0.\]

Let $U^s$ and $U^u$ be open sets in $\C$ containing $\alpha\in\C$, and let $\varphi^s_{a, b} : U^s\to \C^2$ and $\varphi^u_{a, b} : U^u\to \C^2$ be the uniformization of the special pieces $V^s_{31\overline{0}}(a, b)^+$ and $V^u_{\overline{0}23}(a, b)^+$ respectively so that $\varphi^s_{2, 0}(\alpha)=\varphi^u_{2, 0}(\alpha)$ is the unique tangency for $b=0$. Since $\pi_x\circ\varphi^s_{a, 0}(\alpha)=\tilde{q}(a)$ and $\pi_x\circ\varphi^u_{a, 0}(\alpha)=c(a)$ hold, the previous computation implies that 
\[\frac{\partial}{\partial a}\left\{\pi_x\circ\varphi^s_{a, b}(z)-\pi_x\circ\varphi^u_{a, b}(z)\right\}\]
has negative real part for any $z\in \C$ close to $\alpha$ and any $b\in I^+\cap \{b\ne 0\}$ close to zero. This yields that $V^u_{\overline{0}23}(a, b)^+$ makes a tangency with $V^s_{31\overline{0}}(a, b)^+$ at most once when $b$ is fixed near $0$ and $a$ changes. It follows that the degree of $\mathrm{pr}^+ : \mathcal{T}^+\cap \{b\ne 0\}\to I^+\cap \{b\ne 0\}$ is one. The proof for $\mathrm{pr}^- : \mathcal{T}^-_i\cap \{b\ne 0\} \to I^-\cap \{b\ne 0\}$ is similar. This proves (iii) for the case $b\ne 0$.

Now we prove the general case. Since $\mathrm{pr}^+ : \mathcal{T}^+\cap \{b\ne 0\}\to I^+\cap \{b\ne 0\}$ is degree one, it follows from Proposition~\ref{PRP:submanifold} that $\mathcal{T}^+\cap \{b\ne 0\}$ is a complex submanifold of $\mathcal{F}^+\cap \{b\ne 0\}$. Hence, there exists a holomorphic function: 
\[\kappa^+ \, : \, I^+\cap \{b\ne 0\} \longrightarrow \R\]
whose graph coincides with $\mathcal{T}^+\cap \{b\ne 0\}$. Theorem~\ref{THM:tincan} (Tin Can) tells that $\kappa^+$ is locally bounded near $b=0$, hence $b=0$ is a removable singularity of $\kappa^+$. By letting $\kappa^+(0)=2$, we obtain a holomorphic function $\kappa^+$ defined on all of $I^+$ to the $a$-axis whose graph coincides with $\mathcal{T}^+$. It follows that $\mathrm{pr}^+ : \mathcal{T}^+\to I^+$ is proper of degree one and hence $\mathcal{T}^+$ is a complex submanifold of $\mathcal{F}^+$. Similarly we obtain a holomorphic function $\kappa^-_i$ defined on all of $I^-$ to the $a$-axis whose graph coincides with $\mathcal{T}^-_i$. It follows that $\mathrm{pr}^- : \mathcal{T}^-_i\to I^-$ is proper of degree one and hence $\mathcal{T}^-_i$ is a complex submanifold of $\mathcal{F}^-$. This proves all the claims for general case. 
\end{proof}

\subsection{End of the proof}\label{subsection5.4}

In this subsection we investigate the real sections of the tangency loci $\mathcal{T}_{\R}^{\pm}$ and apply it to the proof of the Main Theorem. As a consequence of its proof a characterization is obtained for the H\'enon maps which are hyperbolic horseshoes on $\R^2$ in terms of the special intersections.

Let us first investigate the real locus $\mathcal{T}_{\R}^+\equiv\mathcal{T}^+\cap \mathcal{F}^+_{\R}$.

\begin{prp}
The following properties hold for $\mathcal{T}_{\R}^+$.
\begin{enumerate}
\renewcommand{\labelenumi}{(\roman{enumi})}
\item We have $(a, b)\in \mathcal{T}_{\R}^+$ iff $(a, b)\in\mathcal{F}^+_{\R}$ and $W_{31\overline{0}}^s(a, b)^+$ intersects $W_{\overline{0}23}^u(a, b)^+$ tangentially in $\R^2$.
\item There exists a real analytic function:
\[\kappa_{\R}^+ : (-\varepsilon, 1+\varepsilon)\longrightarrow \R\]
so that $\mathcal{T}_{\R}^+$ coincides with the graph of $\kappa_{\R}^+$.
\end{enumerate}
\label{PRP:real_loci_positive}
\end{prp}

\begin{proof}
(i) If $(a, b)\in \mathcal{T}_{\R}^+$, then $(a, b)\in\mathcal{F}^+_{\R}$ and $V_{31\overline{0}}^s(a, b)^+$ intersects $V_{\overline{0}23}^u(a, b)^+$ tangentially in $\C^2$. If this tangential intersection is not real, then its complex conjugate is also a distinct tangential intersection. This contradicts to the fact that the intersection $V_{31\overline{0}}^s(a, b)^+\cap V_{\overline{0}23}^u(a, b)^+$ consists of two points counted with multiplicity. The converse is obvious.

(ii) Let $\overline{z}\in\C$ denote the complex conjugate of $z\in\C$. We first remark that the complex conjugate of a special piece $V_{\ast}^{u/s}(a, b)$ under $(x, y)\mapsto (\overline{x}, \overline{y})$ in $\C^2$ is $V_{\ast}^{u/s}(\overline{a}, \overline{b})$. Therefore, the tangency loci are invariant under the complex conjugation $(a, b)\mapsto (\overline{a}, \overline{b})$ in $\C^2$.

Take $b_{\ast}\in (-\varepsilon, 1+\varepsilon)$ and consider $(a_{\ast}, b_{\ast})\equiv (\mathrm{pr}^+)^{-1}(b_{\ast})\in\mathcal{T}^+$. If it does not belong to $\mathcal{T}^+_{\R}$, then its complex conjugate belongs to $\mathcal{T}^+$ but different from $(a_{\ast}, b_{\ast})$, and both are mapped to $b_{\ast}$ by $\mathrm{pr}^+$, contradicting to (iii) of Proposition~\ref{PRP:complex_loci}. It follows that $\mathrm{pr}_{\R}^+ : \mathcal{T}^+_{\R}\to (-\varepsilon, 1+\varepsilon)$ is surjective. Since we already know that $\mathrm{pr}_{\R}^+ : \mathcal{T}^+_{\R}\to (-\varepsilon, 1+\varepsilon)$ is injective again by (iii) of Proposition~\ref{PRP:complex_loci}, the locus $\mathcal{T}^+_{\R}$ can be expressed as the graph of a function $\kappa_{\R}^+ : (-\varepsilon, 1+\varepsilon)\to \R$ which is real analytic by (iv) of Proposition~\ref{PRP:complex_loci}. 
\end{proof}

Next, consider the real locus $\mathcal{T}_{\R}^-\equiv\mathcal{T}^-\cap \mathcal{F}^-_{\R}$. Since it consists of two parts $\mathcal{T}_{i, \R}^-\equiv\mathcal{T}^-_i\cap \mathcal{F}^-_{\R}$ ($i=1, 2$) in this case, we need to verify which part corresponds to the tangency locus of $W_{41\overline{0}}^s(a, b)^-$ and $W_{\overline{43}4124}^u(a, b)^-_{\mathrm{inner}}$.

\begin{prp}
The following properties hold for $\mathcal{T}_{i, \R}^-$ ($i=1, 2$).
\begin{enumerate}
\renewcommand{\labelenumi}{(\roman{enumi})}
\item We have $(a, b)\in \mathcal{T}^-_{1, \R}\cup\mathcal{T}^-_{2, \R}$ iff $(a, b)\in\mathcal{F}^-_{\R}$ and $W_{41\overline{0}}^s(a, b)^-$ intersects one of the irreducible components of $W_{\overline{43}4124}^u(a, b)^-$ tangentially in $\R^2$.
\item There exists a real analytic function:
\[\kappa_{i, \R}^- : (-1-\varepsilon, \varepsilon) \longrightarrow \R\]
so that $\mathcal{T}^-_{i, \R}$ coincides with the graph of $\kappa_{i, \R}^-$.
\end{enumerate}
\label{PRP:real_loci_negative}
\end{prp}

\begin{proof}
The proof of this claim is identical to the previous one, hence omitted. 
\end{proof}

\begin{figure}
  \includegraphics[height=5.5cm]{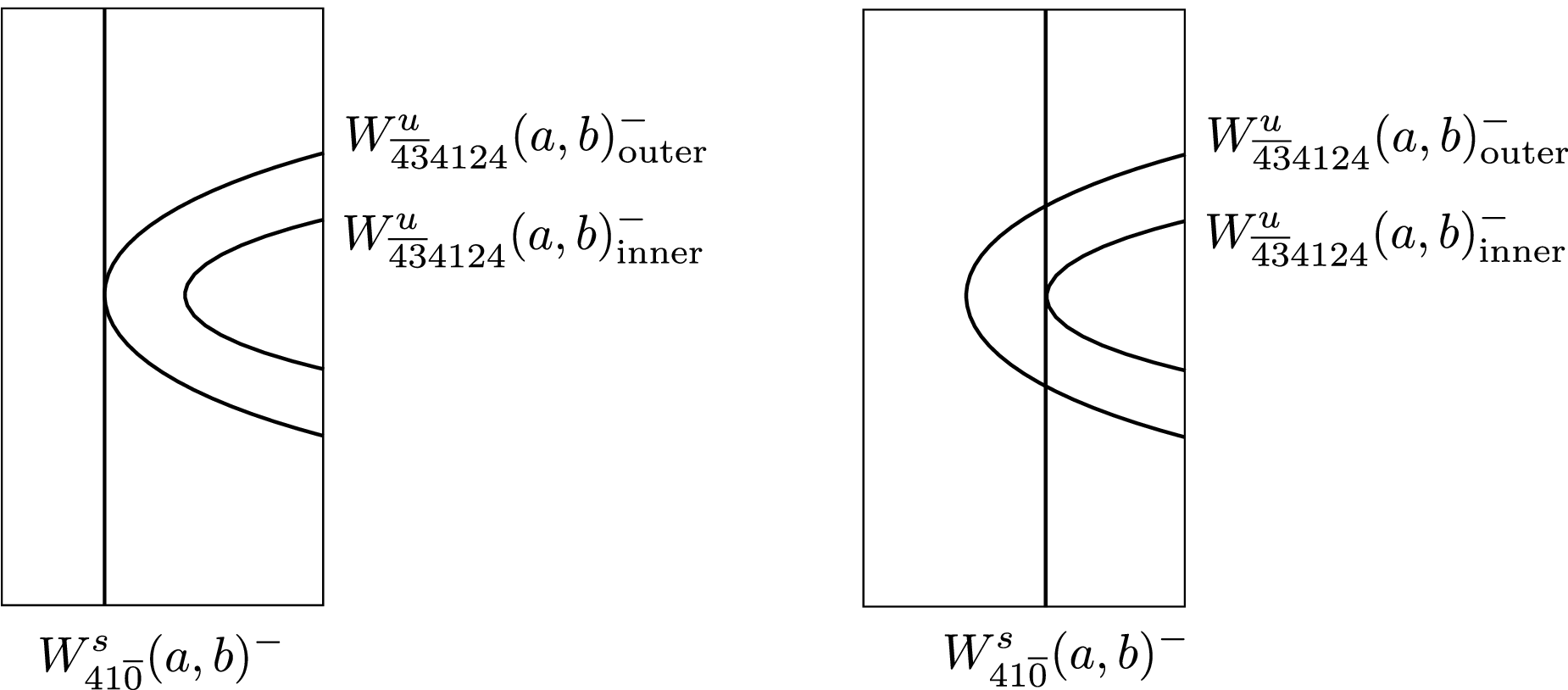}
  \caption{No simultaneous tangencies.}
  \label{FIG:simultaneous}
\end{figure}

Now let us prove the Main Theorem in Section~\ref{section1}.

\begin{proof}[Proof of the Main Theorem]
Consider the case $b<0$. Since the existence of tangency between $W_{41\overline{0}}^s(a, b)^-$ and $W_{\overline{43}4124}^u(a, b)^-_{\mathrm{outer}}$ implies the non-existence of tangency between $W_{41\overline{0}}^s(a, b)^-$ and $W_{\overline{43}4124}^u(a, b)^-_{\mathrm{inner}}$ and vise versa (see Figure~\ref{FIG:simultaneous}), we see $\mathcal{T}^-_{1, \R}\cap \mathcal{T}^-_{2, \R}\cap \{b<0\}=\emptyset$. It follows that $\kappa_{1, \R}^-(b)\ne\kappa_{2, \R}^-(b)$ holds for $-1-\varepsilon<b<0$, hence we may assume $\kappa_{1, \R}^-(b)>\kappa_{2, \R}^-(b)$ for $-1-\varepsilon<b<0$. Let us write $\kappa_{\R}^-(b)\equiv\kappa^-_{1, \R}(b)$  for $-1-\varepsilon<b<\varepsilon$ and put $a_{\mathrm{tgc}}(b)\equiv \kappa_{\R}^-(b)$ for $-1-\varepsilon<b<0$. Since $\kappa_{\R}^-(b)$ is continuous for $-1-\varepsilon<b<\varepsilon$ and $\kappa_{\R}^-(0)=2$, we have $\lim_{b \to -0} a_{\mathrm{tgc}}(b)=2$. Below we show that the function $a_\mathrm{tgc}$ satisfies (i) and (ii) in the Main Theorem and that the H\'enon map $f_{a, b}$ with $a=a_{\mathrm{tgc}}(b)$ has exactly one orbit of heteroclinic tangencies in the case $b<0$. Proof for the case $b>0$ is similar by letting $a_{\mathrm{tgc}}(b)\equiv \kappa_{\R}^+(b)$ for $1+\varepsilon>b>0$ and using Proposition~\ref{PRP:real_loci_positive}, hence omitted.

\begin{figure}
  \includegraphics[height=8.5cm]{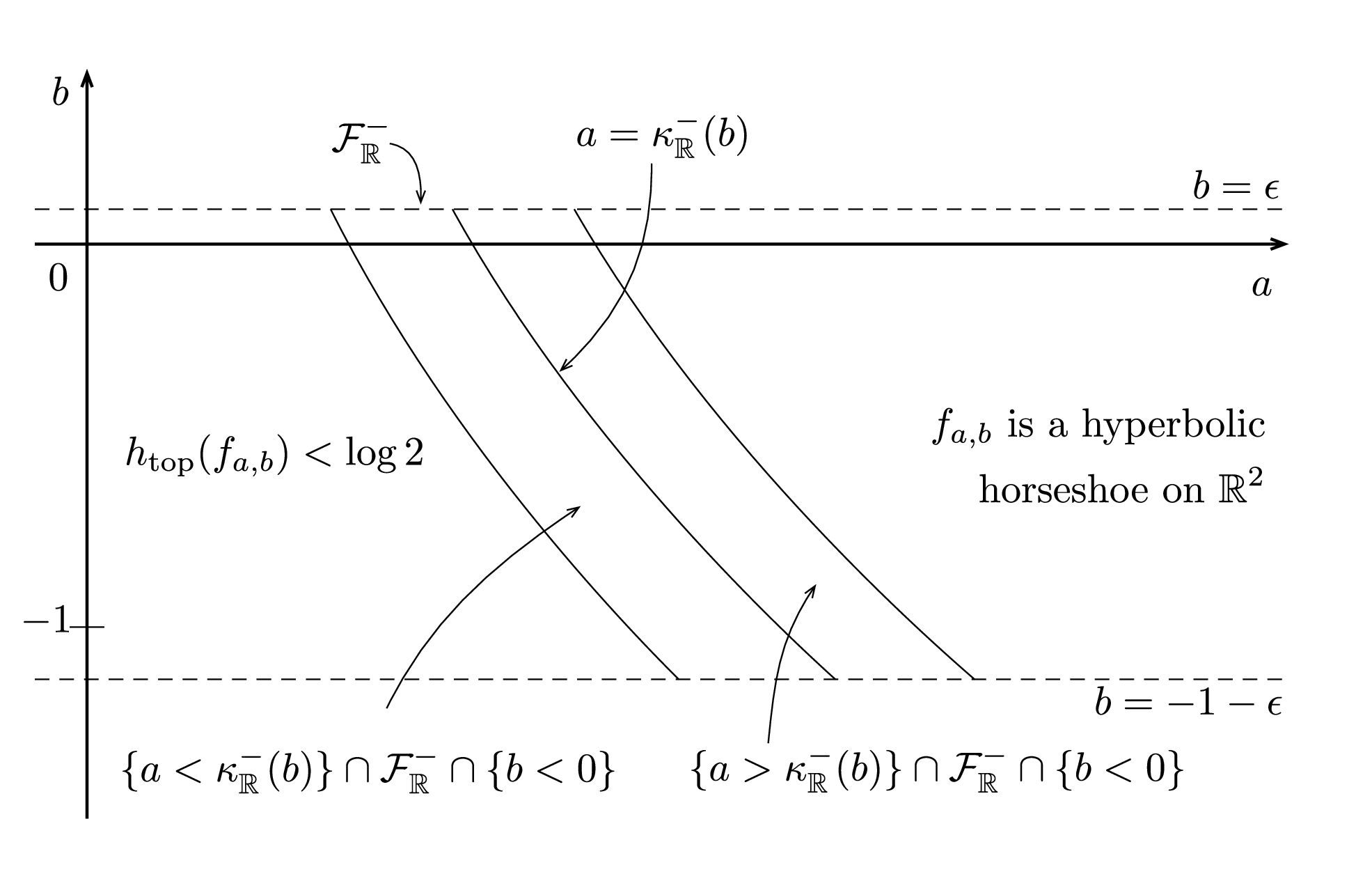}
  \caption{Proof of the Main Theorem.}
\label{FIG:proof_MT}
\end{figure}

First, let us show that the real analytic function $a_{\mathrm{tgc}}$ satisfies (ii) of the Main Theorem. Thanks to (ii) of Proposition~\ref{PRP:real_loci_negative}, $(\mathcal{F}^-_{\R}\cap\{b<0\})\setminus \mathcal{T}_{\R}^-$ consists of two connected components $\{a>\kappa_{\R}^-(b)\}\cap\mathcal{F}^-_{\R}\cap \{b<0\}$ and $\{a<\kappa_{\R}^-(b)\}\cap\mathcal{F}^-_{\R}\cap \{b<0\}$ (see Figure~\ref{FIG:proof_MT}). In each of these components, either the condition $\mathrm{card}(W_{41\overline{0}}^s(a, b)^-\cap W^u_{\overline{43}4124}(a, b)^-_{\mathrm{inner}})=2$ or the condition $\mathrm{card}(W_{41\overline{0}}^s(a, b)^-\cap W^u_{\overline{43}4124}(a, b)^-_{\mathrm{inner}})=0$ holds for all parameters in the component. Since $\{a>\kappa_{\R}^-(b)\}\cap\mathcal{F}^-_{\R}\cap \{b<0\}$ contains a hyperbolic horseshoe parameter by (iii) of Theorem~\ref{THM:quasitrichotomy} (Quasi-Trichotomy), we see that $(a, b)\in\{a>\kappa_{\R}^-(b)\}\cap\mathcal{F}^-_{\R}\cap \{b<0\}$ implies $\mathrm{card}(W_{41\overline{0}}^s(a, b)^-\cap W^u_{\overline{43}4124}(a, b)^-_{\mathrm{inner}})=2$. Similarly, since $\{a<\kappa_{\R}^-(b)\}\cap\mathcal{F}^-_{\R}\cap \{b<0\}$ contains a non-maximal entropy parameter by (i) of Theorem~\ref{THM:quasitrichotomy} (Quasi-Trichotomy), we see that $(a, b)\in\{a<\kappa_{\R}^-(b)\}\cap\mathcal{F}^-_{\R}\cap \{b<0\}$ implies $\mathrm{card}(W_{41\overline{0}}^s(a, b)^-\cap W^u_{\overline{43}4124}(a, b)^-_{\mathrm{inner}})=0$. By combining these, we have
\begin{equation}
a>\kappa_{\R}^-(b) \Longleftrightarrow \mathrm{card}(W_{41\overline{0}}^s(a, b)^-\cap W^u_{\overline{43}4124}(a, b)^-_{\mathrm{inner}})=2
\label{intersection1}
\end{equation} 
and 
\begin{equation}
a\geq\kappa_{\R}^-(b) \Longleftrightarrow \mathrm{card}(W_{41\overline{0}}^s(a, b)^-\cap W^u_{\overline{43}4124}(a, b)^-_{\mathrm{inner}})\geq 1
\label{intersection2}
\end{equation}
for $(a, b)\in\mathcal{F}^-_{\R}\cap \{b<0\}$. Now, the claim (ii) of the Main Theorem for $(a, b)\in\mathcal{F}^-_{\R}\cap\{b<0\}$ follows from Eqn.~(\ref{intersection2}) and Theorem~\ref{THM:maximalentropy} (Maximal Entropy). Together with Theorem~\ref{THM:quasitrichotomy} (Quasi-Trichotomy) for $(a, b)$ outside $\mathcal{F}^-_{\R}\cap\{b<0\}$, we obtain (ii) of the Main Theorem.

Next, let us prove that $a_{\mathrm{tgc}}$ satisfies (i) of the Main Theorem. By (ii) of the Main Theorem, we see $\mathcal{M}_{\R}\cap \mathcal{F}^-_{\R}\cap\{b<0\}=\{a\geq\kappa_{\R}^-(b)\}\cap \mathcal{F}^-_{\R}\cap \{b<0\}$. Since $\mathcal{H}_{\R}$ is an open subset of $\mathcal{M}_{\R}$, this yields $\mathcal{H}_{\R}\cap \mathcal{F}^-_{\R}\cap\{b<0\}\subset \{a>\kappa_{\R}^-(b)\}\cap \mathcal{F}^-_{\R}\cap \{b<0\}$. 

Conversely, take $(a, b)\in\{a>\kappa_{\R}^-(b)\}\cap\mathcal{F}^-_{\R}\cap \{b<0\}$. Then, by Eqn.~(\ref{intersection1}) we have the condition $\mathrm{card}(W_{41\overline{0}}^s(a, b)^-\cap W^u_{\overline{43}4124}(a, b)^-_{\mathrm{inner}})=2$. As in Theorem~\ref{THM:maximalentropy} (Maximal Entropy), this is equivalent to $\htop(f_{a, b}|_{\R^2})=\log 2$. By Theorem 10.1 in~\cite{BLS} this implies $K_{a, b}\subset \R^2$. By Corollary~\ref{COR:no_tangency_negative}, the condition $\mathrm{card}(W_{41\overline{0}}^s(a, b)^-\cap W^u_{\overline{43}4124}(a, b)^-_{\mathrm{inner}})=2$ also yields that there is no tangency between $W^u(p_3)$ and $W^s(p_1)$ when $(a, b)\in\mathcal{F}^-_{\R}\cap \{b<0\}$. Thanks to Theorems 2 and 3 in~\cite{BS1}, this implies the uniform hyperbolicity of $f_{a, b}$ on $K_{a, b}$. Since $\{a>\kappa_{\R}^-(b)\}\cap\mathcal{F}^-_{\R}\cap \{b<0\}$ is connected and contains a hyperbolic horseshoe parameter by Theorem~\ref{THM:quasitrichotomy} (Quasi-Trichotomy), we see that $f_{a, b}$ is a hyperbolic horseshoe on $\R^2$ for $(a, b)\in\{a>\kappa_{\R}^-(b)\}\cap\mathcal{F}^-_{\R}\cap \{b<0\}$ due to its structural stability. Hence the claim (i) of the Main Theorem holds for $(a, b)\in\mathcal{F}^-_{\R}\cap\{b<0\}$. Together with Theorem~\ref{THM:quasitrichotomy} (Quasi-Trichotomy) for $(a, b)$ outside $\mathcal{F}^-_{\R}\cap\{b<0\}$, we obtain (i) of the Main Theorem.

Finally, let us show that the H\'enon map $f_{a, b}$ with $a=a_{\mathrm{tgc}}(b)$ has exactly one orbit of heteroclinic tangencies when $b<0$. By the discussion above, we see that $\mathrm{card}(W_{41\overline{0}}^s(a, b)^-\cap W^u_{\overline{43}4124}(a, b)^-_{\mathrm{inner}})=1$. This implies that the unique point in $W_{41\overline{0}}^s(a, b)^-\cap W^u_{\overline{43}4124}(a, b)^-_{\mathrm{inner}}$ is a heteroclinic tangency of $W^u(p_3)$ and $W^s(p_1)$. Conversely, let $z$ be any point of heteroclinic tangency between $W^u(p_3)$ and $W^s(p_1)$. Since $z\in K_{a, b}$, there is a backward admissible sequence $(i_n)_{n\leq 0}$ different from $\overline{0}$ so that $f^n_{a, b}(z)\in\mathcal{B}^-_{i_n}$ for $n\leq 0$ by Proposition~\ref{PRP:coding_negative}. Thanks to the diagram of admissible transitions $\mathfrak{T}^-$ (see Figure~\ref{FIG:admissible_negative}), we know that there exists $n_0\leq 0$ so that $i_{n_0}=4$ which means $f_{a, b}^{n_0}(z)\in\mathcal{B}^-_{4, \R}$. Again since $\mathrm{card}(W_{41\overline{0}}^s(a, b)^-\cap W^u_{\overline{43}4124}(a, b)^-_{\mathrm{inner}})=1$ and the other pieces of $W^u(p_3)$ and $W^s(p_1)$ in $\mathcal{B}^-_{4, \R}$ intersect at two points (hence they are not tangential), it follows that $f_{a, b}^{n_0}(z)$ is the unique intersection of $W_{41\overline{0}}^s(a, b)^-$ and $W^u_{\overline{43}4124}(a, b)^-_{\mathrm{inner}}$. This implies that $W^u(p_3)$ and $W^s(p_1)$ have exactly one orbit of heteroclinic tangencies.

Argument for $b>0$ is similar, and this finishes the proof of the Main Theorem.
\end{proof}

As a consequence of this proof, we obtain a characterization for a H\'enon map to be a hyperbolic horseshoe on $\R^2$ in terms of the special intersections.

\begin{thm}[\textbf{Hyperbolic Horseshoes}]
When $(a, b)\in \mathcal{F}^+_{\R}\cap\{b>0\}$, $f_{a, b}$ is a hyperbolic horseshoe on $\R^2$ iff $\mathrm{card}(W_{31\overline{0}}^s(a, b)^+\cap W_{\overline{0}23}^u(a, b)^+)=2$. When $(a, b)\in \mathcal{F}^-_{\R}\cap\{b<0\}$, $f_{a, b}$ is a hyperbolic horseshoe on $\R^2$ iff $\mathrm{card}(W_{41\overline{0}}^s(a, b)^-\cap W^u_{\overline{43}4124}(a, b)^-_{\mathrm{inner}})=2$.
\label{THM:hyperbolichorseshoe}
\end{thm}

Compare the above result with Theorem~\ref{THM:maximalentropy} (Maximal Entropy).

In~\cite{BS4} characterizations of $\mathcal{H}_{\R}$ and $\mathcal{M}_{\R}$ similar to Theorem~\ref{THM:maximalentropy} and Theorem~\ref{THM:hyperbolichorseshoe} in a certain subregion of the parameter space (denoted as $\mathcal{W}_{\ast}$ in~\cite{BS4}) have been given in terms of symbolic dynamics with respect to a family of three boxes. We note that both Theorem~\ref{THM:maximalentropy} and Theorem~\ref{THM:hyperbolichorseshoe} hold for all values of $b$, but the results in~\cite{BS4} hold for approximately $-0.5<b<0.4$ (see Appendix~\ref{appendix:comparison}).

\newpage

\section{Proofs Involving Computer-Assistance}\label{section6}

This section is devoted to explaining the ideas of our rigorous numerics and showing how to verify the numerical criteria which are essential to the proof of the Main Theorem. 
We begin with some remarks on the interval arithmetic, the most fundamental machinery in our rigorous numerics, in Subsection~\ref{subsection6.1}. In Subsection~\ref{subsection6.2} we introduce two numerical algorithms based on the interval arithmetic, the interval Krawczyk method and the set-oriented algorithms.
Subsection~\ref{subsection6.3} is devoted to the data structure of our computation and here we explain how we practically handle the system of projective boxes appeared in Theorem~\ref{THM:quasitrichotomy}.
Finally, we present the proofs of lemmas which involve computer-assistance in Subsection~\ref{subsection6.4}.

All algorithms are implemented in \texttt{C/C++} and the entire source code is available at
\begin{center}
\texttt{http://www.isc.chubu.ac.jp/zin\_arai/locus/}
\end{center}
as well as the data necessary for the computation.

\subsection{Interval arithmetic}\label{subsection6.1}
We will not give the precise definition of the interval arithmetic here; instead, we focus on how it works in our setting of the complex H\'enon maps. For the basic and general properties of the interval arithmetic, see \cite{M} for example.

Most of our rigorous verification takes the following form: given a continuous map $f_\lambda$ depending on a parameter $\lambda \in \mathrm{\Lambda}\subset\R^l$ and given sets $X\subset\R^m$ in the domain and $Y\subset\R^n$ in the range of $f_\lambda$, we want to show that $f_\lambda(X) \subset Y$ holds for all $\lambda \in \mathrm{\Lambda}$.\footnote{More generally, $X$ and $Y$ may also depend on $\lambda \in \mathrm{\Lambda}$ and we want to show that $f_{\lambda}(X_{\lambda}) \subset Y_{\lambda}$ holds for all $\lambda$. The following argument can equally be applied to this case with $X$ replaced by $\bigcup_{\lambda \in \mathrm{\Lambda}}X_{\lambda}$ and $Y$ by $\bigcap_{\lambda \in \mathrm{\Lambda}} Y_{\lambda}$.}
In our rigorous computations, $f_\lambda$ will be the H\'enon map $f_{a, b}$ or its higher iterations, or their derivatives. Remark that although the H\'enon map itself is a polynomial map, we need to handle rational maps since we often use projective coordinates.
We denote the union of $f_{\lambda}(X)$ over all $\lambda \in \mathrm{\Lambda}$ by $f_{\mathrm{\Lambda}}(X)$. 

%
%
The fundamental difficulty here is that the set $f_{\mathrm{\Lambda}}(X)$ can not be directly obtained using computers due to numerical errors such as the rounding error.
However, with the help of interval arithmetic, we can find a set that rigorously contains $f_\mathrm{\Lambda}(X)$.
For this, we first enclose $X$ and $\mathrm{\Lambda}$ by rectangular sets in $\R^m$ and $\R^l$ (that is, products of closed intervals; we call them ``cubes'') respectively and then apply interval arithmetic for each component $f^i_{\lambda}$ of the map $f=(f^1_{\lambda}, \dots, f^n_{\lambda})$.
As a consequence we obtain a rectangular set in $\R^n$ containing $f_{\mathrm{\Lambda}}(X)$ rigorously. We denote this set by $F_{\mathrm{\Lambda}}(X)$ and call it an outer approximation of $f_{\mathrm{\Lambda}}(X)$. If $F_{\mathrm{\Lambda}}(X) \subset Y$ holds, then it follows $f_{\mathrm{\Lambda}}(X) \subset Y$, as required.

In practice, it often happens that even when we fail to verify $F_{\mathrm{\Lambda}}(X) \subset Y$, there exist coverings $\{X_i\}$ of $X$ and $\{\mathrm{\Lambda}_j\}$ of $\mathrm{\Lambda}$ by smaller cubes such that we can show $F_{\mathrm{\Lambda}_i}(X_j) \subset Y$ for all pairs of $i$ and $j$. In this case, we still have the same conclusion; namely, $f_{\lambda}(X) \subset Y$ for all $\lambda \in \mathrm{\Lambda}$. Thus, we want to subdivide the domain of the map and the parameter space into pieces as small as our computational power allows.

In fact, for the parameter space, we apply the following subdivision. First we subdivide $\mathcal{F}^{\pm}_{\R}$ using small parallelograms with two edges parallel to the $a$-axis and two other edges parallel to the graph of $a^{\pm}_{\mathrm{aprx}}$. For each parallelogram, we make the smallest rectangle containing it. Finally by taking the product of these rectangles and a subdivision of the $\Im(b)$-axis by small intervals, we have a covering of $\mathcal{F}^{\pm}$ by products of intervals as desired. The size of subdivision elements in $\mathcal{F}^+$ is at most $0.005$, $0.01$ and $0.001$ for the $\Re(a)$-, $\Im(a)$- and $\Re(b)$-directions, respectively. For $\mathcal{F}^-$, it is at most $0.001875$, $0.01$ and $0.0005$. Depending on parameters and conditions to be checked, we sometimes subdivide a subdivision element into further smaller pieces. 
The subdivision of the domain of the map is executed inductively in our algorithms, as we will see in Subsection~\ref{subsection6.4}.

Finally, we remark that the same argument can be applied to verify that $f_\lambda(X) \cap Y = \emptyset$ holds for all $\lambda \in \mathrm{\Lambda}$.

\subsection{Useful algorithms}\label{subsection6.2}
Here we discuss two distinguished numerical algorithms extensively used in our proofs.
One is the interval Krawczyk method, which is used to establish the existence of periodic points with very high accuracy. The other is the set-oriented algorithm, which is introduced for rigorously bounding dynamical objects such as the Julia set, invariant manifolds, etc.

\medskip

\noindent
(i) \textit{Interval Krawczyk method.} 
Below we review the ideas behind the interval Krawczyk method. Basically, it is obtained as a modification of the well-known Newton's root-finding method adapted to the interval arithmetic.

Let $g : \R^n\to \R^n$ be a smooth map. The Newton's method for solving $g(x) = 0$ is given by
\[N_g(x)=x-(Dg(x))^{-1}g(x).\]
In general, however, it is not easy to check that $Dg(U)$ is invertible for a small neighborhood $U$ of $x$ due to the wrapping effect of interval arithmetic. 

To overcome this difficulty, we modify the Newton's method as follows. For an invertible matrix $A$, let us define the modified Newton's method as
\[\widetilde{N}_{g, A}(x)=x-A g(x).\]
If the condition $\widetilde{N}_{g, A}(\mathrm{\Omega})\subset \mathrm{int}(\mathrm{\Omega})$ were verified for the product set $\mathrm{\Omega}\subset \R^n$ of $n$ closed intervals, the Brouwer fixed point theorem implies that there exists $x^{\ast}\in \mathrm{\Omega}$ with $g(x^{\ast})=0$. In practice, $A$ will be a numerical approximation of $(Dg(x))^{-1}$ for some $x \in \mathrm{\Omega}$. The point here is that $Dg(x)$ is not a matrix with interval components; it is just an usual matrix of floating point numbers. We can thus avoid taking the inverse of a matrix with interval components. However, since 
\[\mathrm{diam}(\mathrm{\Omega}-A g(\mathrm{\Omega}))\approx \mathrm{diam}(\mathrm{\Omega})+\mathrm{diam}(A g(\mathrm{\Omega}))>\mathrm{diam}(\mathrm{\Omega}),\]
it turns out that the condition $\widetilde{N}_{g, A}(\mathrm{\Omega})\subset \mathrm{int}(\mathrm{\Omega})$ always fails.

To fix this circumstance, Rudolf Krawczyk introduced the following idea (see equation (13) in page 177 of~\cite{Nm}). Fix a base-point $x_0\in \mathrm{\Omega}$. The interval mean-value theorem yields
\[\widetilde{N}_{g, A}(\mathrm{\Omega}) \subset \widetilde{N}_{g, A}(x_0)+D\widetilde{N}_{g, A}(\mathrm{\Omega}) (\mathrm{\Omega}-x_0) = x_0-A g(x_0)+(I-A \cdot Dg(\mathrm{\Omega})) (\mathrm{\Omega}-x_0),\]
where $I$ is the identity matrix. 

\begin{dfn}
The operator $K_{g, x_0, A}(\mathrm{\Omega}) \equiv x_0-A g(x_0)+(I-A \cdot Dg(\mathrm{\Omega})) (\mathrm{\Omega}-x_0)$ is called the \textit{interval Krawczyk operator} for $g$.
\end{dfn}

Note that $x_0-A g(x_0)$ is a point and $\mathrm{\Omega}-x_0$ is a translation of $\mathrm{\Omega}$. So, if the matrix $A$ is chosen so that $A \cdot Dg(\mathrm{\Omega})$ is close to $I$, we can conclude $\mathrm{diam}(K_{g, x_0, A}(\mathrm{\Omega}))< \mathrm{diam}(\mathrm{\Omega})$. With this operator we obtain 

\begin{prp}
If $K_{g, x_0, A}(\mathrm{\Omega})\subset \mathrm{int}(\mathrm{\Omega})$ holds for some $A$ and $x_0\in \mathrm{\Omega}$, there exists a unique point $x^{\ast}\in\mathrm{\Omega}$ so that $g(x^{\ast})=0$.
\label{PRP:Krawczyk}
\end{prp}

This result is employed to show (i) of Theorem~\ref{THM:quasitrichotomy} (Quasi-Trichotomy). Note that the uniqueness of the solution in $\mathrm{\Omega}$ is also guaranteed. For a proof, see Theorem 5.1.8 of~\cite{Nm}.

The interval Krawczyk method described above can immediately be applied to find a periodic point of a dynamical system $f : \R^n\to\R^n$, since a periodic point $p \in \R^n$ of period $k$ is nothing more than a zero of the equation $p - f^k(p) = 0$ satisfying $p - f^j(p) \ne 0$ for $j = 1, \ldots, k-1$. However, when $k$ is large or when the expansion of the map is strong, it is very difficult to apply the interval Krawczyk method to this equation. This is because the interval containing the true orbit gets expanded significantly in the unstable direction of $f$ and thus the inclusion property of the interval Krawczyk operator is very likely to fail. This is exactly what happens for our case $f = f_{a, b}$ and $k = 7$. For this reason, we transform the equation as follows. Let $p_1, p_2 \ldots, p_k \in \R^n$ be unknowns and consider the set of $k$ equations: 
\[p_2 - f(p_1) = 0, \ p_3 - f(p_2) = 0, \ldots, \ p_1 - f(p_k) = 0.\]
Obviously, the zeros of this system are the fixed points of $f^k$. The new equation is, although its dimension is $k$ times larger than the original equation, usually much easier to solve with the interval Krawczyk method since here we do not take any higher iteration of the map. See \cite{TW} for more detailed discussion on the application of the interval Krawczyk method to dynamical systems.

\medskip

\noindent
(ii) \textit{Set-oriented algorithm.} 
By the set-oriented algorithm we refer to a set of similar algorithms for rigorously enclosing invariant objects of dynamical systems.
In these algorithms, as the name suggests, we compute the time evolution of sets in the phase space instead of computing the orbit of each point \cite{DJ}. 
We combine the idea of the set-oriented algorithm with the interval arithmetic to obtain rigorous enclosures of dynamical objects such as periodic points, the maximal invariant sets and invariant manifolds.

Let $f: \R^n \to \R^n$ be a map and $R \subset \R^n$ a compact set on which we want to know the behavior of $f$. Consider a finite cubical grid on $R$ and assume that $R$ decomposes into small cubes $R = \bigcup_{i \in I} C_i$ where $I$ is the index set. By applying the interval arithmetic, we find a cube $D_i$ such that $f(C_i) \subset D_i$ rigorously holds for each $i \in I$.
The set $D_i$ is not a union of our cubical grid in general. Therefore, we next consider the set of grid elements intersecting with $D_i$.
That is, define a map $\mathtt{F} : I \to 2^I$ by $\mathtt{F}(i) \equiv \{j \in I \mid C_j \cap D_i \ne \emptyset \}$ and call it the \textit{cubical representation} of $f$.
Note that we have
\[f(C_i) \subset D_i \subset \bigcup_{j \in \mathtt{F}(i)} C_j.\]
%

Then we construct a directed graph $\mathcal{G}$ as follows. The set of vertices $V(\mathcal{G})$ of $\mathcal{G}$ is just $I$. We put an arrow from $i\in I$ to $j\in I$ if and only if $j \in \mathtt{F}(i)$. The graph $\mathcal{G}$ can be understood as a combinatorial representation of the dynamics of $f$ and in fact has a very nice property; if $x \in C_i$ and $f(x) \in C_j$ then there must be an arrow of $\mathcal{G}$ from $i$ to $j$. Thus, if there is no arrow from $i$ to itself, then it immediately implies that there is no fixed point of $f$ in $R_i$.
Similarly, if $C_i$ contains a periodic point of period $n$ whose orbit is contained in $R$, then there should be a cycle (closed walk) of consecutive arrows of length $n$ in $\mathcal{G}$.\footnote{Cycles may not be simple; they may contain repetitions of vertices, edges and self-loops.} Therefore, if we want to locate periodic points of period $n$, we remove the vertices $V'$ having no cycle of length $n$ from $V(\mathcal{G})$. Then the set: 
\[\bigcup_{i \in V(\mathcal{G})\setminus V'}C_i\] 
contains all periodic points of period $n$ that is contained in $R$. To have a better approximation, we simply refine the grid by subdividing remaining cubes, reconstruct the directed graph $\mathcal{G}$, and repeat the same procedure.

Set-oriented algorithms can also be applied to approximate maximal invariant sets and invariant manifolds, as we will see in Subsection~\ref{subsection6.4}.

%

\subsection{Numerical data}\label{subsection6.3}

\begin{table}
\[a^+_{\mathrm{aprx}}(1.00)=5.70, \qquad a^-_{\mathrm{aprx}}(-1.00)=6.20,\]
\[a^+_{\mathrm{aprx}}(0.90)=5.15, \qquad a^-_{\mathrm{aprx}}(-0.90)=5.60,\]
\[a^+_{\mathrm{aprx}}(0.80)=4.65, \qquad a^-_{\mathrm{aprx}}(-0.80)=5.04,\]
\[a^+_{\mathrm{aprx}}(0.70)=4.18, \qquad a^-_{\mathrm{aprx}}(-0.70)=4.52,\]
\[a^+_{\mathrm{aprx}}(0.60)=3.76, \qquad a^-_{\mathrm{aprx}}(-0.60)=4.04,\]
\[a^+_{\mathrm{aprx}}(0.50)=3.37, \qquad a^-_{\mathrm{aprx}}(-0.50)=3.61,\]
\[a^+_{\mathrm{aprx}}(0.40)=3.03, \qquad a^-_{\mathrm{aprx}}(-0.40)=3.21,\]
\[a^+_{\mathrm{aprx}}(0.30)=2.72, \qquad a^-_{\mathrm{aprx}}(-0.30)=2.85,\]
\[a^+_{\mathrm{aprx}}(0.20)=2.45, \qquad a^-_{\mathrm{aprx}}(-0.20)=2.53,\]
\[a^+_{\mathrm{aprx}}(0.10)=2.21, \qquad a^-_{\mathrm{aprx}}(-0.10)=2.25,\]
\[a^{\pm}_{\mathrm{aprx}}(0.00)=2.00.\]
\vspace{0.3cm}
\caption{The vertices of the piecewise affine functions $a^{\pm}_{\mathrm{aprx}}$.}
\label{TAB:aprx}
\end{table}

\begin{table}[ht]
  \begin{center}
    \begin{tabular}{cccc}
      \footnotesize
      \begin{minipage}{0.25\hsize}
\begin{verbatim}
tx[0] = 3.58844
ty[0] = 3.58844
tx[1] = 3.41867
ty[1] = 2.42305
tx[2] = -2.93181
ty[2] = 2.48933
tx[3] = -2.60315
ty[3] = 0.4062
tx[4] = 2.59251
ty[4] = -2.42276
tx[5] = 2.42305
ty[5] = -3.24747
tx[6] = 0.97798
ty[6] = -2.04485
tx[7] = 0.4062
ty[7] = -2.93181
\end{verbatim}
      \end{minipage}
      \begin{minipage}{0.25\hsize}
\begin{verbatim}
tx[8] = 0.97798
ty[8] = -2.04485
tx[9] = 0.4062
ty[9] = -2.93181
tx[10] = -2.39628
ty[10] = -0.75464
tx[11] = -2.04485
ty[11] = -2.49658
tx[12] = -2.93181
ty[12] = 2.48933
tx[13] = -2.04485
ty[13] = -2.49658
tx[14] = -3.24747
ty[14] = 2.42305
tx[15] = -2.42276
ty[15] = -2.42276
\end{verbatim}
      \end{minipage}
      \begin{minipage}{0.18\hsize}
\begin{verbatim}
ax[0] = 1.4
ax[1] = 1.4
ax[2] = 1.4
ax[3] = 1.4
ay[0] = 1.2
ay[1] = 1.2
ay[2] = 1.2
ay[3] = 1.2
bx[0] = 0.55
bx[1] = 0.3
bx[2] = 0.45
bx[3] = 0.27
by[0] = 0.23
by[1] = 0.3
by[2] = 0.3
by[3] = 0.6
\end{verbatim}
      \end{minipage}
      \begin{minipage}{0.25\hsize}
\begin{verbatim}
delta_Px[0] = 0.2
delta_Qx[0] = -0.15
delta_Py[0] = 0.1
delta_Qy[0] = -0.4
delta_Px[1] = 0.3
delta_Qx[1] = -0.55
delta_Py[1] = 0.3
delta_Qy[1] = -0.05
delta_Px[2] = 0.32
delta_Qx[2] = -0.22
delta_Py[2] = 0.25
delta_Qy[2] = -0.07
delta_Px[3] = 0.2
delta_Qx[3] = -0.2
delta_Py[3] = 0.1
delta_Qy[3] = -0.2
\end{verbatim}
      \end{minipage}
    \end{tabular}
\vspace{0.5cm}
    \caption{The data for boxes at $(a, b) = (5.7, 1.0)$.}
    \label{TAB:a570b100}
  \end{center}
\end{table}

In this subsection we present numerical data required for the proofs given in Subsection~\ref{subsection6.4}.

First we define approximations $a^{\pm}_{\mathrm{aprx}}$ of the first tangency curve $a_{\mathrm{tgc}} : \R^{\times} \to \R$.
They are defined to be the piecewise affine functions whose vertices are given in Table~\ref{TAB:aprx}.
Although the functions $a^{\pm}_{\mathrm{aprx}}$ are defined on $I^{\pm}_{\R}$ in Subsection~\ref{subsection2.1}, as we will see in the beginning of Subsection~\ref{subsection6.4}, our rigorous verification will be executed only for the case $0 \leq \Re(\pm b) \leq 1$ and therefore we define $a^+_{\mathrm{aprx}}$ only on $\{0\leq b\leq 1\}$ and $a^-_{\mathrm{aprx}}$ on $\{-1\leq b\leq 0\}$.
We remark that the values in Table~\ref{TAB:aprx} need not be rigorous and actually are not.
What we expect for $a^{\pm}_{\mathrm{aprx}}$ is that the actual tangency curve $a_{\mathrm{tgc}}$ stays in the neighborhoods $|a-a^{\pm}_{\mathrm{aprx}}(b)|\leq \chi^{\pm}(b)$.

In the proof of (iii) in Theorem~\ref{THM:quasitrichotomy} (Quasi-Trichotomy), the most fundamental data is a family of projective boxes $\{\mathcal{B}_i^{\pm}\}_i$ defined for all $(a, b) \in \mathcal{F}^{\pm}$. Just to glimpse the idea, we here present the data which is necessary to construct $\{\mathcal{B}_i^{\pm}\}_i$ for a selected parameter value in detail. 

Table~\ref{TAB:a570b100} shows the data of the boxes for $(a, b) = (5.7, 1.0)\in \mathcal{F}^{\pm}$. Each pair $(\verb|tx[k]|, \verb|ty[k]|)$ is computed as the coordinates of the intersection point $t^+_k\in \R^2$ ($0\leq k\leq 15$) in the trellis for $f_{a, b}$ as in Subsection~\ref{subsection2.2}.
Let $\mathcal{Q}^+_i$ be the quadrilateral formed by the four vertices $t^+_{4i}$, $t^+_{4i+1}$, $t^+_{4i+2}$ and $t^+_{4i+3}$ ($0\leq i\leq 3$). Take $L_u\equiv\C\times\{0\}\subset \C^2$ and $L_v\equiv\{0\}\times\C\subset\C^2$ and identify them with $\C$. First, we compute two foci $u$ and $v$ as the unique intersection points of the extensions of two vertical edges of $\mathcal{Q}^+_i$ and that of two horizontal edges of $\mathcal{Q}^+_i$ respectively. These foci together with $L_u$ and $L_v$ define the projective coordinates $(\pi^+_{u, i}, \pi^+_{v, i})$ associated with $\mathcal{Q}^+_i$ (see Figure~\ref{FIG:projective_coordinates}). 

\begin{figure}
  \includegraphics[height=8cm]{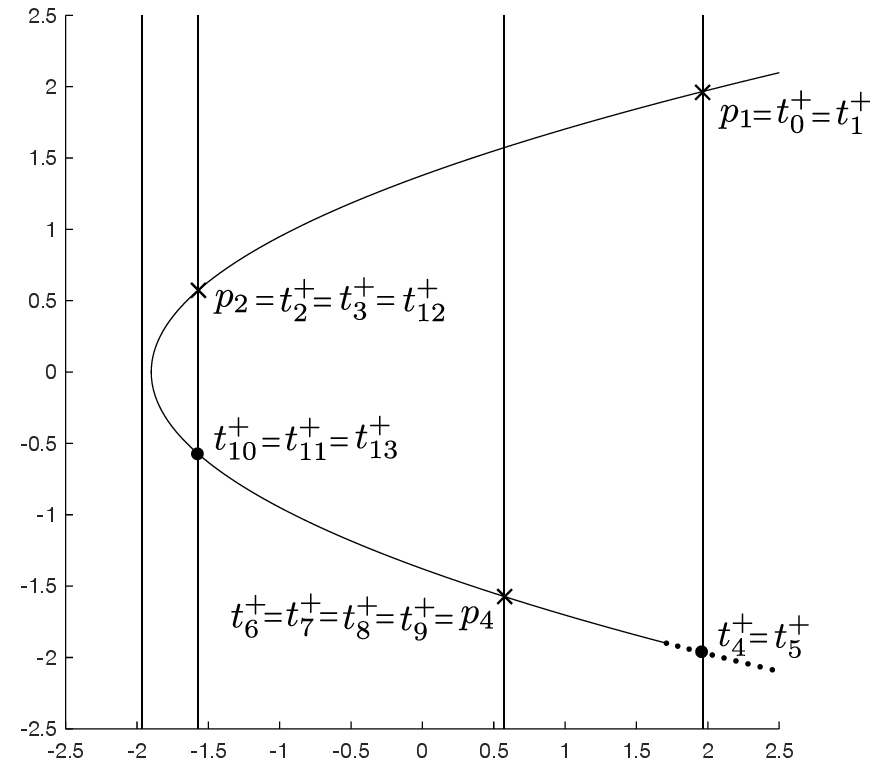} \\
  \includegraphics[height=9cm]{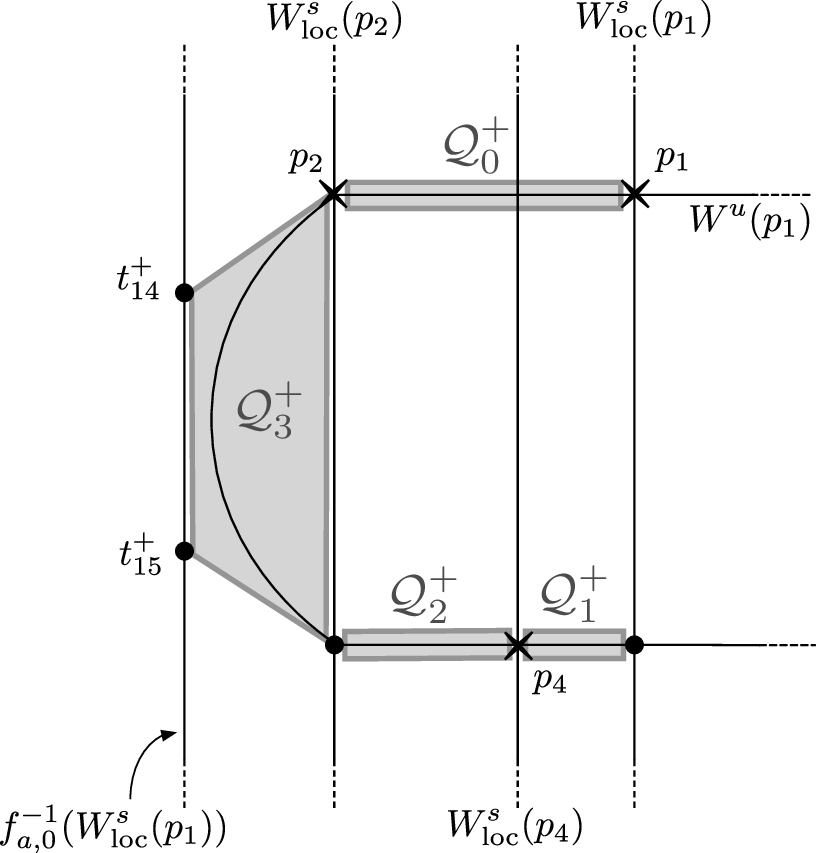}
  \caption{Above: trellis and the quadrilaterals $\{\mathcal{Q}^+_i\}^3_{i=0}$ for $(a, b)=(1.9, 0)$. Below: their cartoon images.}
  \label{FIG:trellis_positive_ghost}
\end{figure}

\begin{figure}
  \includegraphics[height=8cm]{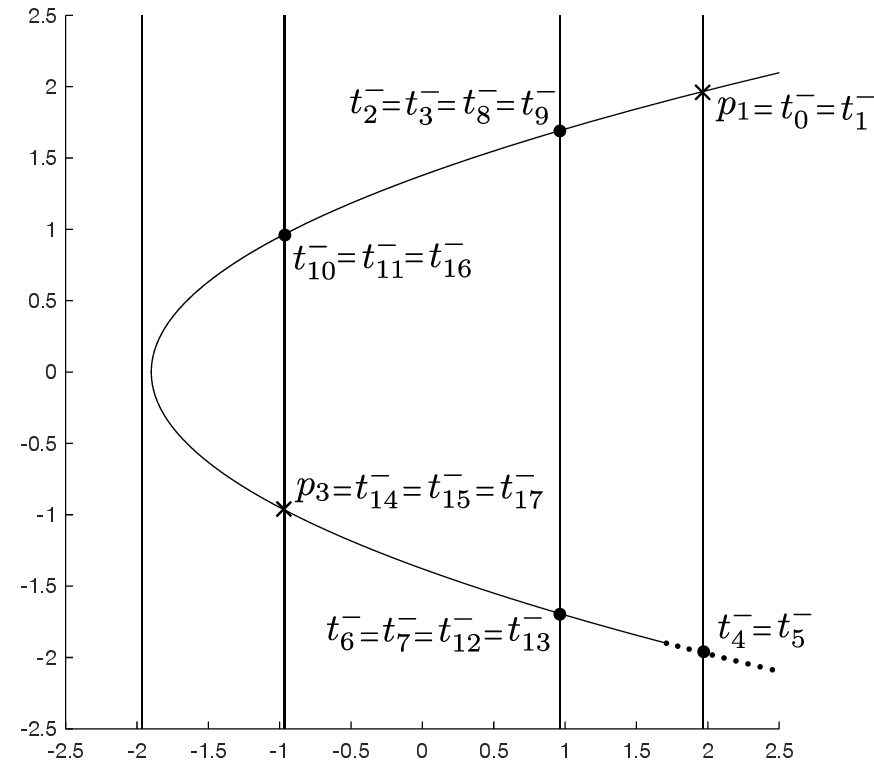} \\
  \includegraphics[height=9cm]{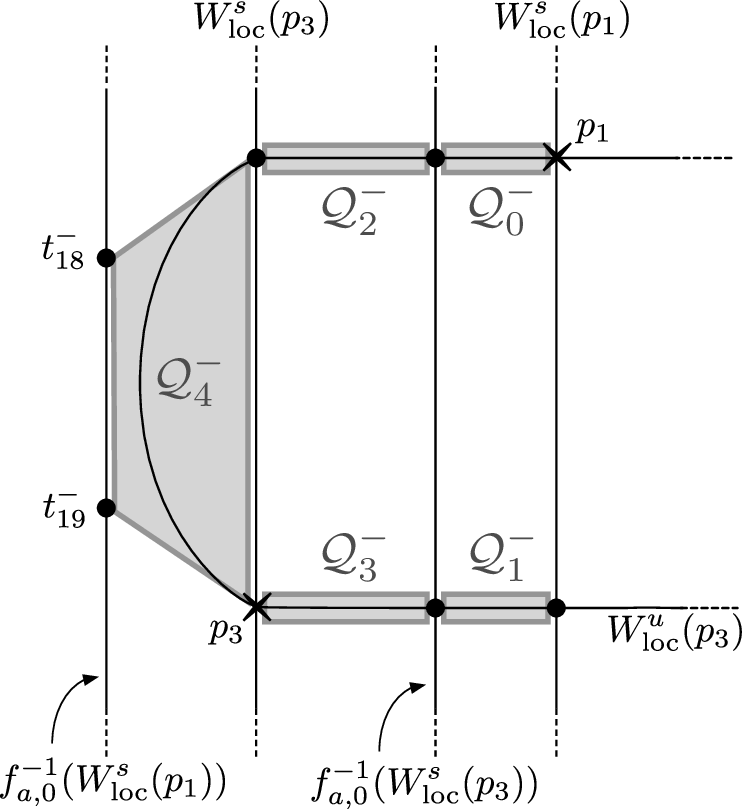}
  \caption{Above: trellis and the quadrilaterals $\{\mathcal{Q}^-_i\}^4_{i=0}$ for $(a, b)=(1.9, 0)$. Below: their cartoon images.}
  \label{FIG:trellis_negative_ghost}
\end{figure}

\begin{rmk}
When $(a, b) \in \mathcal{F}^+_{\R}$ and $b$ is close to zero, the two intersection points $t^+_{14}$ and $t^+_{15}$ may not exist or they may coincide. In this case, instead of using intersection points, we artificially define $t^+_{14}$ and $t^+_{15}$ so that they move continuously with respect to $(a, b) \in \mathcal{F}^+_{\R}$ and define a quadrilateral $\mathcal{Q}^+_3$ that satisfies the same criterions as $\mathcal{Q}^+_3$ for larger $b$ (see Figure~\ref{FIG:trellis_positive_ghost}). For $(a, b) \in \mathcal{F}^-_{\R}$ and $b$ close to zero, the two intersection points $t^-_{18}$ and $t^-_{19}$ may have the same problem. In this case, similarly we define $t^-_{18}$ and $t^-_{19}$ to define $\mathcal{Q}^{-}_4$ (see Figure~\ref{FIG:trellis_negative_ghost}).
\label{RMK:degenerate}
\end{rmk}

To construct a projective box $\mathcal{B}_i^+\equiv D^+_{u, i}\times_{\mathrm{pr}}D^+_{v, i}$ in $\C^2$ associated with $\mathcal{Q}^+_i$, we need to choose appropriate topological disks $D^+_{u, i}\subset L_u$ and $D^+_{v, i}\subset L_v$ (see Definition~\ref{DFN:projective_box_Q} in Subsection~\ref{subsection2.2}). The constants \verb|ax[i]|, \verb|ay[i]|, \verb|bx[i]|, \verb|by[i]|, \verb|delta_Px[i]|, \verb|delta_Qx[i]|, \verb|delta_Py[i]| and \verb|delta_Qy[i]| in Table~\ref{TAB:a570b100} are used to define these topological disks, as we will see below.  We note that the specific values of these constants are not ``canonical'', i.e. they are obtained by trial and error so that the resulted family of boxes $\{\mathcal{B}_i^+\}_{i=0}^3$ satisfies the (BCC).

Now, let us construct $D^+_{u, i}$ and $D^+_{v, i}$. First, the projection of two vertical edges of $\mathcal{Q}^+_i$ to $L_u$ via $\pi^+_{u, i}$ determines two points $q_{x, i}$ and $p_{x, i}$ with $q_{x, i} < p_{x, i}$ (see Figure~\ref{FIG:projective_box}). Set $P_X\equiv p_{x, 0}+\verb|delta_Px[0]|> 0$, $Q_X\equiv q_{x, 3}+\verb|delta_Qx[3]|< 0$, $P_Y\equiv p_{y, 3}+\verb|delta_Py[3]|>0$ and $Q_Y\equiv q_{y, 3}+\verb|delta_Qy[3]|<0$. Note that $P_X$, $Q_X$, $P_Y$ and $Q_Y$ are independent of $i$.

Given two constants $a_{x, i}\equiv\verb|ax[i]|$ and $b_{x, i}\equiv\verb|bx[i]|$, we define the ellipse $E^+_{u, i}$ to be the set of $u\in L_u$ satisfying
\[\left(\Re(u) - \frac{P_X + Q_X}{2} \right) ^2 + \left(\frac{a_{x, i}}{b_{x, i}} \, \Im(u) \right)^2 \leq\left(\frac{P_X - Q_X}{2}\right)^2.\]
Then, given two constants $\delta_{P_x, i}\equiv \verb|delta_Px[i]|$ and $\delta_{Q_x, i}\equiv \verb|delta_Qx[i]|$, we define the topological disk $D^+_{u, i}\subset E^+_{u, i}$ by
\[D^+_{u, i} \equiv E^+_{u, i} \cap \bigl\{u\in\C : q_{x, i} + \delta_{Q_x, i} \leq \Re(u) \leq p_{x, i} + \delta_{P_x, i}\bigr\}.\]
Similarly we define the disk $D^+_{v, i} \subset L_v$ as a part of the ellipse $E^+_{v, i}$ using $p_{y, i}$, $q_{y, i}$, $a_{y, i}=$ \verb|ay[i]|, $b_{y, i}=$ \verb|by[i]|, $\delta_{P_y, i}=$ \verb|delta_Py[i]| and $\delta_{Q_y, i}=$ \verb|delta_Qy[i]| via $\pi^+_{v, i}$. 

\begin{figure}
  \includegraphics[width=7cm]{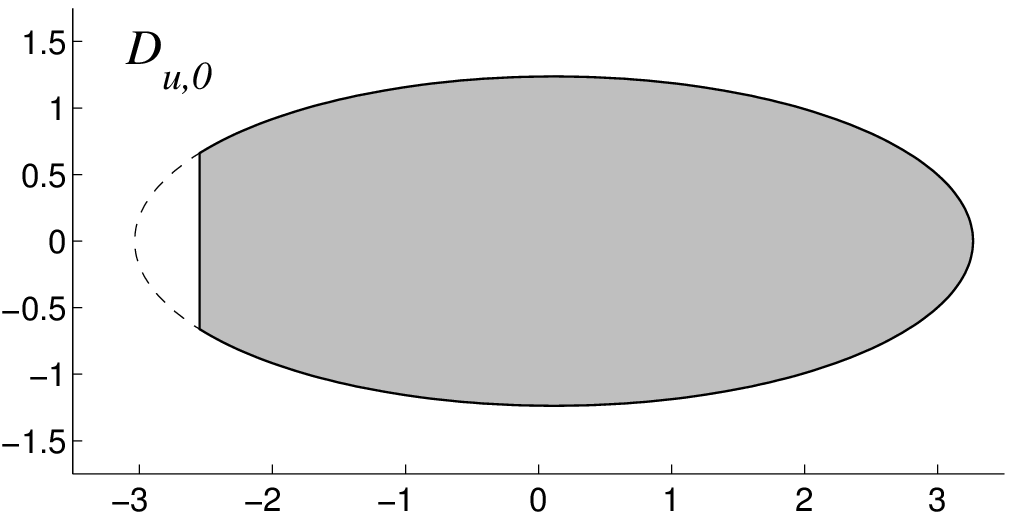}
  \qquad
  \includegraphics[width=7cm]{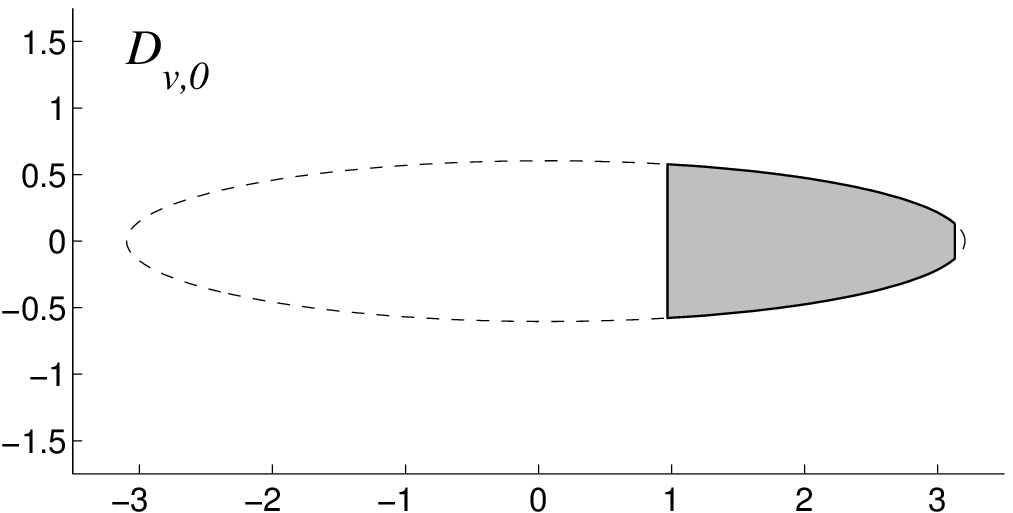}
  \\
  \vspace{3mm}
  \includegraphics[width=7cm]{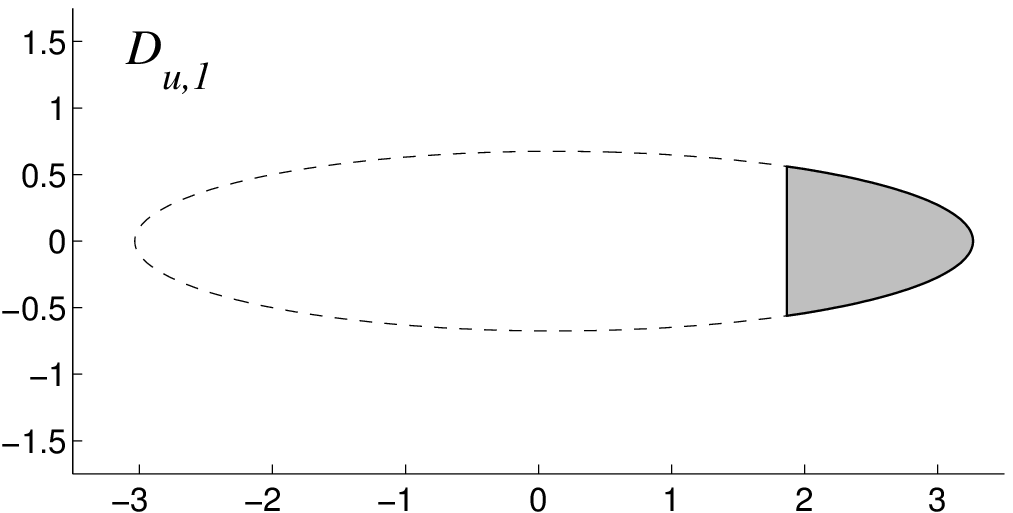}
  \qquad
  \includegraphics[width=7cm]{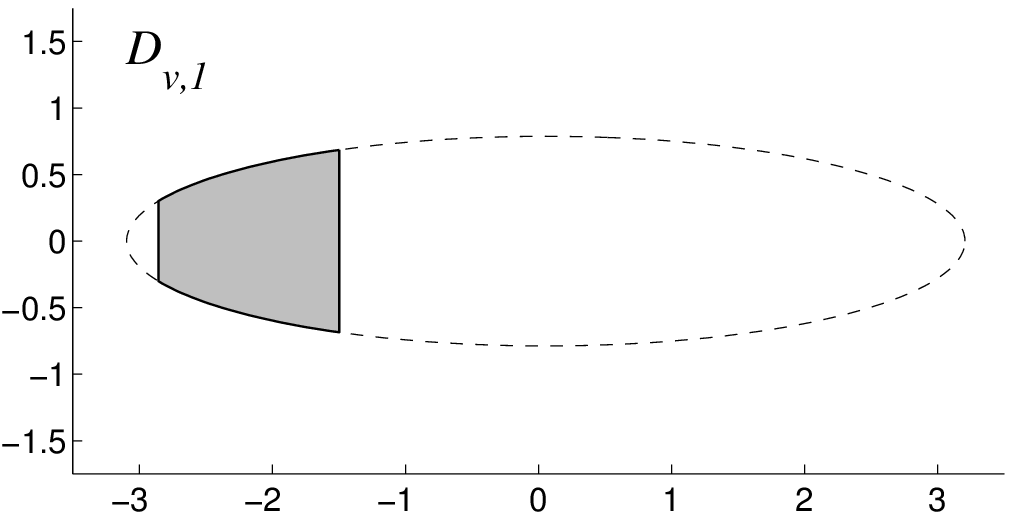}
  \\
  \vspace{3mm}
  \includegraphics[width=7cm]{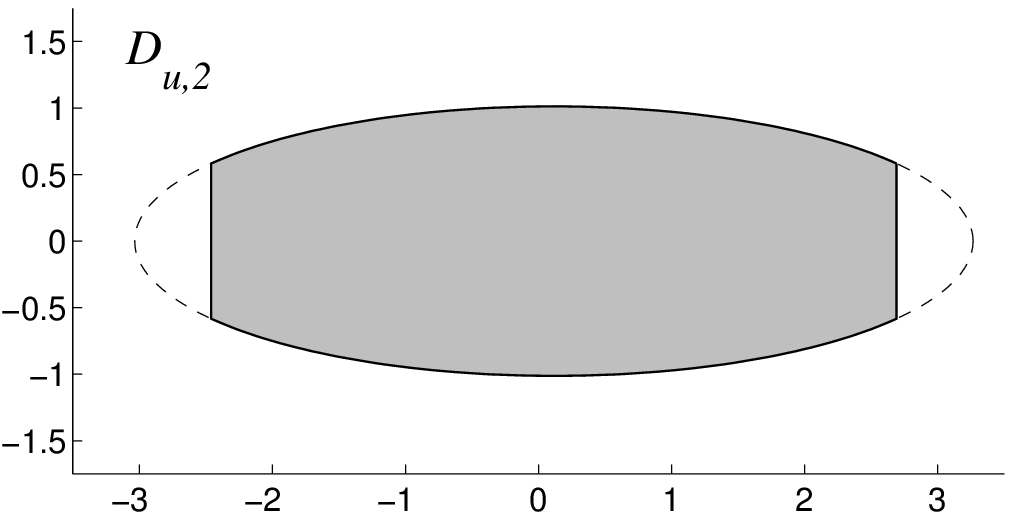}
  \qquad
  \includegraphics[width=7cm]{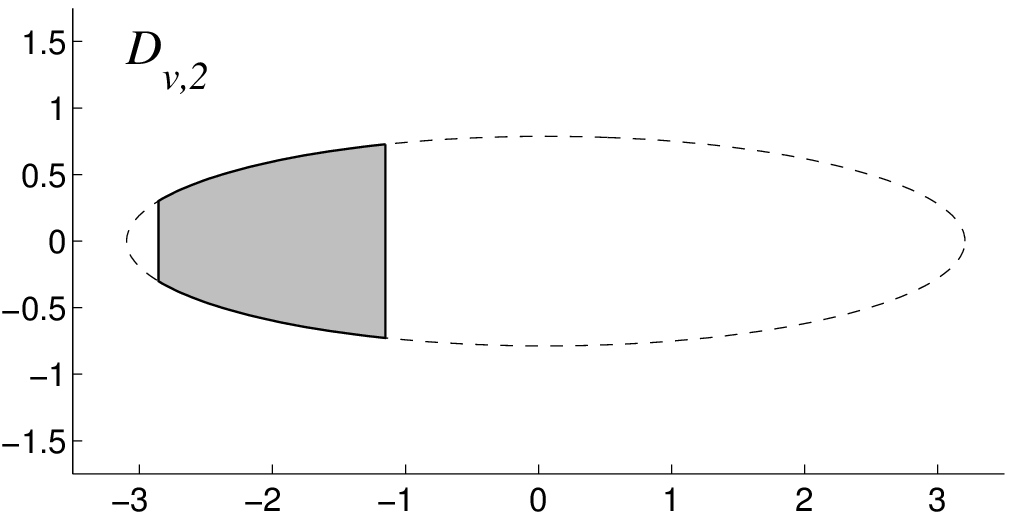}
  \\
  \vspace{3mm}
  \includegraphics[width=7cm]{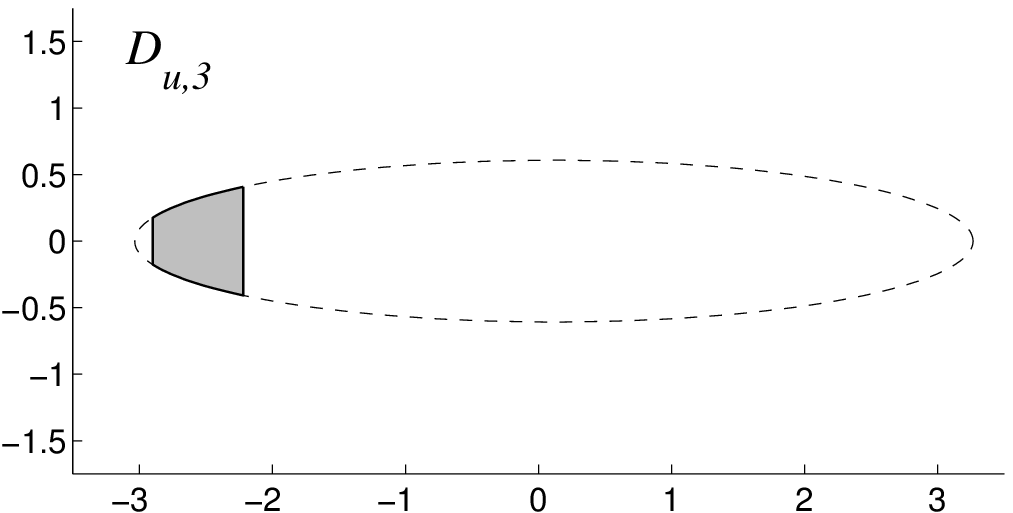}
  \qquad
  \includegraphics[width=7cm]{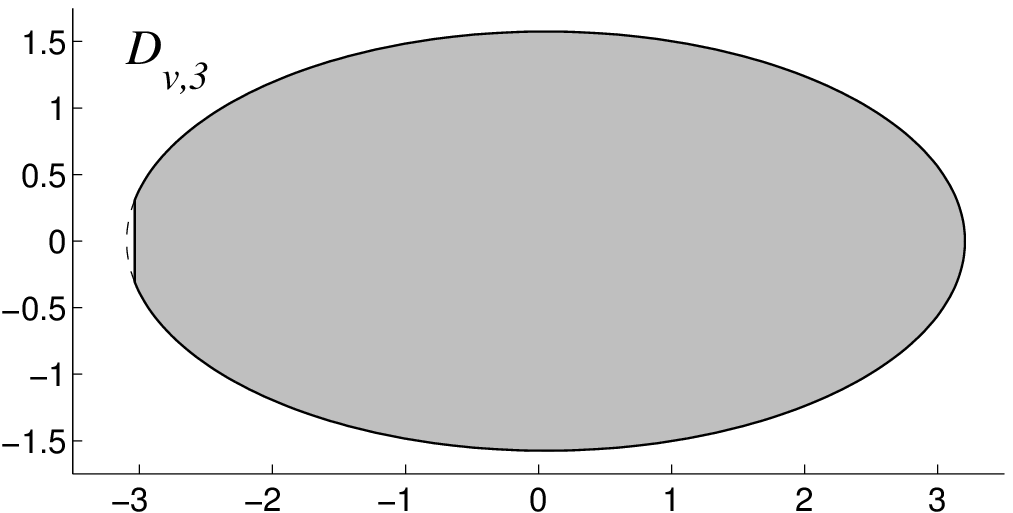}
  \caption{Shapes of $D^+_{u, i}$ and $D^+_{v, i}$ at $(a, b) = (5.7, 1.0)$.}
  \label{FIG:ellipses}
\end{figure}

Finally we take the product of these two topological disks with respect to the projective coordinates $(\pi^+_{u, i}, \pi^+_{v, i})$ to obtain the projective box $\mathcal{B}^+_i\equiv D^+_{u, i}\times_{\mathrm{pr}} D^+_{v, i}$ associated with $\mathcal{Q}^+_i$ (see Figure~\ref{FIG:projective_box} again). Figure~\ref{FIG:ellipses} shows the actual shapes of $D^+_{u, i}$ and $D^+_{v, i}$ so that $\{\mathcal{B}^+_i\}_{i=0}^3$ satisfies the (BCC), which implies the (CMC) for $(a, b)=(5.7, 1.0)$ as in (iii) of Theorem~\ref{THM:quasitrichotomy} (Quasi-Trichotomy). The construction of $\mathcal{B}^{\pm}_i$ for all $(a, b)\in\mathcal{F}^{\pm}$ will be explained in the next subsection where we prove Lemma~\ref{LMM:CMC}.

\subsection{Proofs of lemmas}\label{subsection6.4}

In this subsection we present the proofs of lemmas which require computer-assistance.

First we remark that we run computer assisted proofs only for the case $0 \leq \Re(b) \leq 1$ and $\Im(b) = 0$ nevertheless lemmas holds for all $b \in I^{\pm}$ (or, $I^{\pm}_\R$). This is because of the continuity of the map and the box systems, and the fact that all the statements involving rigorous interval arithmetic are verified with certain amounts of margin. That is, when our program verifies a statement $f_{\lambda}(X) \subset Y$ (see Subsection~\ref{subsection6.1}), it in fact guarantees that a small neighborhood of $f_{\lambda}(X)$, which is larger than $f_{\lambda}(X)$ at least by the smallest positive floating point number, is contained in $Y$ and therefore the continuity of the map implies that the same inclusion holds for all $\lambda'$ close enough to $\lambda$. Since the number of statements we verify is finite, we can choose small $\varepsilon$ and $\delta$ so that our lemmas hold for all $b \in I^{\pm}$.

\begin{proof}[Proof of Lemma~\ref{LMM:period7}]
To prove this, we show that for all $(a, b) \in \R\times I^{\pm}_{\R}$ with $-(b+1)^2/4 \leq a \leq a^{\pm}_{\mathrm{aprx}} (b) - \chi^{\pm}(b)$, there exists a periodic point of period $7$ in $\C^2 \setminus \R^2$.

The verification process goes as follows. We first construct a covering of the bounded set in the parameter space:
\[\bigl\{(a, b) \in \R\times I^{\pm}_{\R} : \ -(b+1)^2/4 \leq a \leq a^{\pm}_{\mathrm{aprx}} (b) - \chi^{\pm}(b)\bigr\}\]
by small rectangles of the form $A \times B$ where $A$ and $B$ are closed intervals. For each small rectangle $A \times B$, we select a parameter value $(a, b) \in A \times B$. We then use the conventional Newton's method to numerically find a candidate $\{(x_1, y_1), (x_2, y_2), \ldots, (x_7, y_7) \}$ of periodic orbit of period $7$ with respect to $f_{a, b}$ such that $(x_i, y_i) \in \C^2 \setminus \R^2$ for $i = 1, 2, \ldots, 7$. Next, we verify the inclusion assumption $K_{g, x_0, A}(\mathrm{\Omega})\subset \mathrm{int}(\mathrm{\Omega})$ in Proposition~\ref{PRP:Krawczyk} for a small rectangle $\mathrm{\Omega} \subset \C^{2 \times 7} \setminus \R^{2 \times 7}$ containing $(x_1, y_1, x_2, y_2, \ldots, x_7, y_7)$. This establish the existence of a fixed point $(x_{\ast}, y_{\ast})$ of $f^7_{a, b}$ in $\C^2 \setminus \R^2$ for all $(a, b) \in A \times B$. It easy to check that $(x_{\ast}, y_{\ast})$ is not fixed by $f_{a, b}$ and thus we conclude $(x_{\ast}, y_{\ast})$ is a periodic point of period $7$. 

For example, at the parameter value $(a, b) = (5.6, 1.0)$, we prove the existence of a periodic orbit of period $7$ such that
{\footnotesize
  \begin{align*}
    x_1 &\in [-2.81703, -2.80968] + [-0.044259, -0.036907]i, & y_1 &\in [-0.17505, -0.167697] + [0.233134, 0.240487]i\\
    x_2 &\in [2.48102, 2.48837] + [-0.012138, -0.004786]i, & y_2 &\in [-2.81703, -2.80968] + [-0.044259, -0.036907]i\\
    x_3 &\in [3.38331, 3.39066] + [-0.005141, 0.002212]i, & y_3 &\in  [2.48102, 2.48837] + [-0.012138, -0.004786]i\\
    x_4 &\in [3.38331, 3.39066] + [-0.005141, 0.002212]i, & y_4 &\in [3.38331, 3.39066] + [-0.005141, 0.002212]i\\
    x_5 &\in [2.48102, 2.48837] + [-0.012138, -0.004786]i, & y_5 &\in [3.38331, 3.39066] + [-0.005141, 0.002212]i\\
    x_6 &\in [-2.81703, -2.80968] + [-0.044259, -0.036907]i, & y_6 &\in [2.48102, 2.48837] + [-0.012138, -0.004786]i\\
    x_7 &\in [-0.17505, -0.167697] + [0.233134, 0.240487]i, & y_7 &\in [-2.81703, -2.80968] + [-0.044259, -0.036907]i.
  \end{align*}
}
Since the imaginary part of $x_1$ is non-zero, this orbit is not contained in $\R^2$.
Similarly, for $(a, b) = (5.6, 1.0)$, we prove the existence of a periodic orbit of period $7$ such that
{\footnotesize
  \begin{align*}
    x_1 &\in [3.2245, 3.23185] + [-0.005848, 0.001505]i, & y_1 &\in [-3.05792, -3.05057] + [0.011763, 0.019117]i\\
    x_2 &\in [1.26319, 1.27055] + [-0.002260, 0.005093]i, & y_2 &\in [3.2245, 3.23185] + [-0.005848, 0.001505]i\\
    x_3 &\in [-1.27055, -1.26319] + [-0.002260, 0.005093]i, & y_3 &\in [1.26319, 1.27055] + [-0.002260, 0.005093]i\\
    x_4 &\in [-3.23185, -3.2245] + [-0.005848, 0.001505]i, & y_4 &\in [-1.27055, -1.26319] + [-0.002260, 0.005093]i\\
    x_5 &\in [3.05057, 3.05792] + [0.011763, 0.019117]i, & y_5 &\in [-3.23185, -3.2245] + [-0.005848, 0.001505]i\\
    x_6 &\in [-0.003676, 0.003677] + [0.088467, 0.095821]i, & y_6 &\in [3.05057, 3.05792] + [0.011763, 0.019117]i\\
    x_7 &\in [-3.05792, -3.05057] + [0.011763, 0.019117]i, & y_7 &\in [-0.003676, 0.003677] + [0.088467, 0.095821]i,
  \end{align*}
}
which is also not contained in $\R^2$.
\end{proof}

\begin{proof}[Proof of Lemma~\ref{LMM:CMC}]
In Subsection~\ref{subsection6.3} we explained how to construct a family of projective boxes $\{\mathcal{B}^{\pm}_i\}_i$ which satisfies the (CMC) for a selected parameter $(a, b)=(5.7, 1.0)$. To extend this result to all over the parameters $(a, b) \in\mathcal{F}^{\pm}$, we proceed as follows. 

First, we choose $33$ (resp. $65$) ``sample parameters'' $(a, b)$ in $\mathcal{F}^+_{\R}$ (resp. $\mathcal{F}^-_{\R}$) of the form:
\[(a, b) = (a^{\pm}_{\mathrm{aprx}}(0.1\times k) + 0.1\times j, 0.1\times k)\in\mathcal{F}^{\pm}_{\R},\]
where $k$ and $j$ are integers. For each choice of sample parameter $(a, b) \in\mathcal{F}^+$ (resp. $(a, b) \in\mathcal{F}^-$) we carefully look at numerically drawn pictures of the trellis generated by $f_{a, b}$ and extract the coordinates of the $12$ (resp. $14$) intersection points $t^+_k$ (resp. $t^-_k$) appeared in Figure~\ref{FIG:trellis_positive_zero} (resp. Figure~\ref{FIG:trellis_negative_zero}). These points define quadrilaterals $\mathcal{Q}^{\pm}_i$ and their associated projective coordinates $\{(\pi^{\pm}_{u, i}, \pi^{\pm}_{v, i})\}_i$. 
Next, by trial and error, we find appropriate ``non-canonical constants'' and determine the topological disks $D^{\pm}_{u, i}$ and $D^{\pm}_{v, i}$ as in Subsection~\ref{subsection6.3} so that the projective boxes $\mathcal{B}^{\pm}_i=D^{\pm}_{u, i}\times_{\mathrm{pr}}D^{\pm}_{v, i}$ satisfies the (BCC) for each sample parameter $(a, b)\in \mathcal{F}_{\R}^{\pm}$ (see Definition~\ref{DFN:BCC}). We then linearly interpolate the coordinates of the intersection points $t^{\pm}_k$ (i.e. the data $\verb|tx[k]|$ and $\verb|ty[k]|$ in Table~\ref{TAB:a570b100}) and the data of the ``non-canonical constants'' for the topological disks $D^{\pm}_{u, i}$ and $D^{\pm}_{v, i}$ (i.e. the other data shown in Table~\ref{TAB:a570b100}) to all $(a, b)\in \mathcal{F}_{\R}^{\pm}$. For a complex parameter $(a, b) \in \mathcal{F}^{\pm} \setminus \mathcal{F}^{\pm}_{\R}$, the same boxes are used as the ones for its real part $(\Re(a), \Re(b))\in\mathcal{F}^{\pm}_{\R}$.
This defines a family of boxes $\{\mathcal{B}^{\pm}_i\}_i$ with respect to the family of projective coordinates $\{(\pi^{\pm}_{u, i}, \pi^{\pm}_{v, i})\}_i$ for all $(a, b)\in \mathcal{F}^{\pm}$.

Given a family of projective boxes, the verification of the boundary compatibility condition (BCC) is rather straightforward with the interval arithmetic. For example, to verify the (BCC) for the transition $(0, 2) \in \mathfrak{T}^{+}$, the absolute values of $\verb|delta_Px[0]|$ and $\verb|delta_Qx[0]|$ should be large enough so that the image of $\pi^+_{u, 2} \circ f(\partial^v \mathcal{B}^{+}_0)$ does not intersect with $D^+_{u, 2}$. However, if these values are too large, then it is likely that the (BCC) for the transition $(1, 0)$ fails, in turn. Therefore, we must choose adequate values of $\verb|delta_Px[i]|$ and $\verb|delta_Qx[i]|$ carefully so that the (BCC) holds for all possible transitions in $\mathfrak{T}^{+}$.
In practice, we divide $\mathcal{F}^+$ into 1,600,000 cubes (resp. $\mathcal{F}^-$ into 80,000,000 cubes), and for each such small cube we check the conditions in the lemma for all projective boxes corresponding to the cube. 

Another issue that we need to pay attention is the precision of the coordinate change. While the H\'enon map itself is defined in the Euclidean coordinate, the (BCC) is described in the projective coordinate. Therefore, the verification of the (BBC) involves the rigorous interval arithmetic for the coordinate change between them. This becomes problematic when the foci $u$ and $v$ are too close to the projective box (see Figure~\ref{FIG:projective_coordinates}), because then a small divisor appears in the coordinate change form the Euclidean one to the projective one, resulting loss of precision in the interval arithmetic. Therefore, again we must carefully adjust the values of $\verb|tx[i]|$ and $\verb|ty[i]|$ so that the foci are far enough from the boxes. 
\end{proof}

\begin{proof}[Proof of Lemma~\ref{LMM:bounding_K}]
  Roughly saying, we start with a cubical covering $\mathtt{C} \equiv \{C_i\}$ of $\mathcal{D}$
  and inductively remove cubes in $\mathtt{C}$ which does not intersect with
  $f(|\mathtt{C}|) \cap f^{-1}(|\mathtt{C}|)$ until $|\mathtt{C}|$ is
  contained in $\mathcal{B}^{\pm}$.
  Here we mean by $|\mathtt{C}|$ the union of all cubical sets in $\mathtt{C}$.
  However, since $f^{-1}$ is not defined when $b = 0$, we avoid using it by introducing
  appropriate ``flags'' for cubes.
  
  More precisely, we fix a cubical covering of $\mathcal{F}^{\pm}$ and for each
  cube, we check the statement of the Lemma as follows. Choose one of these
  parameter cubes and $\mathtt{F}$ be the cubical representations of $f_{a, b}$ on it.
  For each cube $C \in \mathtt{C}$ we assign two flags $C_{f}$ and $C_b \in \{\mathtt{true}, \mathtt{false}\}$
  that indicates the possibility for $C$ having intersection with
  $f(|\mathtt{C}|)$ and $f^{-1}(|\mathtt{C}|)$, respectively.
  We then run the following algorithm:

  \noindent\hrulefill
  \begin{algorithmic}
    \State{$\mathtt{C}$ := a cubical covering of $\mathcal{D}$}
    \State{$\mathtt{C}'$ = $\emptyset$}
    \While{$\mathtt{C}' \ne \mathtt{C}$}
    \State{$\mathtt{C}' = \mathtt{C}$}
    \State{Set $C_f = C_b = \mathtt{false}$ for all $C \in \mathtt{C}$}
    \Foreach{$c \in \mathtt{C}$}
    \If{$\mathtt{F}(C) \cap \mathtt{C} \ne \emptyset$}
    \State{Set $C_b = \mathtt{true}$}
    \State{Set $\widetilde{C}_f = \mathtt{true}$ for all $\widetilde{C} \in \mathtt{F}(C) \cap \mathtt{C}$}
    \EndIf
    \EndForeach
    \State{$\mathcal{C} = \{C \in \mathtt{C} \mid C_f = C_b = \mathtt{true}\}$}
    \If{$|\mathcal{C}| \subset \mathcal{B}^{\pm}$}
    \Return{true}
    \EndIf
    \EndWhile
    \State\Return{false}
  \end{algorithmic}
  \noindent\hrulefill

  If the algorithm returns true, then the statement of 
  Lemma~\ref{LMM:bounding_K} holds for all parameter values on the chosen parameter cube, 
  with $N$ being the number of ``while'' loops executed. 
  Otherwise, we subdivide each cubes in
  $\mathcal{C}$ and then run the algorithm again.
  We have checked that the algorithm returns true for all parameter cubes.
\end{proof}

\begin{proof}[Proof of Lemma~\ref{LMM:allowed_positive}]
Let $(a, b)\in\mathcal{F}^+$. Then, with the help of computer-assistance, we verify
\begin{enumerate}
\renewcommand{\labelenumi}{(\roman{enumi})}
\item $\mathcal{B}^+_i\cap\mathcal{B}^+_j\cap K_{a, b}=\emptyset$ for $(i, j)=(0, 1), (0, 2), (1, 3)$,
\item $f(\mathcal{B}^+_1\cap K_{a, b})\cap (\mathcal{B}^+_1\cup\mathcal{B}^+_2)=\emptyset$,
\item $f(\mathcal{B}^+_3\cap K_{a, b})\cap (\mathcal{B}^+_0\cup\mathcal{B}^+_3)=\emptyset$,
\item $f((\mathcal{B}^+_0\cap K_{a, b})\setminus \mathcal{B}^+_3)\cap \mathcal{B}^+_1=\emptyset$,
\item $f((\mathcal{B}^+_3\cap K_{a, b})\setminus (\mathcal{B}^+_0\cup\mathcal{B}^+_2))\cap \mathcal{B}^+_2=\emptyset$,
\item $f((\mathcal{B}^+_2\cap K_{a, b})\setminus (\mathcal{B}^+_1\cup\mathcal{B}^+_3))\cap (\mathcal{B}^+_0\cup\mathcal{B}^+_1)=\emptyset$,
\item $f((\mathcal{B}^+_1\cap K_{a, b})\setminus \mathcal{B}^+_2)\cap \mathcal{B}^+_3=\emptyset$.
\end{enumerate}

By (i), we see that $\mathcal{B}^+_I$ is empty for $I=\{0, 1\}, \{0, 2\}, \{1, 3\}$. By (ii), the arrows $\{1, 2\} \to \{2, 3\}$, $\{1, 2\} \to \{2\}$, $\{1, 2\} \to \{1, 2\}$, $\{1, 2\} \to \{1\}$, $\{1\} \to \{2, 3\}$, $\{1\} \to \{2\}$, $\{1\} \to \{1, 2\}$ and $\{1\} \to \{1\}$ are not allowed. By (iii), the transitions $\{0, 3\}\to \{0\}$, $\{0, 3\}\to \{0, 3\}$, $\{0, 3\} \to \{3\}$, $\{0, 3\} \to \{2, 3\}$, $\{3\} \to \{0\}$, $\{3\} \to \{0, 3\}$, $\{3\} \to \{3\}$, $\{3\} \to \{2, 3\}$, $\{2, 3\} \to \{0\}$, $\{2, 3\} \to \{0, 3\}$, $\{2, 3\} \to \{3\}$, $\{2, 3\} \to \{2, 3\}$ are not allowed. By (iv), the transitions $\{0\} \to \{1, 2\}$ and $\{0\} \to \{1\}$ are not allowed. By (v), the transitions $\{3\} \to \{2\}$ and $\{3\} \to \{1, 2\}$ are not allowed. By (vi), the transitions $\{2\} \to \{0\}$ and $\{2\} \to \{0, 3\}$ are not allowed. By (vii), the transitions $\{1\} \to \{0, 3\}$ and $\{1\} \to \{3\}$ are not allowed. 
\end{proof}

\begin{proof}[Proof of Lemma~\ref{LMM:allowed_negative}]
Let $(a, b)\in\mathcal{F}^-$. Then, with the help of computer-assistance, we verify
\begin{enumerate}
\renewcommand{\labelenumi}{(\roman{enumi})}
\item $\mathcal{B}^-_i\cap\mathcal{B}^-_j\cap K_{a, b}=\emptyset$ for $(i, j)=(0, 1), (0, 3), (0, 4), (1, 2), (1, 4), (2, 3)$,
\item $f((\mathcal{B}^-_0\cup \mathcal{B}^-_1)\cap K_{a, b})\cap (\mathcal{B}^-_1\cup \mathcal{B}^-_3)=\emptyset$,
\item $f((\mathcal{B}^-_2\cup \mathcal{B}^-_3)\cap K_{a, b})\cap (\mathcal{B}^-_0\cup \mathcal{B}^-_1)=\emptyset$,
\item $f(\mathcal{B}^-_4\cap K_{a, b})\cap (\mathcal{B}^-_0\cup \mathcal{B}^-_2)=\emptyset$,
\item $f((\mathcal{B}^-_0\cap K_{a, b})\setminus\mathcal{B}^-_2)\cap \mathcal{B}^-_4=\emptyset$,
\item $f((\mathcal{B}^-_1\cap K_{a, b})\setminus\mathcal{B}^-_3)\cap \mathcal{B}^-_4=\emptyset$,
\item $f((\mathcal{B}^-_2\cap K_{a, b})\setminus(\mathcal{B}^-_0\cup\mathcal{B}^-_4))\cap (\mathcal{B}^-_2\cup\mathcal{B}^-_3)=\emptyset$,
\item $f((\mathcal{B}^-_3\cap K_{a, b})\setminus(\mathcal{B}^-_1\cup\mathcal{B}^-_4))\cap (\mathcal{B}^-_2\cup\mathcal{B}^-_3)=\emptyset$,
\item $f((\mathcal{B}^-_4\cap K_{a, b})\setminus(\mathcal{B}^-_2\cup\mathcal{B}^-_3))\cap \mathcal{B}^-_4=\emptyset$.
\end{enumerate}

The rest of the proof is same as Lemma~\ref{LMM:allowed_positive}, hence omitted.
\end{proof}

\begin{proof}[Proof of Lemma~\ref{LMM:OCC_negative1}]
  We fix a cubical covering of $\mathcal{F}^{-}$ and for each
  cube we check the statement of Lemma~\ref{LMM:OCC_negative1} as follows. 
  Choose one of these parameter cubes and $F$ and $F^{2}$ be the cubical
  outer approximations of $f_{a, b}\circ\iota_v(u)$ and $f_{a, b}^{2}\circ\iota_v(u)$ on it.
  Denote the outer approximation of
  \[g(u, v)\equiv\frac{\partial}{\partial u}\left\{\pi^-_{u, 3} \circ f_{a, b}^2\circ\iota_v(u)\right\}\]
  by $G$.
  We then run the following algorithm:

  \noindent\hrulefill
  \begin{algorithmic}
    \State{$\mathtt{D}_v$ := a cubical coverings of $D_{v, 3}^{-}$}
    \State{$\mathtt{D}_u$ := a cubical coverings of $D_{u, 3}^{-}$}
    \State{$\mathtt{C} := \mathtt{D}_v \times \mathtt{D}_u$}
    \While{$\mathtt{C} \ne \emptyset$}
    \State{Subdivide cubes in $\mathtt{C}$}
    \Foreach{$C \in \mathtt{C}$}
    \If{($0 \not \in G(C)$) or ($|C| \cap \mathcal{B}^-_3 = \emptyset$) or ($F(C) \cap \mathcal{B}^-_4 = \emptyset$) or
    ($F^2(C) \cap \mathcal{B}^-_3 = \emptyset$)}
    \State{remove $C$ from $\mathtt{C}$}
    \EndIf
    \EndForeach
    \EndWhile
    \State\Return{true}
  \end{algorithmic}
  \noindent\hrulefill

  If the algorithm returns true, then it implies that 
  $\mathcal{C} = \emptyset$ holds at some subdivision level
  and therefore the statement of Lemma~\ref{LMM:OCC_negative1} holds 
  for all parameter values on the chosen parameter cube.
  Otherwise, the algorithm does not terminates.
  We have checked that the algorithm returns true for all parameter cubes.
\end{proof}

\begin{proof}[Proof of Lemma~\ref{LMM:OCC_negative2}]
The proof of this claim is similar to the previous one, hence omitted. 
\end{proof}

\begin{proof}[Proof of Lemma~\ref{LMM:tangency_positive}]
  We fix a cubical covering of $\partial^v \mathcal{F}^{+}$ and for each
  cube we check the statement of Lemma~\ref{LMM:tangency_positive} as follows.
  Notations $F, F^2$ and $G$ are the same as in the proof of Lemma~\ref{LMM:OCC_negative1}.
  We then run the following algorithm:

  \noindent\hrulefill
  \begin{algorithmic}
    \State{$\mathtt{D}_v$ := a cubical coverings of $\widehat{D}_{v, 0}^+$}
    \State{$\mathtt{D}_u$ := a cubical coverings of $D_{u, 0}^{+}$}
    \State{$\mathtt{V}$ := a cubical coverings of $\mathcal{V}^s_{31\overline{0}}(a, b)^+$}
    \State{$\mathtt{C} := \mathtt{D}_{v} \times \mathtt{D}_{u}$}
    \While{$\mathtt{C} \ne \emptyset$}
    \State{Subdivide cubes in $\mathtt{C}$}
    \State{Refine the covering $\mathtt{V}$}
    \Foreach{$C \in \mathtt{C}$}
    \If{\scalebox{0.925}{($0 \not \in G(C)$) or ($|C| \cap \mathcal{B}^+_0 = \emptyset$) or ($F(C) \cap \mathcal{B}^+_2 = \emptyset$) or
    ($F^2(C) \cap (|V|\times_{\mathrm{pr}}D_{v, 3}^+) = \emptyset$)}}
    \State{remove $C$ from $\mathtt{C}$}
    \EndIf
    \EndForeach
    \EndWhile
    \State\Return{true}
  \end{algorithmic}
  \noindent\hrulefill

  Note that in the step of refining $\mathtt{V}$ in the ``while'' loop, 
  we need to construct a tight outer approximation of the set
  \[\mathcal{V}^s_{\mathrm{loc}}(p_1) = \mathcal{B}^+_0\cap f^{-1}(\mathcal{B}^+_0)\cap \cdots \cap f^{-N+1}(\mathcal{B}^+_0)\cap f^{-N}(\mathcal{B}^+_0).\]
  However, this can be done with the same algorithm as in the proof of Lemma~\ref{LMM:bounding_K}
  (In fact, if we ignore $C_f$ and define $\mathcal{C} = \{C \in \mathtt{C} \mid C_b = \mathtt{true}\}$ in the algorithm, then it computes a covering of the local stable manifold).

  If the algorithm returns true, then the statement of 
  Lemma~\ref{LMM:tangency_positive} holds for all parameter values 
  on the chosen parameter cube.
  Otherwise, the algorithm does not terminates.
  We have checked that the algorithm returns true for all parameter cubes.
\end{proof}

\begin{proof}[Proof of Lemma~\ref{LMM:tangency_negative}]
The proof of this lemma is similar to the previous one, hence omitted. 
\end{proof}

\newpage

\begin{appendix}

\section{Regularity of Loci Boundary}\label{appendix:regularity}

In this appendix we collect some basic definitions and facts on complex subvarieties (analytic subsets) which are essential in the proof of the Main Theorem. Moreover, we take this opportunity to quote a proof of Lemma 1.1 in~\cite{BS0}, which is in fact missing in its published version~\cite{BS2}. We refer to~\cite{C} for the generalities on complex subvarieties. 

\begin{figure}
  \includegraphics[height=6cm]{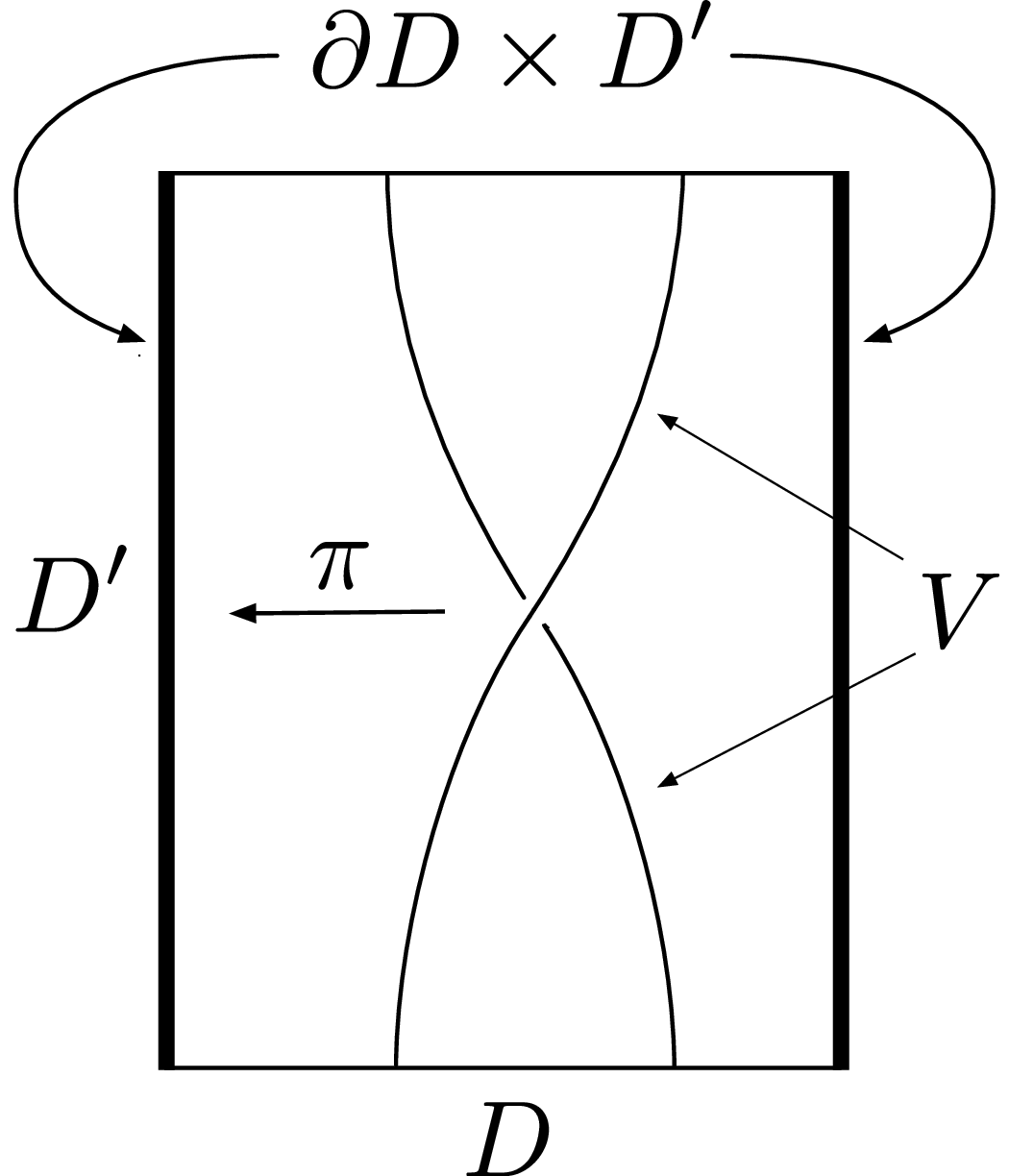}
  \caption{Properness of the projection.}
\label{FIG:proper}
\end{figure}

Below $X$ and $Y$ are assumed to be Hausdorff and locally compact topological spaces. We start with a simple criterion for a projection to be proper, which is used in the proof of Proposition~\ref{PRP:complex_loci} in Subsection~\ref{subsection5.2}. For a proof, see (3) on page 29 of~\cite{C}. 

\begin{lmm}
Let $D\subset X$ and $D'\subset Y$ be subsets with $\overline{D}$ compact and let $V$ be a closed subset in $D\times D'$. Let $\pi : D\times D'\to D'$ be the projection. Then, the restriction of the projection $\pi : V\to D'$ is proper iff $\overline{V}\cap (\partial D\times D')=\emptyset$, where the closure of $V$ is taken in $X\times Y$ (see Figure~\ref{FIG:proper}).
\label{LMM:proper}
\end{lmm}

Let $\mathrm{\Omega}\subset\C^n$ be a domain. Recall the following notion.

\begin{dfn}
A subset $V\subset \mathrm{\Omega}$ is called a \textit{complex subvariety} (or an \textit{analytic subset}) of $\mathrm{\Omega}$ if for each point $z\in V$ there exist a neighborhood $U$ of $z$ and finitely many holomorphic functions $f_i$ ($i=1, \dots, N$) on $U$ so that $V\cap U$ is the set of common zeros of $f_i$.
\label{DFN:subvariety}
\end{dfn}

The next fact is also crucial in the proof of Proposition~\ref{PRP:complex_loci} in Subsection~\ref{subsection5.2}.

\begin{prp}
Let $D\subset \C^n$ and $D'\subset \C^m$ be open subsets and let $\pi : D\times D'\to D'$ be the projection. Assume that $V\subset D\times D'$ is an analytic subset and $\pi : V\to D'$ is proper of degree one. Then, $V$ is a complex submanifold in $D\times D'$ and $\pi : V\to D'$ is biholomorphic.
\label{PRP:submanifold}
\end{prp}

This follows from the well-known Weierstrass' preparation theorem. See Proposition 3 on p.32 of~\cite{C} for a proof. 

Now we prove that the complex tangency loci $\mathcal{T}^{\pm}$ form complex subvarieties. Consider a holomorphic family of biholomorphic maps $f_{\lambda} : \C^2\to\C^2$ defined for $\lambda\in \mathrm{\Lambda}\subset\C^N$. Fix $\lambda_0\in \mathrm{\Lambda}$ and assume that $f_{\lambda_0}$ has two saddle points $p^s_{\lambda_0}, p^u_{\lambda_0}\in\C^2$. Let $p^s_{\lambda}, p^u_{\lambda}$ be their continuations and let $V^s(p^s_{\lambda}; f_{\lambda})$ and $V^u(p^u_{\lambda}; f_{\lambda})$ be their stable and unstable manifolds for $f_{\lambda}$ respectively. Assume that $V^s(p^s_{\lambda_0}; f_{\lambda_0})$ and $V^u(p^u_{\lambda_0}; f_{\lambda_0})$ intersect tangentially and let $z_0$ be a such intersection point. Let $\psi^{s/u}(\, \cdot\, , \lambda) : \C\to \C^2$ be the uniformizations of $V^{s/u}(p^{s/u}_{\lambda}; f_{\lambda})$ such that $\psi^{s/u}(0, \lambda_0)=z_0$. Since $z_0$ is an isolated point of $V^s(p^s_{\lambda_0}; f_{\lambda_0})\cap V^u(p^u_{\lambda_0}; f_{\lambda_0})$ with respect to their leaf topology, there exists $\varepsilon>0$ so that
\begin{equation}
\inf_{(\zeta^s, \zeta^u)\in X} \mathrm{dist}(\psi^s(\zeta^s, \lambda), \psi^u(\zeta^u, \lambda))\geq \delta>0,
\end{equation}
holds for $\lambda=\lambda_0$, where 
\[X\equiv \bigl\{(\zeta^s, \zeta^u)\in\C^2 : |\zeta^s|\leq \varepsilon,\ |\zeta^u|=\varepsilon\bigr\}\cup\bigl\{(\zeta^s, \zeta^u)\in\C^2 : |\zeta^s|= \varepsilon,\ |\zeta^u|\leq\varepsilon\bigr\}.\]
Since $X$ is compact, there exists a neighborhood $U$ of $\lambda_0$ so that (A.1) holds for all $\lambda\in U$. 

By writing as $\psi^{s/u}=(\psi^{s/u}_1, \psi^{s/u}_2)$, the two tangent vectors $\partial_{\zeta}\psi^s(\zeta^s, \lambda)$ and $\partial_{\zeta}\psi^u(\zeta^u, \lambda)$ are parallel iff
\begin{equation}
\partial_{\zeta}\psi_1^s(\zeta^s, \lambda)\cdot\partial_{\zeta}\psi_2^u(\zeta^u, \lambda)=\partial_{\zeta}\psi_2^s(\zeta^s, \lambda)\cdot\partial_{\zeta}\psi_1^u(\zeta^u, \lambda)
\end{equation}
holds. Then,
\[M\equiv\bigl\{(\zeta^s, \zeta^u, \lambda)\in\C^2\times U : |\zeta^s|, |\zeta^u|<\varepsilon,\ \psi(\zeta^s, \lambda)=\psi(\zeta^u, \lambda) \mbox{ and (A.2) hold}\bigr\}\]
forms a complex subvariety of $\{\zeta^s\in\C : |\zeta^s|<\varepsilon\}\times \{\zeta^u\in\C : |\zeta^u|<\varepsilon\}\times U$. Let $\pi : (\zeta^s, \zeta^u, \lambda)\mapsto \lambda$ be the projection to $U$ and set 
\[\mathcal{T}(z_0, \lambda_0)\equiv \pi(M).\]
Thus, $\mathcal{T}(z_0, \lambda_0)$ is the locus of parameters $\lambda$ near $\lambda_0$ for which $V^s(p^s_{\lambda}; f_{\lambda})$ has a tangential intersection with $V^u(p^u_{\lambda}; f_{\lambda})$ near $z_0$ in the leaf topology. Now we are ready to state Lemma 1.1 of~\cite{BS0} as

\begin{prp}
The locus $\mathcal{T}(z_0, \lambda_0)$ is a complex subvariety of $U$.
\label{PRP:T_is_subvariety}
\end{prp}

\begin{proof}
Thanks to Lemma~\ref{LMM:proper}, the projection $\pi : M \to U$ is proper. Since a proper projection of a subvariety is again a subvariety by Theorem in page 29 of~\cite{C}, we know that $\mathcal{T}(z_0, \lambda_0)=\pi(M)$ is a subvariety of $U$. 
\end{proof}

\newpage

\section{Comparison of Box Systems}\label{appendix:comparison}

Recall that in Subsection~\ref{subsection2.2} we have employed a $4$-box system $\{\mathcal{B}^+_i\}_{i=0}^3$ for $(a, b)\in\mathcal{F}^+$ based on a trellis formed by the invariant manifolds of the saddle fixed point $p_1$ and the saddle periodic points $p_2$ and $p_4$ of period two. It is of course possible to construct a $5$-box system $\{\mathcal{B}'_i\}_{i=0}^4$ for $(a, b)\in\mathcal{F}^+$ in a similar manner to the case $(a, b)\in\mathcal{F}^-$ based on a trellis formed by the invariant manifolds of the two saddle fixed points $p_1$ and $p_3$. However, when $b$ is close to $1$, the fixed point $p_3$ is relatively close to the $y$-axis and thus the expansion and the contraction at this point are relatively weak compared to the case $(a, b)\in\mathcal{F}^-$. In fact, the multipliers are $\lambda^u(p_3)\approx -2.8$ and $\lambda^s(p_3)\approx -0.35$ for $(a, b)=(5.7, 1)$, but $\lambda^u(p_3)\approx -5.2$ and $\lambda^s(p_3)\approx 0.19$ for $(a, b)=(6.2, -1)$. Presumably due to this fact, our numerical experiments suggest that it seems impossible for the ``neighboring transitions'' $f : \mathcal{B}'_3\cap f^{-1}(\mathcal{B}'_4)\to \mathcal{B}'_4$ and $f : \mathcal{B}'_4\cap f^{-1}(\mathcal{B}'_3)\to \mathcal{B}'_3$ around $p_3$ to verify the (BCC) when $b$ is close to $1$. On the other hand, our $4$-box system $\{\mathcal{B}^+_i\}_{i=0}^3$ avoids such neighboring transitions and we were able to verify the (BCC) with this system. This is the main advantage of choosing the $4$-box system $\{\mathcal{B}^+_i\}_{i=0}^3$. 

\begin{figure}
\begin{center}
\includegraphics[width=15.5cm]{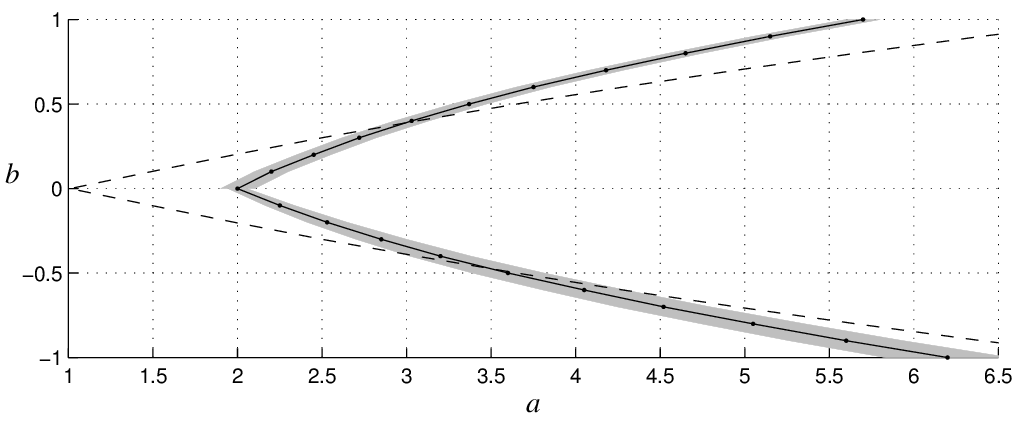}
\end{center}
\caption{Comparison of $\mathcal{F}^{\pm}_{\R}$ (shaded), the $3$-box system~\cite{BS2} (dashed) and the graphs of $a^{\pm}_{\mathrm{aprx}}$ (solid).}
\label{FIG:comparison}
\end{figure}

Next we discuss the $3$-box system introduced in~\cite{BS2}. For $(a, b)\in \R\times\R^{\times}$, let 
\[R \equiv \frac{1+|b|+\sqrt{(1+|b|)^2+4a}}{2}\]
and
\begin{align*}
D_0 & \equiv \bigl\{x\in\C : 0<|x|<R, -\pi/2 < \arg x < \pi/2\bigr\}, \\
D_1 & \equiv \bigl\{x\in\C : |x|<R\bigr\}\cap p_c^{-1}\bigl(\bigl\{x\in\C: \mathrm{Re}\, (x) < |b|R\bigr\}\bigr), \\
D_2 & \equiv \bigl\{x\in\C : 0<|x|<R, \pi/2 < \arg x < 3\pi/2\bigr\}.
\end{align*}
Then, the $3$-box system in~\cite{BS2} is defined as $\mathcal{B}_i \equiv D_i\times \{y\in\C : |y|<R\}\subset \C^2$ for $i=0, 1, 2$. Put $\alpha\equiv\sqrt{|b|R+a}$ so that $[-\alpha, \alpha]=\R\cap D_1$. Then, a sufficient condition for the (BCC) is $a > \sqrt{|b|R+a} + |b|R$ (this condition looks close to optimal). 

The shaded region\footnote{At $b=1$ we verified the (BCC) for $5.60\leq a \leq 5.80$ as in (iii) of Theorem~\ref{THM:quasitrichotomy} (Quasi-Trichotomy). However, we were not able to verify the (BCC) for $a<5.60$ even with our $4$-box system for several choices of non-canonical constants explained in Subsection~\ref{subsection6.3}.}, the dashed and solid lines in Figure~\ref{FIG:comparison} are the regions $\mathcal{F}^{\pm}_{\R}$, the curve $a = \sqrt{|b|R+a}+|b|R$ and the graph of the function $a^{\pm}_{\mathrm{aprx}}$, respectively. The figure illustrates that the $3$-box system of~\cite{BS2} works only for $-0.5 < b < 0.4$ near $\partial \mathcal{H}^{\pm}_{\R}=\partial \mathcal{M}^{\pm}_{\R}$.

\end{appendix}

\end{document}